\input amstex
\input xy
\xyoption{all}
\documentstyle{amsppt}
\document

\magnification 1000

\def\gen{\frak{g}}

\def\ten{\frak{t}}
\def\aen{\frak{a}}
\def\ben{\frak{b}}

\def\Fen{\frak{F}}

\def\hen{\frak{h}}

\def\len{\frak{l}}

\def\qen{\frak{q}}
\def\cen{\frak{c}}

\def\sen{\frak{s}}
\def\uen{\frak{u}}

\def\zen{\frak{z}}

\def\Aen{\frak{A}}

\def\Cen{\frak{C}}

\def\Een{\frak{E}}
\def\Fen{\frak{F}}

\def\Sen{\frak{S}}

\def\Xen{\goth{X}}

\def\a{{\alpha}}
\def\g{{\gamma}}
\def\o{{\omega}}
\def\l{{\lambda}}
\def\b{{\beta}}
\def\eps{{\varepsilon}}

\def\lub{{{}^B}}

\def\un{{{\bold 1}}}

\def\bb{{\bold b}}
\def\lb{{\bold b}}

\def\gb{{\bold g}}

\def\lb{{\bold l}}

\def\nb{{\bold n}}

\def\qb{{\bold q}}

\def\nub{{\nu,\kappa}}

\def\tb{{\bold t}}
\def\ub{{\bold u}}

\def\zb{{\bold z}}
\def\Ab{{\bold A}}

\def\Bb{{\bold B}}

\def\wu{{\underline w}}

\def\Hb{{\bold H}}

\def\Rb{{\bold R}}

\def\Sb{{\bold S}}

\def\Lb{{\bold L}}

\def\Vb{{\bold V}}

\def\f{{fg}}

\def\c{{\roman c}}

\def\s{{\roman s}}

\def\W{{\roman W}}

\def\rhoaf{{\rhob}}

\def\Aut{\text{Aut}}
\def\Ext{\text{Ext}}

\def\Hom{{\roman{Hom}}}

\def\Res{\text{Res}}

\def\dim{{\roman{dim}}}

\def\Ind{{\Gamma}}

\def\End{{\roman{End}}}
\def\eu{{\roman{eu}}}
\def\ch{{\roman{ch}}}
\def\Id{{\roman{id}}}
\def\Irr{{\roman{Irr}}}
\def\Oplus{\ts\bigoplus}

\def\Prod{\ts\prod}

\def\tr{\roman{tr}}

\def\Ker{\roman{Ker}\,}

\def\GL{{\roman{GL}}}

\def\re{{\roman{re}}}

\def\re{{{\text{re}}}}

\def\Aloc{{{C}}}
\def\Alocn{{C_n}}

\def\Rnl{{\Rb_\nl}}

\def\Alocnl{{C_{n,\ell}}}
\def\nl{{{n,\ell}}}

\def\Alocul{{C_{1,\ell}}}

\def\reg{{{\heartsuit}}}
\def\KZ{{{\roman{KZ}}}}

\def\proj{{{\roman{proj}}}}
\def\projb{{{\bold{proj}}}}

\def\op{{{\roman{op}}}}

\def\CI{{{\goth C}}}

\def\CIloc{{{\CI}}}
\def\V{{{\Een}}}
\def\W{{{\Fen}}}

\def\VVkap{{{\V}}}
\def\VVKZ{{{\V_{\KZ}}}}

\def\VVCC{{{\V_{\CC}}}}

\def\dom{{\sss\geqslant 0}}
\def\ddom{{\sss> 0}}
\def\alc{{{\Aen}}}

\def\domm{{\sss >0}}
\def\Rloc{{\Rb}}
\def\Rlocn{{\Rb_n}}
\def\Rlocnl{{\Rb_{n,\ell}}}
\def\Bnl{{\Bb_{n,\ell}}}

\def\vk{{{}^{\sss\varkappa}\!}}
\def\dg{{{}^{\sss\dag}\!}}
\def\ddg{{{}^{\sss\ddag}\!}}
\def\sp{{{}^{\sss\sharp}\!}}
\def\s{{{}_S}}

\def\CC{{\Bbb C}}

\def\NN{{\Bbb N}}

\def\PP{{\Bbb P}}
\def\QQ{{\Bbb Q}}
\def\RR{{\Bbb R}}

\def\ZZ{{\Bbb Z}}

\def\Ac{{\Cal A}}
\def\Bc{{\Cal B}}
\def\Cc{{\Cal C}}

\def\Gc{{\Cal G}}
\def\Hc{{\Cal H}}

\def\Mc{{\Cal M}}

\def\Oc{{\Cal O}}

\def\Pib{{{\hat\Pi}}}
\def\rhob{{\hat\rho}}
\def\Pc{{\Cal P}}

\def\Sc{{\Cal S}}

\def\mod{{\text{mod}}}
\def\modb{{\bold{mod}}}

\def\tob{{\buildrel b\over\to}}

\def\and{{\text{and}}}

\def\ds{\displaystyle}
\def\ts{\textstyle}

\def\sss{\scriptscriptstyle}
\def\qed{\hfill $\sqcap \hskip-6.5pt \sqcup$}        
\overfullrule=0pt                                    

\def\7dag{{{\!\!\!\!\!\!\!\dag}}}
\def\6dag{{{\!\!\!\!\!\!\dag}}}
\def\5dag{{{\!\!\!\!\!\dag}}}
\def\4dag{{{\!\!\!\!\dag}}}
\def\3dag{{{\!\!\!\dag}}}
\def\2dag{{{\!\!\dag}}}
\def\1dag{{{\!\dag}}}

\def\pro{{\lim\limits_{\lla}}}

\def\la{{\langle}}
\def\ra{{\rangle}}
\def\lla{{\longleftarrow}}
\def\lra{{{\longrightarrow}}}

\newdimen\Squaresize\Squaresize=14pt
\newdimen\Thickness\Thickness=0.5pt
\def\Square#1{\hbox{\vrule width\Thickness
          \alphaox to \Squaresize{\hrule height \Thickness\vss
          \hbox to \Squaresize{\hss#1\hss}
          \vss\hrule height\Thickness}
          \unskip\vrule width \Thickness}
          \kern-\Thickness}
\def\Vsquare#1{\alphaox{\Square{$#1$}}\kern-\Thickness}

\nologo

\topmatter
\title
Cyclotomic double affine Hecke algebras and affine parabolic
category $\Oc$, I
\endtitle
\rightheadtext{} \abstract Using the orbifold KZ connection we
construct a functor from an affine parabolic category $\Oc$ of type
$A$ to the category $\Oc$ of a cyclotomic rational double affine
Hecke algebra $\Hb$. Then we give several results concerning this
functor.
\endabstract
\author M. Varagnolo, E. Vasserot\endauthor
\address D\'epartement de Math\'ematiques,
Universit\'e de Cergy-Pontoise, 2 av. A. Chauvin, BP 222, 95302
Cergy-Pontoise Cedex, France, Fax : 01 34 25 66 45\endaddress \email
michela.varagnolo\@math.u-cergy.fr\endemail
\address D\'epartement de Math\'ematiques,
Universit\'e Paris 7, 175 rue du Chevaleret, 75013 Paris, France,
Fax : 01 44 27 78 18
\endaddress
\email vasserot\@math.jussieu.fr\endemail
\thanks
2000{\it Mathematics Subject Classification.} Primary ??; Secondary
??.
\endthanks
\endtopmatter
\document

\head Introduction \endhead

Fix positive integers $n$, $\ell$.  Let $W$ be the wreath product of
$\Sen_n$ and $\ZZ/\ell\ZZ$. The cyclotomic rational double affine
Hecke algebra $\Hb$ is a deformation of the semi-direct product
$\CC[\CC^{2n}]\rtimes W$. Its category $\Oc$ 
is a quasi-hereditary cover of the
Ariki-Koike algebra, see \cite{R}. One important problem is to
compute the dimension of the simple modules of the category $\Oc$ of $\Hb$, 
or,
equivalently, the Jordan-H\"older multiplicities of the standard
modules. So far the main approach to the representation theory of
double affine Hecke algebras associated with complex reflection
groups is geometric and uses D-modules on quiver varieties. However
a dimension formula for simple modules seems out of reach yet by
these techniques. It is expected that the Jordan-H\"older
multiplicities of the standard modules are values at one of some
affine parabolic Kazhdan-Lusztig polynomial, see \cite{R, sec. 6.5}.
We'll call this the dimension conjecture. See Section 8 below for
details. These multiplicites are encoded in a combinatorial object
called the level $\ell$ Fock space.

If $\ell=1$ the dimension conjecture is proved. It follows from
\cite{R} and \cite{VV}, or from \cite{Su2} and \cite{V}. In this case
there is another algebraic approach to the algebra $\Hb$ due to
Suzuki. He constructed a functor from Kazhdan-Lusztig's category of
modules over the type $A^{(1)}$ affine Lie algebra to the category
$\Oc$ of $\Hb$. We give a proof that this functor is an equivalence
in Section A.5 below. This functor is constructed via affine
coinvariants over the configuration space of $\PP^1$ and the
Knizhnik-Zamolodchikov connection.

In this paper we construct a similar functor for any $\ell$. The new
ingredient is the space of orbifold affine coinvariants over the
configuration space of the stack $[\PP^1/(\ZZ/\ell\ZZ)]$ and the
corresponding Knizhnik-Zamolodchikov connection. A priori this space
of coinvariants involves a choice of a twisted affine Lie algebra.
Choosing an inner twist of the type $A^{(1)}$  affine Lie algebra,
we get a functor $\Een$ from an affine parabolic category $\Oc$ to
the category $\Oc$ of $\Hb$. Then we study $\Een$, in particular its
behavior on standardly filtered modules. We do not prove the
dimension conjecture. Contrarily to the case $\ell=1$ mentioned
above, the functor $\Een$ is not an equivalence of quasi-hereditary
categories in general. However, we expect the functor $\Een$ to be
an important tool to prove it. In particular $\Een$ should behave
nicely on indecomposable projective modules, as explained in this
paper. We'll come back to this elsewhere.

Let us now describe the structure of the paper. The first and second
sections are reminders on DAHA's (=double affine Hecke algebras) and
Suzuki's functor. In the third one we consider the case $\ell\neq
1$. Using the orbifold Knizhnik-Zamolodchikov connection we define a
functor taking a smooth module over a twisted affine Lie algebra to
a $\Hb$-module and we compute the image of parabolic Verma modules.
In the fourth section we compare the space of twisted affine
coinvariants and the space of non-twisted ones when the affine Lie
algebra is equipped with an inner twist. In the sixth section we
define the functor $\Een$. It goes from the affine parabolic
category $\Oc$ of type $A^{(1)}$ to the category $\Oc$ of $\Hb$. We
prove that $\Een$ preserves the posets of standard modules and is
exact on standardly filtered modules in section seven. We conjecture
that it preserves the set of indecomposable projective modules. In
the last section we compare $\Een$ with what one expects from the
dimension conjecture.

\vskip2cm

\head Contents \endhead

\item{0.} Notations
\item{1.} The cyclotomic rational DAHA and the Dunkl operators
\item{2.} The affine category $\Oc$
\item{3.} Twisted affine coinvariants
\item{4.} Untwisting the space of twisted affine coinvariants
\item{5.} Complements on the category $\hat\Oc_\nub$
\item{6.} Definition of the functor $\Een$
\item{7.} The functor $\Een$ is exact on standardly filtered modules
\item{8.} The affine parabolic category $\Oc$ and the Fock space
\item{A.} Appendix
\itemitem{A.1.} Proof of Proposition 5.8
\itemitem{A.2.} Proof of Proposition 7.2 and Corollary 7.3
\itemitem{A.3.} Reminder on induction
\itemitem{A.4.} Reminder on the Fock space
\itemitem{A.5.} The functor $\Een$ for $\ell=1$
\itemitem{A.6.} Proof of Proposition 8.7

\vskip2cm

\head 0. Notations\endhead

\subhead 0.1\endsubhead
First, let us gather a few basic notations on categories.
The categories we'll consider are all $\CC$-linear, i.e., they are additive and the $\Hom$
sets are $\CC$-vector spaces. A category is  {\it Artinian} if the Hom sets
are finite dimensional $\CC$-vector spaces and every object has a
finite length.
Write $\Ac^\f$ for the full subcategory
of an Abelian category $\Ac$ consisting of the objects of finite
length and $\Ac^\proj\subset\Ac^\f$ for the full exact subcategory
of projective objects.
Given a set $I$ of objects of the Abelian category $\Ac$, 
we denote by $\Ac^I$ the exact full subcategory of
$I$-filtered objects, i.e., of objects $M$ with a finite 
filtration such that
each successive quotient is isomorphic to an object of $I$.

Write $[M:S]$ for the Jordan-H\"older multiplicity of a simple
module $S$ in an object $M$ of finite length.

By a {\it quasi-hereditary category} we mean a highest weight category in
the sense of \cite{CPS} which is equivalent, as a highest weight
category, to the category of finitely generated modules of a
quasi-hereditary (finite dimensional) $\CC$-algebra. 
In other words, it is an Artinian Abelian category with a finite poset $\Delta_\Ac$
of standard modules satisfying
the following axioms :

\vskip1mm
\item{$(a)$} 
we have $\End_\Ac(M)=\CC$ for each $M\in\Delta_\Ac$,
\vskip1mm
\item{$(b)$} 
if $\Hom_\Ac(M_1,M_2)\neq 0$ and $M_1,M_2\in\Delta_\Ac$ then $M_1\leqslant M_2$,

\vskip1mm
\item{$(c)$} 
if $\Hom_\Ac(M,N)= 0$ for each $M\in\Delta_\Ac$ then $N=0$,

\vskip1mm
\item{$(d)$} 
if $M\in\Delta_\Ac$ there is a projective object $P\in\Ac$ and an epimorphism
$P\to M$ whose kernel is $\Delta_\Ac$-filtered with subquotients $>M$.

\vskip1mm

\noindent
See \cite{CPS, thm.~3.6}, \cite{Do, appendix A} for more details.
An equivalence of highest weight categories is an equivalence of
categories which restricts to a bijection between both sets of
standard modules. We'll abbreviate
$\Ac^\Delta=\Ac^{\Delta_\Ac}$.

We'll write $[\Bc]$ for the Grothendieck group of an Abelian or an
exact category $\Bc$. Let $[M]$ denote the class in $[\Bc]$ of an
object $M$.
Note that the obvious
embedding $\Ac^\Delta\subset\Ac^\f$ yields a group isomorphism
$[\Ac^\Delta]=[\Ac^\f]$.

\subhead 0.2\endsubhead Let $R$ be a commutative Noetherian ring
with 1, and let $\Ab$ be an $R$-algebra. Write $\Ab$-$\modb$ for the
module category of $\Ab$, $\Ab$-$\projb$ for $\Ab$-$\modb{}^\proj$
and $\Ab$-$\modb^\f$ for
$(\Ab$-$\modb)^\f$. Let $\Irr(\Ab)$ be the set of isomorphism classes
of simple objects of $\Ab$-$\modb^\f$. To any $R$-algebra
homomorphism $\phi:\Ab\to\Bb$ we associate the functor
$$\phi:\Bb\text{-}\modb\to\Ab\text{-}\modb,\quad
M\mapsto{}^\phi M,$$ where ${}^\phi M$ is the twist of $M$ by
$\phi$.

\subhead 0.3\endsubhead Let $M$ be a $\CC$-vector space and $R$ be a
commutative $\CC$-algebra. Assume that $R$ is the functions algebra
of a $\CC$-variety $R$. We write 
$$M[X]=M_R=M\otimes R.$$
Given an automorphism $F$ of a set $M$, let $M^F$ be
the fixed points set. If $F$ is a $R$-linear automorphism $F$ of
$M_R$ we may abbreviate $M_R^F=(M_R)^F$.

\vskip2cm

\head 1. The cyclotomic rational DAHA and the Dunkl
operators\endhead

\subhead 1.1.~The complex reflection group
$G(\ell,1,n)$\endsubhead Let $D_\ell\subset\CC^\times$ be the group
consisting of the $\ell$-th roots of unity. Fix a generator of $D_\ell$
once for all. We'll denote it by $\eps$. Let $\Sen_n$ be the
symmetric group on $n$ letters and $W$ be the semi-direct product
$\Sen_n\ltimes (D_\ell)^n$, where $(D_\ell)^n$ is the Cartesian product of $n$
copies of $D_\ell$. To avoid any confusion we may write
$$W_A=W,\quad A=\{1,2,\dots n\}.$$ Let $\eps_i\in (D_\ell)^n$ be the element with $\eps$ at the
$i$-th place and to 1 at the other ones. 
If $a\in\ZZ$ we may identify $\ZZ/a\ZZ$ with the set $\{1,2,\dots a\}$ in the
obvious way, hoping it will not create any confusion.
Set
$$\Lambda=\{1,2,\dots,\ell\}\simeq\ZZ/\ell\ZZ.$$
For each $p\in\Lambda$
and each $i\neq j$, we write $s_{i,j}^{(p)}$ for
$s_{i,j}\eps_i^p\eps_j^{-p}$.

\subhead 1.2.~The cyclotomic rational DAHA\endsubhead
Fix a basis $(x,y)$ of $\CC^2$. Let $x_i$,
$y_i$ denote the elements $x,y$ respectively in the $i$-th summand
of $(\CC^2)^{\oplus n}$. The group $W$ acts on $(\CC^2)^{\oplus n}$
such that for distinct $i$, $j$, $k$ we have
$$\eps_i(x_i)=\eps^{-1} x_i,
\quad \eps_i(x_j)=x_j, \quad \eps_i(y_i)=\eps y_i, \quad
\eps_i(y_j)=y_j,$$
$$s_{i,j}(x_i)=x_j,
\quad s_{i,j}(y_i)=y_j, \quad s_{i,j}(x_k)=x_k, \quad
s_{i,j}(y_k)=y_k.$$ Fix  $k\in\CC$ and $\g_p\in\CC$ for $0\neq
p\in\Lambda$. We'll write $\g$ for the $\ell$-tuple $(\g_p)$.
The CRDAHA (=cyclotomic rational DAHA) is the quotient $\Hb_{k,\g}$
of the smash product of $\CC W$ and the tensor algebra of
$(\CC^2)^{\oplus n}$ by the relations
$$[y_i,x_i]=1-k\sum_{j\neq i}\sum_{p}s_{i,j}^{(p)}
-\sum_{p\neq 0}\g_p\eps_i^p,$$
$$[y_i,x_j]=k\sum_{p}\eps^{p}s_{i,j}^{(p)}
\quad \roman{if}\ i\neq j,$$
$$[x_i,x_j]=[y_i,y_j]=0.$$
This presentation is
the same as in \cite{EG}. We'll use another presentation where the
parameters are $h$, $h_p$ with $p\in\Lambda$ and $\sum_ph_p=0$.
Set $H=(h_1,h_2,\dots h_{\ell-1})$. The corresponding algebra is
denoted by the symbol $\Hb_{h,H}$. It is isomorphic to the algebra $\Hb_{k,\g}$
with $k=-h$ and $\g_p=-\sum_{p'\in\Lambda}\eps^{-pp'}h_{p'}$.
Our parameter $h_p$ is the same as the parameter $H_p$ in \cite{G}
and it is equal to $h_{H_0,p}-h_{H_0,p-1}$ with respect to the
parameters $h_{H_0,p}$, $p\in\Lambda$, in \cite{R}.

\subhead 1.3.~The Dunkl operators\endsubhead The subalgebras of
$\Hb_{h,H}$ generated by $\{x_1,\dots x_n\}$ and $\{y_1,\dots y_n\}$
are free commutative. We'll write $\Rb$ for the first one and
$\Rb^*$ for the second one. Note that $\Rb=\CC[\CC^n]$. There is a
$\CC$-linear representation of $\Hb_{h,H}$ on $\Rb$ such that
$x_i\mapsto x_i,$ $w\mapsto w$ and $y_i\mapsto \bar y_i$ with
$$\bar y_i=\partial_{x_i}+
k\sum_{j\neq i}\sum_{p}{1\over x_i-\eps^{-p}x_j} (s_{i,j}^{(p)}-1)+
\sum_{p\neq 0}{\g_p\over x_i-\eps^{-p}x_i} (\eps_i^p-1).
$$
The operators $\bar y_i$ are called the {\it Dunkl operators}.

\subhead 1.4.~Combinatorics\endsubhead
Let $\Cc_{m,\ell}$ be the set
of compositions of $m$ with $\ell$ parts, i.e., the set of tuples
$\nu=(\nu_1,\nu_2,\dots\nu_\ell)\in\NN^\ell$ with sum $|\nu|=m$. Let
$\Cc_{m,\ell,n}$ be the subset of compositions whose parts are all
$\geqslant n$. Set
$$J=\{1,2,\dots,m\}.$$
To each $\nu\in\Cc_{m,\ell}$ we associate the
following partition 
$$\gathered
J=J_{\nu,1}\sqcup J_{\nu,2}\sqcup\dots\sqcup J_{\nu,\ell},
\quad J_{\nu,p}=\{i_p,i_p+1,\dots j_p\},\vspace{2mm}
i_p=1+\nu_1+\cdots+\nu_{p-1},\quad j_p=i_{p+1}-1,
\quad p\in\Lambda.
\endgathered$$
We may write $J_p=J_{\nu,p}$ if there is no risk of confusion.
Next, we write
$$\gathered
\CC^\nu_\dom=
\{\l\in\CC^m;\l_i-\l_{i+1}\in\NN,\ \forall i\neq j_1,j_2,\dots\},
\vspace{2mm}
\quad \ZZ^\nu_\dom= \ZZ^m\cap\CC^\nu_\dom,
\quad \ZZ^\nu_\ddom= \rho+\ZZ^ \nu_\dom.
\endgathered
$$
The elements of $\CC^\nu_\dom$ are called the {\it $\nu$-dominant
weights}.
The elements of $\ZZ^\nu_\dom$ are called the {\it $\nu$-dominant
integral weights}.
Let $\Pc_n$ be the set of partitions of $n$,
i.e., the set of non-increasing sequences $\l$ of integers
$\l_1,\l_2,\dots> 0$ with sum $n$. We write ${}^t\!\l$ for the
transposed partition, $|\l|$ for the weight of $\l$, $n(\l)$ for the
integer $\sum_i\l_i(i-1)$ and $l(\l)$ for its length, i.e., for the
number of parts in $\l$.
Let $\Pc=\bigsqcup_n\Pc_n$ be the set of all partitions.
Let $\Pc_n^\ell$ be the set of $\ell$-partitions of $n$. It is the
set of $\Lambda$-tuples $\l=(\l_p)$ of partitions with
$\sum_p|\l_p|=n$. 
Let $\Pc^\ell=\bigsqcup_n\Pc_n^\ell$  be the set of all $\ell$-partitions.
Given any $\ell$-tuple $\nu=(\nu_p)$ in $\NN^\ell$ we set
$$\Pc_{n,\nu}^\ell=\{\l\in\Pc_n^\ell;l(\l_p)\leqslant \nu_p\}.$$
The transpose of a
$\ell$-partition $\l$ is given by
$${}^t\l=({}^t\l_\ell,\dots{}^t\l_2,{}^t\l_1).$$
Any $\ell$-partition $\l\in\Pc_{n,\nu}^\ell$ may be viewed as an
element in $\NN^\nu_{\dom}$ by adding zeroes on the right of each
partition $\l_p$ such that $l(\l_p)<\nu_p$. This yields a bijection
$$\NN^\nu_\dom=\ZZ^\nu_\dom\cap\NN^m=\bigsqcup_n\Pc_{n,\nu}^\ell.\leqno(1.1)$$
Finally, for any $\ell$-tuple
$\nu=(\nu_1,\nu_2,\dots\nu_\ell)$ we'll write
$\nu_p^\circ=(\nu^\circ)_p$ and $\nu^\bullet_p=(\nu^\bullet)_p$
where
$$\nu^\circ=(\nu_\ell,\nu_{\ell-1},\dots\nu_1),\quad
\nu^\bullet=(\nu_{\ell-1},\nu_{\ell-2},\dots\nu_1,\nu_\ell).
$$

\subhead 1.5.~Representations of $W$\endsubhead For
$p\in\Lambda$ there is an unique character
$$\chi_p:\ D_\ell\to\CC^\times,\
\eps\mapsto\eps^p.$$
We set
$$\Irr(\CC\Sen_n)=\{\Xen_\l;\l\in\Pc_n\},
\quad \Irr(\CC W)=\{\Xen_\l;\l\in\Pc_n^\ell\}.\leqno(1.2)$$ If
$\l\in\Pc_n$ then $\Xen_\l$ is defined as in (2.1) below. If
$\l\in\Pc_n^\ell$ then $\Xen_\l$ defined as follows. Any composition
$\mu\in\Cc_{n,\ell}$ can be regarded as a partition
$A=A_{\mu,1}\sqcup\dots\sqcup A_{\mu,\ell}$ as above. We consider
the subgroups
$$\Sen_\mu=\Prod_p\Sen_{A_{\mu,p}}\subset\Sen_n,\quad
W_\mu=\prod_pW_{A_{\mu,p}}\subset W=W_A.$$ Let $w_\mu$ be the
longest element in $\Sen_\mu$. Write $w_0$ for $w_{(n)}$. Fix a
$\ell$-partition $\l=(\l_p)$ in $\Pc_n^\ell$. The tuple
$\mu=(\mu_p)$, with $\mu_p=|\l_p|$ for each $p$, belongs to
$\Cc_{n,\ell}$. Let $\Xen_{\l_p}\chi_{p-1}^{\otimes\mu_p}$ be the
representation of $W_{A_{\mu,p}}$ which is the tensor product of the
$\Sen_{A_{\mu,p}}$-module $\Xen_{\l_p}$ and the one-dimensional
$(D_\ell)^{\mu_p}$-module $\chi_{p-1}^{\otimes\mu_p}$. Then the $W$-module
$\Xen_\l$ is given by
$$\Xen_\l=\Ind_{W_\mu}^{W}\bigl(
\Xen_{\l_1}\chi_\ell^{\otimes\mu_1} \otimes
\Xen_{\l_2}\chi_1^{\otimes\mu_2} \otimes\cdots
\Xen_{\l_\ell}\chi_{\ell-1}^{\otimes\mu_\ell} \bigr),\leqno(1.3)$$ where the
symbol $\Gamma$ denotes the induction.

\subhead 1.6.~The highest weight category of $\Hb_{h,H}$\endsubhead
Let $\Hc_{h,H}$ be the category of
$\Hb_{h,H}$-modules which are locally nilpotent over $\Rb^*$. The
category $\Hc_{h,H}$ is quasi-hereditary by \cite{GGOR}.
The standard
modules of $\Hc_{h,H}$ are the induced modules
$$\Delta_{\l,h,H}=\Ind^{\Hb_{h,H}}_{W\ltimes\Rb^*}(\Xen_\l),\quad
\l\in\Pc_n^\ell.$$ Here $\Xen_\l$ is viewed as a $W\ltimes\Rb^*$-module such
that $y_1,\dots y_n$ act trivially. Let $S_{\l,h,H}$, $P_{\l,h,H}$
denote the top and the projective cover of $\Delta_{\l,h,H}$.

The category $\Hc_{h,H}^{fg}$ consists of the
$\Hb_{h,H}$-modules which are locally nilpotent over $\Rb^*$ and
finitely generated over $\Rb$. The Grothendieck group of
$\Hc_{h,H}^{fg}$ is spanned by the set
$\{[S_{\l,h,H}];\l\in\Pc^\ell_n\}$ and by the set
$\{[\Delta_{\l,h,H}];\l\in\Pc^\ell_n\}$.

The algebra $\Hb_{h,H}$ is given the inner
$\ZZ$-grading such that $x_i$, $y_i$, $w$ have the degrees $1$,
$-1$, $0$ respectively. The sum $\eu=\sum_ix_iy_i+\eu_0$, with
$$\eu_0=
h\sum_{i<j}\sum_{p\in\Lambda}(1-s_{i,j}^{(p)})+
\sum_{i=1}^n\sum_{p,p'=1}^{\ell-1}
\eps^{-pp'}(h_1+\dots+h_{p'})\eps_i^p,$$ is an Euler element for the
grading, i.e., an element $x\in\Hb_{h,H}$ is of degree $i$ iff we have
$[\eu, x]=ix$.

For each $\l\in\Pc_n^\ell$ we denote by $\theta_\l$ the scalar by
which $\eu_0$ acts on $\Xen_\l$. We have the following formula
\cite{R, prop.~6.2}
$$\theta_\l=\ell\sum_{p=2}^\ell|\l_p|(h_1+\cdots +h_{p-1})
-h\ell\sum_{p}\bigl(n(\l_p)-n({}^t\l_p)\bigr)+\theta_0,$$ where
$\theta_0$ is a constant independant of $\l$. 
The partial order $\succcurlyeq$ on the
set of standard modules is the unique order relation such that
\cite{GGOR, sec.~2.5}
$$\Delta_{\mu,h,H}\succ\Delta_{\l,h,H}\iff\theta_\l-\theta_\mu\in\ZZ_\ddom.
\leqno(1.4)$$

\subhead 1.7.~From local systems to $\Hb_{h,H}$-modules\endsubhead Let
$\Rlocnl=\CC[\Alocnl]$, where
$\Alocnl\subset\CC^n$ is the complement
of the hypersurface
$$x_1x_2\cdots x_n\prod_p\prod_{i\neq j}(x_i-\eps^px_j)=0. $$
Note that
$\Alocul=\CC^\times$. For each $\Rb$-module $M$ we write
$$M_{n,\ell}=M\otimes_\Rb\Rlocnl.$$
In particular we have the $\CC$-algebra
$$\Bb=\CC W\ltimes\Rb,\quad\Bb_{n,\ell}=\CC W\ltimes\Rlocnl.\leqno(1.5)$$ 
We'll abbreviate
$\Hb_{h,H,\nl}=(\Hb_{h,H})_\nl$. The algebra $\Hb_{h,H,\nl}$
does not depend on the choice of the parameters $h,H$, up to
canonical isomorphisms. See \cite{GGOR, thm.~5.6} for details. We'll
need the following basic result.

\proclaim{1.8. Proposition} (a) Let $M_{n,\ell}$ be a
$\Bb_{n,\ell}$-module with an integrable $W$-equivariant connection
$\nabla=\sum_i\nabla_i dx_i$. Set
$$\bar y_i=\nabla_i+
k\sum_{j\neq i}\sum_{p}{1\over
x_i-\eps^{-p}x_j}s_{i,j}^{(p)}+\sum_{p\neq 0}{\g_p\over
x_i-\eps^{-p}x_i}\eps_i^p.$$ The assignment $y_i\mapsto\bar y_i$
yields a $\CC$-linear representation of $\Hb_{h,H}$ on $M_{n,\ell}$.

(b) Let $M$ be a $\Bb$-module with an integrable $W$-equivariant
connection. Assume that $M$ is torsion free as a $\Rb$-module and that
the operators $\bar y_1, \bar y_2,\dots \bar y_n$ on $M_{n,\ell}$
preserve the subset $M$.
Then $M$ is a $\Hb_{h,H}$-submodule of $M_{n,\ell}$.
\endproclaim

\noindent{\sl Proof :} Part $(b)$ is obvious. Let us concentrate on
Part $(a)$. Set $\nabla'$ equal to $\sum_i\nabla'_i dx_i$, where
$$\nabla'_i=\nabla_i+
k\sum_{j\neq i}\sum_{p}{1\over x_i-\eps^{-p}x_j}+\sum_{p\neq
0}{\g_p\over x_i-\eps^{-p}x_i}.$$ We have
$$\bar y_i=\nabla'_i+
k\sum_{j\neq i}\sum_{p}{1\over
x_i-\eps^{-p}x_j}(s_{i,j}^{(p)}-1)+\sum_{p\neq 0}{\g_p\over
x_i-\eps^{-p}x_i}(\eps_i^p-1).$$ Thus it is enough to check that
$\nabla'$ is an integrable $W$-equivariant connection and to apply
the same argument as for the Dunkl operators with $\nabla_i'$ instead of
$\partial_{x_i}$. The $W$-equivariance of $\nabla'$ is obvious. If
$i\neq j$ a direct computation yields
$$\aligned
[\nabla'_i,\nabla'_j] &=[\partial_{x_i},k\sum_{r\neq
j}\sum_{p}{1\over x_j-\eps^{-p}x_r}+\sum_{p\neq 0}{\g_p\over
x_j-\eps^{-p}x_j}]+\cr &\quad +[k\sum_{r\neq i}\sum_{p}{1\over
x_i-\eps^{-p}x_r}+\sum_{p\neq 0}{\g_p\over
x_i-\eps^{-p}x_i},\partial_{x_j}]\cr
&=k[\partial_{x_i},\sum_{p}{1\over x_j-\eps^px_i}]+
k[\sum_{p}{1\over x_i-\eps^{-p}x_j},\partial_{x_j}] \cr &=0.
\endaligned$$

\qed

\vskip2cm

\head 2.  ~The affine category $\Oc$ \endhead

According to Suzuki \cite{Su1} the Knizhnik-Zamolodchikov connection
gives a functor from Kazhdan-Lusztig's category of modules over the
affine Lie algebra to the category $\Oc$ of the rational DAHA (for
$\ell=1$). This functor uses affine coinvariants. In this section
we briefly review this construction. Since the results here are not
new, we do not give proofs.

\subhead 2.1.~Lie algebras\endsubhead
Fix once for all an integer $m>0$. We set
$\gen=\gen\len_m(\CC)$ and $G=\GL_m(\CC)$.
To avoid any confusion we may write
$$\gen_J=\gen_m=\gen,\quad G_J=G_m=G$$
Let $U(\gen)$ be the enveloping algebra of $\gen$.
For any $g\in G$, $\xi\in\gen$ let $g\xi$ be the adjoint action of
$g$ on $\xi$ and let ${}^t\xi$ be the transpose of $\xi$. Let
$$\ben\subset\gen,\quad\ten\subset\gen,\quad T\subset G$$
be the Borel Lie subalgebra consisting of upper
triangular matrices and
the maximal
tori consisting of the diagonal matrices.
Let $(\epsilon_i)$, $(\check\epsilon_i)$, $i\in J$, be the canonical
bases of $\ten^*$, $\ten$. 
There is a unique $G$-invariant pairing on $\gen$ such that
$\la\check\epsilon_i:\check\epsilon_j\ra=\delta_{i,j}$ for each
$i,j\in J$.
We have canonical isomorphisms
$\ten=\ten^*=\CC^J$ taking $\check\l$, $\l$ to the tuples
$(\check\l_i)$, $(\l_i)$ with $\check\l_i=\check\l(\epsilon_i)$ and
$\l_i=\l(\check\epsilon_i)$ respectively. The elements of $\ten$,
$\ten^ *$ are called coweights and weights respectively. The weights
in $\ZZ^J$ are called {\it integral weights}. If the weight $\l$ is
given then the symbol $\check\l$ denote the coweight
$\check\l=\sum_{j\in J}\l_j\check\epsilon_j,$ and vice-versa.
We put
$$\rho=(m,\dots,2,1),\quad\a_i=\epsilon_i-\epsilon_{i+1},\quad i\in I,\quad
I=\{1,2,\dots m-1\}.$$ Let $\Pi\subset\ZZ^J$ be the set of roots,
$\Pi^+\subset\Pi$ the set of positive roots containing $\{\a_i;i\in
I\}$ and $\ZZ\Pi$ be the root lattice. Let $e_{k,l}$, $k,l\in J$, be the 
canonical basis vectors of $\gen$. For each simple root $\a_i$
we write $e_{i}=e_{i,i+1}$, $f_{i}=e_{i+1,i}$ for the corresponding root
vectors in $\gen$.
Write $$L(\l)=L(\gen,\l)$$ for the simple $\gen$-module
with highest weight $\l$ (relative to the Borel Lie subalgebra
$\ben$). If $\l$ is integral and dominant we define the $\Sen_n$-module
$$\Xen_\l= H_0(\gen, (\Vb^*)^{\otimes n}\otimes L(\l)),\leqno(2.1)$$
where $H_0$ denotes the space of coinvariants. It vanishes unless $\l\in\Pc_n$.
Compare (1.2).
Let $\Vb$ be the vectorial representation of $G$ and let $\Vb^*$
be the dual module. We can identify $\Vb$ with $\ten$ as $\CC$-vector spaces.
Hence we may view $(\epsilon_j)$ as a basis of $\Vb^*$.
For each $\nu\in\Cc_{m,\ell}$ and $p\in\Lambda$,
let $$\Vb^*_p\subset\Vb^*$$ be the subspace
spanned by the vectors $\epsilon_j$, $j\in J_{\nu,p}$.

\subhead 2.2.~Affine Lie algebras\endsubhead Let $t$ be a formal
variable. For each integer $r$ we set
$$\gathered
\gb=\gen\otimes\CC[t,t^{-1}],\quad \gb_{\sss\geqslant r}=\gen\otimes
t^{r}\CC[t],\quad \gb_{\sss\leqslant r}=\gen\otimes
t^{r}\CC[t^{-1}],\quad\bb=\ben\oplus \gb_\ddom.
\endgathered$$
Let $\hat\gb$ be the central extension of $\gb$ by $\CC$ associated
with the cocycle $(\xi\otimes f,\zeta\otimes
g)\mapsto\la\xi:\zeta\ra\Res_{t=0}(gdf)$. Write $\un$ for the
canonical central element of $\hat\gb$. We abbreviate
$\hat\gb_\dom=\gb_\dom\oplus\CC\un$ and $\hat\bb=\bb\oplus\CC\un$,
the trivial central extensions.
The element $\partial=t\partial_t$ acts on $\gb$ in the obvious way,
yielding a derivation of the Lie algebra $\hat\gb$ such that
$\partial(\un)=0$. We put $\tb=\CC\partial\oplus\ten\oplus\CC\un$,
and $\tilde\gb=\CC\partial\oplus\hat\gb$,
$\tilde\bb=\CC\partial\oplus\hat\bb$.
For any commutative $\CC$-algebra $R$ with 1
we set $\gen_R=\gen\otimes R$, $\hat\gb_R=\hat\gb\otimes R$, etc. We
regard $\gen_R$, $\hat\gb_R$ and $\tilde\gb_R$ as $R$-Lie algebras.

The adjoint $\tb$-action on $\hat\gb$ is diagonalisable.
An element of $\tb^*$ is called an {\it affine weight}.
Let $\hat\Pi\subset\tb^*$
be the set of roots of $\hat\gb$, let $\hat\Pi^+\subset\hat\Pi$ be
the set of roots of the pro-nilpotent radical $\nb$ of
$\hat\bb$, and let $\hat\Pi_\re\subset\hat\Pi$ be the set of real
roots. 
Let $\delta,\omega_0\in\tb^*$ be the linear forms given by
$$\delta(\partial)=\omega_0(\un)=1,\quad
\omega_0(\CC\partial\oplus\ten)=\delta(\ten\oplus\CC\un)=0.$$
The simple roots in $\hat\Pi^+$ are
$$\hat\a_i,\quad i\in \hat I,\quad
\hat I=\{0,1,\dots, m-1\}.$$
From now on we'll use the canonical isomorphism
$\tb^*=\CC^{m+2}=\CC\times\CC^m\times\CC$ such that
$\hat\a_i\mapsto(0,\a_i,0)$ if $i\neq 0$, $\o_0\mapsto (0,0,1)$, and
$\delta\mapsto(1,0,0)$. When there is no risk of confusion we'll
abbreviate $\a_i=\hat\a_i$. Recall that $\a_0=\delta-\epsilon_1+\epsilon_m$.
Let $\la\ :\ \ra$ denote also the
symmetric bilinear form on $\tb^*$ such that
$$\la\o_0:\delta\ra=1,\quad
\la\epsilon_i:\delta\ra=\la\epsilon_i:\o_0\ra=0,\quad
\la\epsilon_i:\epsilon_j\ra=\delta_{i,j}.$$
An affine weight $\l$ such that $\l(\partial)=0$ is called
a {\it classical affine weight}. Let $\tb'\subset\tb^*$ be the set of classical affine weights.

\subhead 2.3.~Enveloping algebras\endsubhead
For any commutative $\CC$-algebra $R$ with 1
let $U(\gen_R)$, $U(\hat\gb_R)$ be the
enveloping algebras of $\gen_R$, $\hat\gb_R$ over $R$. Given an
element $\kappa\in R$ let $U(\hat\gb_R)\to\hat\gb_{R,\kappa}$ be the
quotient by the two-sided ideal generated by the element $\un-c$
where $$c=\kappa-m.$$ We call $c$ the {\it level}. A
$\hat\gb_{R,\kappa}$-module is the same as a $\hat\gb_R$-module such that
the element $\un$ acts as the multiplication by $c$. If $R=\CC$
we'll abbreviate
$\hat\gb_\kappa=\hat\gb_{\CC,\kappa}$, etc.

\subhead 2.4.~Smooth $\hat\gb_R$-modules\endsubhead A
$\hat\gb_{R,\kappa}$-module $M$ is {\it almost smooth}
if each element of $M$ is annihilated by
$\gb_{R,{\sss\geqslant r}}$ for a large enough integer $r$. Let
$\Cc(\hat\gb_{R,\kappa})$ for the category of almost smooth
$\hat\gb_{R,\kappa}$-modules. We may abbreviate
$$\Cc_{R,\kappa}=\Cc(\hat\gb_{R,\kappa}).$$
For each integer $r\geqslant 1$ let
$Q_{R,r}\subset\hat\gb_{R,\kappa}$ be the subspace generated by the
products of $r$ elements of $\gb_{R,{\sss\geqslant 1}}$. Let also
$Q_{R,0}=R$. Given a $\hat\gb_{R,\kappa}$-module $M$ let
$M(r)\subset M$ be the annihilator of $Q_{R,r}$ and let
$M(-r)\subset M$ be the annihilator of $Q_{R,-r}=\sharp(Q_{R,r})$.
Set
$$M(\infty)=\bigcup_{r\geqslant 0} M(r).$$
Note that $M(r)$ is a
$\hat\gb_{R,\dom}$-submodule of $M$ and that $M(\infty)$
is a $\hat\gb_{R,\kappa}$-submodules of $M$. The
$\hat\gb_{R,\kappa}$-module $M$ is called {\it smooth} if
$M=M(\infty)$. All smooth modules are almost smooth. See \cite{KL,
lem.~1.10} for details.

\subhead 2.5.~The Sugawara operator\endsubhead
Put $\xi^{(a)}=\xi\otimes t^a$ for
$\xi\in\gen$. Let $R^\times \subset R$ be the set of units. 
From now on we'll always assume
that $\kappa\in R^\times$. For each integer $b\in\ZZ$ the formal
sum
$$\Lb_b={1\over 2\kappa}\sum_{a\geqslant -b/2}\sum_{k,l\in J}
e_{k,l}^{(-a)}e_{l,k}^{(a+b)}
+{1\over 2\kappa}\sum_{a<-b/2}\sum_{k,l\in
J}e_{k,l}^{(a+b)}e_{l,k}^{(-a)}$$ is called the Sugawara operator.
It lies in a completion of $\hat\gb_{R,\kappa}$.
For any $\hat\gb_{R,\kappa}$-module $M$ and any integer $b$
the Sugawara operator $\Lb_b$ acts on $M(\infty)$ and we have
$$[\Lb_b,\xi^{(a)}]=-a\xi^{(a+b)}.\leqno(2.2)
$$
Let $\Omega=\sum_{k,l\in J}e_{k,l}e_{l,k}$
be the Casimir element in $U(\gen_R)$.
Then
$\Lb_0-\Omega/2\kappa$ acts trivially on the subset $M(1)$.

\subhead 2.6.~Duality for $\hat\gb_R$-modules\endsubhead For each
$\hat\gb_R$-module $M$ we define the $\hat\gb_R$-modules $\sp M$,
$\dg M$, $M^*$, and $M^d$ as follows

\vskip1mm

\item{$\bullet$} $\sp M$ is the $\hat\gb_R$-module
equal to the twist of $M$ by the automorphism
$\sharp:\hat\gb_R\to\hat\gb_R$ such that
$\xi^{(r)}\mapsto(-1)^r\xi^{(-r)}$ and $\un\mapsto-\un$,

\vskip1mm

\item{$\bullet$} $\dg M$ is the $\hat\gb_R$-module equal to
the twist of $M$ by the automorphism $\dag:\hat\gb_R\to\hat\gb_R$
such that $\xi^{(r)}\mapsto-{}^t\xi^{(r)}$ and $\un\mapsto\un$,

\vskip1mm

\item{$\bullet$} $M^*=\Hom_R(M,R)$
with the $\hat\gb_R$-action such that
$(\xi^{(r)}\varphi,m)=-(\varphi,\xi^{(r)}m)$ and
$(\un\varphi,m)=-(\varphi,\un m)$,

\vskip1mm

\item{$\bullet$} $M^d$ is the set of $\gen_{R,\nu}$-finite elements
of $M^*$, i.e., it is the sum of all $\gen_{R,\nu}$-submodules of
$M^*$ which are finitely generated as $R$-modules.

\vskip2mm

The functors $M\mapsto \dg M, \sp M, M^*, M^d$
commute to each other.
We set
$\ddag=\dag\circ\sharp$, 
$DM=(\sp M^*)(\infty)$, and $\dg DM=(\ddg M^*)(\infty)$.
As a $R$-module $\dg DM$ consists of those $R$-linear forms $M\to R$
which are zero on $Q_{R,-r}M$ for some $r\geqslant 1$.

For each $\gen$-module $M$ 
we define the $\gen$-modules $\dg M$,
$M^*$, $M^d$ and the duality functor $\dg D$ 
as in the affine case.

\subhead 2.7.~The (parabolic) category $\Oc$ of $\gen$ and $\hat\gb$\endsubhead 
In this subsection we set 
$R=\CC$. The adjoint $\tb$-action  on $\tilde\gb$ preserves the Lie subalgebra
$\hat\gb$.
By a parabolic Lie subalgebra of $\hat\gb$
we'll mean a $\tb$-diagonalizable Lie subalgebra
$\hat\qb\subset\hat\gb_\dom$ containing a conjugate of $\hat\bb$.
Fix a Levi subalgebra $\hat\lb\subset\hat\qb$. Let
$\Oc(\hat\gb_{\kappa},\hat\qb)$ be the category of the
$\hat\gb_{\kappa}$-modules which are $\hat\lb$-semisimple and
$\hat\qb$-locally finite. 
We abbreviate
$$\hat\Oc_\kappa=\Oc(\hat\gb_\kappa)=\Oc(\hat\gb_\kappa,\hat\bb)$$
with $\hat\lb=\ten\oplus\CC\un$.  Fix a composition $\nu\in\Cc_{m,\ell}$.
Let $\qen_{\nu}\subset\gen$
be the standard parabolic Lie subalgebra  with block
diagonal Levi subalgebra
$$\hen_{\nu}=\gen_{\nu_1}\oplus\dots\oplus\gen_{\nu_\ell}.$$
Let $\uen_{\nu}$ be the nilpotent radical of $\qen_{\nu}$.
Let $\hat\qb_{\nu}$ be the parabolic Lie subalgebra of
$\hat\gb$ with Levi subalgebra
$\hen_{\nu}\oplus\CC\un$, and let $\ub_{\nu}$
be the pronilpotent radical of $\hat\qb_{\nu}$. 
We abbreviate
$$\hat\Oc_{\nu,\kappa}=\Oc(\hat\gb_\kappa,\hat\qb_\nu),\quad
\hat\Oc_{\dom,\kappa}=\Oc(\hat\gb_\kappa,\hat\gb_\dom).$$
In the same way let $\Oc(\gen,\qen_\nu)$ be the category of the
$\gen$-modules which are $\hen_\nu$-semisimple and
$\qen_\nu$-locally finite. 
We abbreviate
$$\Oc=\Oc(\gen,\ben),\quad
\Oc_\nu=\Oc(\gen,\qen_\nu),\quad
\Oc_\dom=\Oc(\gen,\gen).
$$
So $\Oc_\dom$ is the category of all semisimple $\gen$-modules.
Let 
$\qen_{\nu}'=w_0(\qen_\nu)$,
$\hat\qb_{\nu}'=w_0(\hat\qb_\nu)$,
and $\hat\bb'=w_0(\hat\bb)$.

\subhead 2.8.~Induction and generalized Weyl modules\endsubhead  
Set $R=\CC$. For any $\qen_{\nu}$-module $M$ we set
$$M_\nu=\Gamma_{\qen_{\nu}}^{\gen}(M),$$ the induced $\gen$-module. 
An $\hen_{\nu}$-module $M$ may be viewed as a $\qen_{\nu}$-module 
such that $\uen_\nu$
acts as zero, 
yielding again a $\gen$-module $M_\nu$.
Similarly, for any $\hat\qb_{\nu}$-module $M$ of level $c$
we set 
$$M_{\nu,\kappa}=\Gamma_{\hat\qb_{\nu}}^{\hat\gb}(M).$$ 
If $M$ is an $\hen_{\nu}$-module we define the 
$\hat\gb_{\kappa}$-module $M_\nub$ in the
obvious way.

If $\ell=1$ then we have
$\hen_{\nu}=\gen$, 
$\hat\qb_{\nu}=\hat\gb_{\dom}$ and we write
$$M_\kappa=M_\nub. $$
If there is an
integer $r>0$ such that $Q_{r}$ acts by zero on the $\hat\gb_{\dom,\kappa}$-module $M$,
then $M_\kappa$ is called a {\it
generalized Weyl module}. 
More precisely, if $M$ is a $\hat\gb_{\dom,\kappa}$-module
with a finite
filtration by $\hat\gb_{\dom,\kappa}$-modules such that the
successive quotients are annihilated by $Q_{1}$ and lie in
$\Oc_{\nu}^{fg}$ (as $\gen$-modules) then $M_{\kappa}$ belongs to
the category $\hat\Oc_{\nu,\kappa}^{fg}$, and it is called a {\it
generalized Weyl module of type $\nu$}. Further, if $M$ is a simple $\gen$-module 
in $\Oc_\nu$ which is regarded as
a $\hat\gb_{\dom,\kappa}$-module such that $\gb_\ddom$ acts as zero, then $M_\kappa$ is called a 
{\it Weyl module of type $\nu$}.

For each weight $\l\in\ten^*$ let $L(\hen_\nu,\l)$ 
be the simple $\hen_\nu$-module with
highest weight $\l$ (relative to the Borel Lie subalgebra
$\ben\cap\hen_\nu$ of $\hen_\nu$). 
We have the induced $\gen$-module
$$M(\l)_\nu=L(\hen_\nu,\l)_\nu.$$ The top
of $M(\l)_\nu$ is $L(\l)$.
Consider the classical affine weight 
$$\hat\l=\l+c\omega_0.\leqno(2.3)$$ 
We define a $\hat\gb_\kappa$-module by setting
$$M(\hat\l)_\nu=L(\hen_\nu,\l)_\nub.$$
Let $L(\hat\l)$ be the top of $M(\hat\l)_\nu$.
It is the simple $\hat\gb_\kappa$-module with highest weight $\hat\l$.
Recall that $\ten^*$ is identified with $\CC^m$ and that
the elements of $\CC^\nu_\dom\subset\CC^m$ are called 
$\nu$-dominant weights.
If $\l$ is $\nu$-dominant then
$M(\hat\l)_\nu$ is a generalized Weyl module of type $\nu$.
More precisely, if $\l\in\CC^\nu_\dom$ then
$M(\l)_\nu\in\Oc_{\nu}$,
$M(\hat\l)_\nu\in\hat\Oc_{\nu,\kappa}$ 
and they are both called {\it parabolic Verma modules}.
If $\ell=1$ we'll abbreviate 
$$
M(\hat\l)=M(\hat\l)_\nu\in\hat\Oc_{\dom,\kappa}.$$
Although $\hat\Oc_{\nub}$ is not a highest weight category, we'll adopt
the following abuse of language : an object $M$ is said to be
{\it standardly filtered} if it is equipped with a finite 
filtration by submodules such that the successive quotients are 
parabolic Verma modules. Let
$\hat\Oc_{\nub}^{\Delta}\subset\hat\Oc^\f_\nub$ be the full subcategory of
the standardly filtered modules. Let
$\Delta_{\hat\Oc_{\nub}}$ be the set of parabolic Verma modules in 
$\hat\Oc_\nub$.


\proclaim{2.9.~Proposition} Assume that $\kappa\notin\QQ_\dom$.

(a) All modules in $\hat\Oc_{\nu,\kappa}$ are smooth. 
A $\hat\gb_\kappa$-module lies in $\hat\Oc^\f_{\nub}$ 
iff it is a quotient of a generalized
Weyl module of type $\Oc_\nu^\f$.

(b) The category $\hat\Oc_{\nub}^\f$ consists of the finitely generated
smooth $\hat\gb_\kappa$-modules $M$ such that
$M(r)\in\Oc_\nu^\f$ for each $r>0$.

(c) The category $\hat\Oc^\f_{\nub}$ is Abelian, 
any object has a finite length,
and Hom sets are finite dimensional.

(d)
If $M$ lies in 
$\hat\Oc^{fg}_{\nu,\kappa}$ 
then we have $\dg DM=\ddg M^d.$
The functor $\dg D$ takes 
$\hat\Oc^{fg}_{\nu,\kappa}$ to itself, and
it yields an involutive anti-auto-equivalence 
which fixes the simple modules. 
\endproclaim

\noindent{\sl Proof :} 
Clearly, any module in $\hat\Oc_{\nu,\kappa}$ is smooth, and
$M(r)$ lies in $\Oc_{\nu,\kappa}$ for each $M$ in $\hat\Oc_{\nu,\kappa}$
and each $r>0$. Now, assume that $M$ lies in $\hat\Oc_{\nu,\kappa}^\f$. 
There is an integer $r>0$ and a $\hat\gb_{\dom,\kappa}$-submodule
$V\subset M(r)$ such that $M$ is generated by $V$ as a 
$\hat\gb_{\kappa}$-module
and $V\in\Oc_\nu^\f$ as a $\gen$-module.
Then $M$ is a quotient of the generalized Weyl module $V_\kappa$.
Recall that $\Oc_\nu^\f$ is an Abelian category, and that a
subquotient of a $\gen$-module which lies in $\Oc_\nu^\f$ 
lies again in $\Oc_\nu^\f$.
Thus $V$ has a finite filtration by $\hat\gb_{\dom,\kappa}$-submodules such that the successive
quotients are annihilated by $Q_1$ and lies in $\Oc_\nu^\f$.
So $V_\kappa$ is a generalized Weyl module of type $\nu$.
Next, observe that a subquotient of a $\hat\gb$-module which lies in
$\hat\Oc_{\nu,\kappa}$ lies again in $\hat\Oc_{\nu,\kappa}$.
Thus a quotient of a generalized
Weyl module of type $\Oc_\nu^\f$ lies in $\hat\Oc^\f_{\nub}$.
This proves $(a)$.
Part $(b)$ follows from $(a)$ and \cite{Y, thm.~3.5$(1)$}.

Part $(c)$ follows from parts $(a)$, $(b)$.
Indeed, since $\hat\Oc_{\nu,\kappa}$ is obviously Abelian, to prove that 
$\hat\Oc_{\nu,\kappa}^\f$ is Abelian it is enough to prove that 
any submodule of a $\hat\gb$-module in $\hat\Oc_{\nu,\kappa}^\f$ is 
finitely generated.  This follows from \cite{Y, thm.~3.5$(3)$}.
Next, 
we must check that any object has a finite length and that the Hom
sets are finite dimensional. The first claim follows from loc.~cit.,
because $\Oc_\nu^\f$ is Artinian. The second claim is proved as in
\cite{KL, prop.~I.2.29}. It follows from Frobenius reciprocity and
the following three facts : any object in $\hat\Oc^\f_{\nub}$ is a
quotient of a generalized Weyl module of type $\nu$, the Hom
spaces in $\Oc^\f_{\nu}$ are finite dimensional and
$M(r)\in\Oc^\f_\nu$ for each $M\in\hat\Oc^\f_{\nub}$ and each integer
$r>0$.

Part $(d)$ is left to the reader, compare \cite{KL, sec.~1,2}.

\qed

\vskip3mm

Set 
$\hat\Oc'_{\nu,\kappa}=\Oc(\hat\gb_\kappa,\hat\qb'_\nu).$
The functor $\dag$ yields an involutive equivalence
$$\hat\Oc_{\nu,\kappa}\to\hat\Oc'_{\nu,\kappa}.$$
The functor $D$ takes $\hat\Oc^{fg}_{\dom,\kappa}$ into itself, and it
yields an involutive anti-auto-equivalence of
$\hat\Oc^{fg}_{\dom,\kappa}$. 
Note that the modules $\Vb_\kappa$ and
$\Vb_\kappa^*:=(\Vb^*)_\kappa$ are parabolic Verma modules
in $\hat\Oc_{\dom,\kappa}$ such that
$$D(\Vb_{\kappa})=\Vb^*_{\kappa},\quad \Vb^*_{\kappa}=\dg\,\Vb_{\kappa}.$$

From now on we'll always assume that $\kappa\notin\QQ_\dom$.

\subhead 2.10.~The category $\Oc$ of $\hat\gb$ versus $\tilde\gb$\endsubhead 
Set $R=\CC$.
Let $\hat\qb$ be a parabolic Lie subalgebra of $\hat\gb$.
We write $\tilde\qb=\hat\qb\oplus\CC\partial$ and 
$\tilde\lb=\hat\lb\oplus\CC\partial$.
Let 
$\Oc(\tilde\gb_\kappa,\tilde\qb)$ be the category of the
$\tilde\gb_\kappa$-modules which are $\tilde\lb$-semisimple,
$\tilde\qb$-locally finite.
We abbreviate
$$\tilde\Oc_\kappa=\Oc(\tilde\gb_\kappa)=\Oc(\tilde\gb_\kappa,\tilde\bb),
\quad\tilde\Oc_{\nu,\kappa}=\Oc(\tilde\gb_\kappa,\tilde\qb_\nu),\quad
\tilde\Oc_{\dom,\kappa}=\Oc(\tilde\gb_\kappa,\tilde\gb_\dom).$$
Let $\hat\Omega=\partial+\Lb_0$ be the generalized Casimir operator of
$\tilde\gb$, as in \cite{K1, sec.~2.5}. Given a parabolic Lie
subalgebra $\hat\qb\subset\hat\gb$, let
$$\Oc(\tilde\gb_\kappa,\tilde\qb)^0\subset\Oc(\tilde\gb_\kappa,\tilde\qb)$$
be the full subcategory of the modules on which $\hat\Omega$ acts
locally nilpotently. Forgetting the action of $\partial$ yields in
equivalence of categories
$$\Oc(\tilde\gb_\kappa,\tilde\qb)^0\to\Oc(\hat\gb_\kappa,\hat\qb).$$
A quasi-inverse takes the $\hat\gb$-module $M$ to the unique
$\tilde\gb$-module $\tilde M$ which is equal to $M$ as a
$\hat\gb$-module and such that $\partial$ acts as the
semisimplification of the operator $-\Lb_0$. See \cite{S, prop.~8.1}
for details. We identify $M$ and $\tilde M$ if there is no
ambiguity. For any affine weight
$\l\in\tb^*$ the $\l$-weight space
of $M$ (which is identified with $\tilde M$) is
$$M_{\l}=\{x\in\tilde M;hx=\l(h)x,\ \forall h\in\tb\}.$$
For each weight $\l\in\ten^*$ we write
$$\hat\l=\l+c\o_0,\quad
\tilde\l=\hat\l+z_\l\delta,\quad z_\l=-\la\l:2\rho+\l\ra/2\kappa.
$$ We'll
say that the affine weights $\hat\l$, $\tilde\l$ are $\nu$-dominant if
the weight $\l$ is $\nu$-dominant.
The formula for $\Lb_0$ shows that the functor $$\hat\Oc_{\nu,\kappa}\to\tilde\Oc_{\nu,\kappa},
\quad M\mapsto\tilde M$$  takes
$L(\hat\l)$ to the simple module $L(\tilde\l)$ with the ($\nu$-dominant) highest
weight $\tilde\l$ and $M(\hat\l)_\nu$ to the parabolic Verma module
$M(\tilde\l)_\nu$ with the same highest weight. Here $\tilde\l$ is given by the formula above.

\subhead 2.11.~The formal loop Lie algebra\endsubhead 
Fix a commutative $\CC$-algebra $R$ with 1.
Given a finite set $S$ we fix formal
variables $t_i$, $i\in S$, and we set
$$R((t_S))=\Oplus_{i\in S}R((t_i)).$$
Each $f(t)\in R((t))$
yields an element in the $i$-th factor of $R((t_S))$
denoted by
$$f(t)_{[i]}=f(t_i).$$
We set
$$
\Gc_R=\gen\otimes R((t)),\quad
\Gc_{R,S}=\gen\otimes R((t_S)).
$$
As above, if $R=\CC$ we simply forget the subscript $R$ everywhere.
Let $\hat\Gc_{R,S}$ (resp.~$\hat\Gc_R$) be the
central extension of $\Gc_{R,S}$ (resp.~$\Gc_R$) by $R$ associated with the
cocycle $(\xi\otimes f,\zeta\otimes g)\mapsto
\la\xi:\zeta\ra\sum_i\Res_{t_i=0}(gdf).$ Write $\un$ for the
canonical central element of $\hat\Gc_{R,S}$.
Let
$$U(\hat\Gc_{R,S})\to\hat\Gc_{R,S,\kappa}$$
be the quotient by the two-sided ideal
generated by $\un-c$.
Now, let $M_i$, $i\in S$, be almost smooth $\hat\gb_{R,\kappa}$-modules.
We'll use
the following notation : if $\g$ is an operator on $M_i$ then
$\g_{(i)}$ is the operator on  $$\s M=\bigotimes_{i\in S}M_i$$ given by
the action of $\g$ on the $i$-th factor. We'll abbreviate
$$e_{k,l,(i)}=(e_{k,l})_{(i)},\quad \Lb_{b,(i)}=(\Lb_b)_{(i)}.\leqno(2.4)$$
The assignment
$$\xi\otimes f(t)_{[i]}\mapsto\xi\otimes f(t)_{(i)}\leqno(2.5)$$ yields a
representation of $\hat\Gc_{R,S,\kappa}$ on $\s M$.

\subhead 2.12.~The space of affine coinvariants\endsubhead 
Fix a point $\infty\in\PP^1$ and a
coordinate $z$ on $\PP^1-\{\infty\}$. We'll identify
$\PP^1-\{\infty\}$ with $\CC$ and $\CC[\CC]$ with $\CC[z]$. Let
$0\in\PP^1$ be the point such that $z=0$ and 
$\CC^\times=\{z\in\CC;z\neq 0\}$.
Fix an $S$-tuple $x=(x_i;i\in S)$ of
distinct points in $\PP^1$. We set $z_i=z-x_i$ if $x_i\neq\infty$
and $z_i=-z^{-1}$ else.  So $z_i$ is a coordinate on $\PP^1$
centered on $x_i$. Write $$\PP^1_x=\PP^1-\{x_i;i\in S\}.$$ Composing
the expansion of rational functions on $\PP^1$ at $x_i$ with the
assignment $z_i\mapsto t_i$ we get the inclusion
$$\iota_x:\CC[\PP^1_x]\to \CC((t_S)).$$
By the residue theorem the embedding
$$\gen[\PP^1_x]\to\Gc_S,\quad \xi\otimes
f\mapsto\xi\otimes\iota_x(f)$$ lifts to a Lie algebra embedding
$$\gen[\PP^1_x]\to\hat\Gc_S.\leqno(2.6)$$
If $M_i$, $i\in S$, are almost smooth $\hat\gb_\kappa$-modules then
$\hat\Gc_S$ acts on $\s M$. Thus $\gen[\PP^1_x]$ acts also on
$\s M$ through (2.6). 

\proclaim{2.13.~Definition} The space of affine coinvariants of the $M_i$'s
is $\la M_i;i\in S\ra_x=H_0(\gen[\PP^1_x],\s M).$
\endproclaim

\noindent 
Given a point $x_0$ in $\PP^1_x$ we set
$$\hat S=\{0\}\cup S, \quad \hat x=(x_0,x).$$
The following properties are well-known, see e.g., 
\cite{KL, prop.~9.15,9.16,9.18}.

\proclaim{2.14.~Proposition} Assume that
$M_i\in\Cc_\kappa$ for each $i\in S$.

(a) If $M_i=N_{i,\kappa}$ is a generalized Weyl module for each
$i$ then there is a natural isomorphism $\la M_i;i\in
S\ra_x=H_0(\gen,\s N).$

(b) If $M_i\in\hat\Oc_{\dom,\kappa}^\f$ for each $i$
then $\la M_i;i\in S\ra_x$ is finite dimensional.

(c) If $M_0=M(c\o_0)$ there is a natural isomorphism $\la
M_i;i\in \hat S\ra_{\hat x}=\la M_i;i\in S\ra_x.$

(d) The canonical vector space isomorphisms $M_i\to\dg M_i$, $i\in
S$, yield an isomorphism $\la M_i;i\in S\ra_x\to\la \dg
M_i;i\in S\ra_x.$
\endproclaim

\subhead 2.15.~Remark\endsubhead Note that $M(c\o_0)$ is simple, i.e.,
it is equal to $L(c\o_0)$, if $\kappa\in\QQ_{\sss <0}$ by 
\cite{KL, prop.~2.12$(b)$}. 

\vskip3mm

In the rest of Section 2 we'll assume that $S=A\cup\{n+1\}$ and
$x_{n+1}=\infty$.
Now, we allow the tuple $(x_1,\dots x_n)$
to vary in the set $\Alocn$, where $$\Alocn\subset\CC^n$$ is the
complement of the big diagonal. So we may view the $x_i$'s as
regular functions on $\Alocn$, i.e., as elements of the algebra
$\Rlocn=\CC[\Alocn]$.  Let $R\subset\Rb_n(z)$ be the
$\Rb_n$-submodule spanned by the rational functions
$$z_i^{-a}, z^b\ \roman{with}\ i\in A,\ a>0,\ b\geqslant 0.$$
It is a $\Rb_n$-algebra. The assignment
$$x_1,x_2,\dots x_n,z\mapsto x_1,x_2,\dots x_n,x_{n+1}$$
yields a $\Rb_n$-algebra isomorphism $R\to\Rb_{n+1}$. The
expansion at $\{x_i;i\in S\}$ of a rational function
in $R$ yields $\Rlocn$-linear maps
$$R\to\Rlocn((t_S)),\quad\gen_R\to\hat\Gc_{\Rlocn,S}.$$ For 
almost smooth $\hat\gb_\kappa$-modules $M_i$, $i\in S$, the $\Rlocn$-module
$\s M_\Rlocn=(\s M)_\Rlocn$ is equipped with a
$\hat\Gc_{\Rlocn,S}$-action. So we get a $\gen_R$-action on $\s
M_\Rlocn$ such that the element $\xi\otimes z^a$ acts as the sum
$$\sum_{i\in A}(\xi\otimes(t+x_i)^a)_{(i)}+(\xi\otimes(-t^{-1})^a)_{(n+1)}.\leqno(2.7)$$
Consider the $\Rb_n$-module $$\la M_i;i\in
S\ra=H_0(\gen_R,\s M_\Rlocn).$$ 
The following is well-known, see e.g., \cite{KL,
sec.~9.13, prop.~12.12}.

\proclaim{2.16.~Proposition} If
$M_i\in\hat\Oc_{\dom,\kappa}^\f$ for each $i\in S$,
then $\la M_i;i\in S\ra$ is a projective $\Rlocn$-module of finite
rank whose specialization at the point $x\in\Alocn$ is equal to $\la M_i;i\in S\ra_x$. 
\endproclaim

\subhead 2.17.~The space of affine coinvariants of $T(M)_\Rlocn$
\endsubhead 
Now we fix a module
$M\in\Cc_\kappa$. Let $R\subset\Rb_n(z)$ be as above.
We'll apply the construction of
affine coinvariants to the $\hat\gb_\kappa$-modules
$M_1=\dots=M_n=\Vb^*_{\kappa}$ and $M_{n+1}=M.$ Write
$$\Bb_n=\CC\Sen_n\ltimes\Rb_n.$$
There is an obvious representation of $\Bb_n$ on $\s
M_\Rlocn$, such that the group $\Sen_n$ switches the modules
$M_1,M_2,\dots M_n$ and acts on $\Rlocn$, 
and $\Rlocn$ acts by multiplication. This action
centralizes the $\gen_{R}$-action, see formula (2.7). Thus $\la M_i;i\in S\ra$ is
equipped with a representation of $\Bb_n$. To unburden notation
we'll set $$T(M)=(\Vb^*)^{\otimes n}\otimes M.$$ 
The canonical inclusion $\Vb^*\subset
\Vb^*_\kappa$ yields an embedding $$T(M)\subset\s M.$$ Equip the
$\Rlocn$-module 
$T(M)_\Rlocn=T(M)\otimes\Rlocn$ 
with the representation of
$\gen[\CC]$ such that
$$\xi\otimes z^a\mapsto
\sum_{i\neq n+1}\xi_{(i)}\otimes x_i^a+(\sharp\xi^{(a)})_{(n+1)}\otimes 1, 
\quad a\geqslant 0.\leqno(2.8)$$

\proclaim{2.18.~Proposition} The inclusion $T(M)\subset\s M$ yields
a $\Rlocn$-module isomorphism
$$H_0(\gen[\CC],T(M)_\Rlocn)\simeq\la \Vb^*_\kappa,\dots
\Vb^*_\kappa,M\ra.$$ If $\l$ is a dominant weight there is an isomorphism of 
$\Bb_n$-modules
$$H_0(\gen[\CC],T(M(\hat\l))_\Rlocn)\simeq\Gamma^{\Bb_n}_{\Sen_n}
(\Xen_\l).$$
\endproclaim

\noindent{\sl Proof :} We do not give a proof here, since it is rather standard.
Note that part $(a)$ is a particular case of formula (3.15) below (set $\ell=1$),
and part $(b)$ is a particular case of Proposition 3.8$(b)$ below (set again $\ell=1$).
So the proof can be recovered from its twisted version given below.

\qed

\vskip3mm

\subhead 2.19.~The local system of affine coinvariants of $T(M)_\Rlocn$
\endsubhead Let $M$ be as above. Recall that $S=A\cup\{n+1\}$.
We define differential operators on the
$\Rlocn$-module $\s M_\Rlocn$ by the formula
$$\nabla_i=\partial_{x_i}+\Lb_{-1,(i)},\quad\forall i\in A.$$ The
operators $\nabla_i$ commute to each other and give commuting
operators on the $\Rloc_n$-module $T(M)_\Rlocn$, see \cite{BF,
lem.~13.3.7} for instance. The connection
$\nabla=\sum_i\nabla_idx_i$ is called the KZ connection. Let us
compute it. 

For each $i\in S$ we have the $\CC$-linear operator
of $T(M)$ given by
$$e^{(a)}_{k,l,(i)}=(e_{k,l}^{(a)})_{(i)},$$ 
where $e_{k,l}^{(a)}$
acts on $\Vb^*_\Rlocn=\Vb^*\otimes\Rlocn$ 
as the operator $e_{k,l}\otimes x_i^{a}$. 
Next, for $i,j\in A$ the Casimir tensor $\sigma$ yields the 
$\CC$-linear operator on $T(M)$ given by
$$\sigma_{i,j}=\sum_{k,l}e_{k,l,(j)}e_{l,k,(i)}.$$ 
We define the $\Rlocn$-linear operators $\g_i$ on $T(M)_\Rlocn$
by the formula
$$\gathered
\g_{i}=\sum_{j\neq i}\g_{i,j}+\g_{i,n+1},\cr
\g_{i,j}={1\over \kappa}{\sigma_{i,j}\over x_i-x_j},\quad
\g_{i,n+1}={1\over \kappa}\sum_{a>0}\sum_{k,l}
(-1)^{a-1}e^{(a)}_{k,l,(n+1)}e^{(a-1)}_{l,k,(i)}.
\endgathered$$ 
The
following is standard, see \cite{BF, sec.~13.3.8} for instance.

\proclaim{2.20.~Proposition} Under the identification in Proposition
2.18 we have $\partial_{x_i}+\g_i=\nabla_i$.
These operators normalize $\gen[\CC]$ and yield an integrable
$\Sen_n$-equivariant connection on the $\Rloc_n$-module
$H_0(\gen[\CC],T(M)_\Rlocn)$.
\endproclaim

\noindent Note that Propositions 1.8, 2.20 yield a representation of
$\Hb_{h,H}$ on $H_0(\gen[\CC],T(M)_\Rlocn)$.

\vskip2cm

\head 3. Twisted affine coinvariants\endhead

In this section $\ell$ is any integer $>0$. We'll use the
orbifold Knizhnik-Zamolodchikov connection over the configuration
space of the stack $[\PP^1/D_\ell]$. It yields a functor
taking modules over a twisted affine Lie algebra to
$\Hb_{h,H}$-modules. This functor is a generalization of Suzuki's
functor for any $\ell$. As in the case $\ell=1$ we can easily
compute the image of parabolic Verma modules.

\subhead 3.1.~The twisted affine Lie algebra\endsubhead  
Equip the $G$-module $\Vb^*$ 
with a representation of the group $D_\ell$. Since $\Vb^*$ is
the dual of the vectorial representation of $G$, there is a unique
element $g\in G$ which acts on $\Vb^*$ as the generator of $D_\ell$ does.
Let $H\subset G$ be the centralizer of $g$.
We set
$$\gen=\bigoplus_{p\in\Lambda}\gen_p,\quad
\gen_p=\{\xi\in\gen;g\xi=\eps^p\xi\},\quad
\hen=\gen_0.\leqno(3.1)$$
Let $F$ be the automorphism of $\hat\gb$ given
by $$F:\ \hat\gb\to\hat\gb,\ \xi^{(a)}\mapsto\eps^a(g\xi)^{(a)},\
\un\mapsto\un.$$
The twisted affine Lie algebra is the fixed points set
$\hat\gb^F\subset\hat\gb$. We'll say that a $\hat\gb^F$-module is
of level $c$ if the element $\un$ acts as $(c/\ell)\,\Id$. Let
$\hat\gb_\kappa^F$ be the quotient of the enveloping algebra
$U(\hat\gb^F)$ by the two-sided ideal generated by the element
$\un-c/\ell$.

\subhead 3.2.~Modules over the twisted affine Lie algebra\endsubhead
We have the triangular decomposition
$$\hat\gb^F
=\gb_{\sss <0}^F\oplus (\hen\oplus\CC\un)\oplus\gb_{\sss
>0}^F=\gb_{\sss <0}^F\oplus\hat\gb^F_{\dom} .$$
The middle term is a reductive Lie algebra. So the categories
$$\hat\Oc_\kappa^F=\Oc(\hat\gb_\kappa^F),\quad
\hat\Oc^F_{\dom,\kappa}=\Oc(\hat\gb^F_\kappa,\hat\gb^F_{\dom}),
\quad\Cc_\kappa^F=\Cc(\hat\gb^F_\kappa)$$ are defined as in the untwisted case.
The definitions above still make sense by replacing $F$ by the
automorphism
$$F':\ \hat\gb\to\hat\gb,\ \xi^{(a)}\mapsto\eps^{-a}(g\xi)^{(a)},\
\un\mapsto\un.$$ They are indicated by the upperscript $F'$ instead
of $F$. Note that
$$\sharp(\hat\gb_{-\kappa+2m}^{F'})=\hat\gb_\kappa^F,
\quad
\dag(\hat\gb_{\kappa}^{F'})=\hat\gb_\kappa^F.$$ 
In particular
the functor $\dg D$ yields an involutive anti-auto-equivalence
of $\hat\Oc_{\dom,\kappa}^{F,\f}$.
For each $\hat\gb_{\dom}^F$-module $M$ of level $c$ let
$M_\kappa^F$ be the induced $\hat\gb_\kappa^F$-module. The generalized Weyl
$\hat\gb^F_\kappa$-modules are defined as in the untwisted case. 
If $M$ is an $\hen_\nu$-module then it can be equipped with
the unique representation of $\hat\gb_{\dom}^F$ of level $c$ such
that $\gb_{\sss >0}^F$ acts trivially, yielding an induced
$\hat\gb_\kappa^F$-module. We denote it again by $M^F_\kappa$.
For each $\l\in\ten^*$ we put
$$\hat\l=\l+(c/\ell)\o_0.
\leqno(3.2)$$
Compare (2.3).
Consider the $\hat\gb^F_\kappa$-module
$$M(\hat\l)^F=L(\hen_\nu,\l)_\kappa^F.$$
Let $L(\hat\l)^F$ be the top of $M(\hat\l)^F$.
Finally, consider the $\CC$-algebra automorphism
$$\flat:\ \CC[\CC]\to\CC[\CC],\quad
f(z)\mapsto f(\eps^{-1} z),$$ and let $F$ denote also the automorphism
$g\otimes\flat$ of $\gen[\CC]$.

\subhead 3.3.~The space of twisted affine coinvariants of $T(M)$\endsubhead
Recall that for any vector space $M$ we have
$$T(M)=(\Vb^*)^{\otimes n}\otimes M.$$
Consider the following endomorphisms of $T(M)$ 
$$g_i=g_{(i)},\quad
\sigma_{i,j}^{(p)}=\sum_{k,l}(g^pe_{k,l})_{(j)}e_{l,k,(i)}.$$ 
Since $\Vb$ is a $D_\ell$-module, the group $W$
acts on the vector space $(\Vb^*)^{\otimes n}$ in such a way that
$s_{i,j}^{(p)}$, $\eps_i$ act as $\sigma_{i,j}^{(p)}$, $g_{i}$
respectively. This
action gives rise to a representation of $W$ on $T(M).$
The Lie algebra $\hen$ acts on $\Vb^*$ in the obvious way.
Given an $\hen$-module $M$ let $\hen$ act diagonally on $T(M)$.
The representation of $W$ on $T(M)$ factors to a representation of
$W$ on $$\Xen(M)=H_0(\hen,T(M)).\leqno(3.3)$$ The diagonal $\hen$-action and
the $W$-action on $(\Vb^*)^{\otimes n}$ satisfy the double centralizer
property, see \cite{KP, sec.~6} for instance. So $\Xen$ yields a map
$$\Xen:\Irr(U(\hen))\to\Irr(\CC W),$$
where the symbol $\Irr$ denotes the set of isomorphism classes
of finite dimensional modules.
Now, recall the algebras $\Bb$, $\Bnl$ in (1.5).
Given $w\in W$ we'll use the symbol $w$ for the $w$-action on the
$\Bnl$-module $\Gamma^{\Bb_{n,\ell}}_W(T(M))$ and
the symbol $\wu$ for the
operator $\Id\otimes w$ on 
$$T(M)_\Rlocnl=T(M)\otimes\Rlocnl$$
such that $w$ acts on $\Rlocnl$ as in Section 1.2. We'll identify
$$T(M)_\Rb=\Gamma^\Bb_W(T(M)),\quad
T(M)_\Rlocnl=\Gamma^{\Bb_{n,\ell}}_W(T(M))\leqno(3.4)$$ in the
obvious way.
Next, if $M$ is a $\hat\gb^F_\kappa$-module we equip $T(M)_{\Rb_\nl}$ with
the unique $\Rb_\nl$-linear representation of the Lie algebra
$\gen[\CC]^F$ such that
$$\xi\otimes z^a\mapsto
\sum_{i\neq n+1}\xi_{(i)}\otimes x_i^a+(\sharp\xi^{(a)})_{(n+1)}\otimes 1.
\leqno(3.5)$$ 
Note that
$\sharp\xi^{(a)}\in\hat\gb^F$ because
$\xi\otimes z^a\in\gen[\CC]^F$.
This action preserves the subspace $T(M)_\Rb$. 

\proclaim{3.4.~Definition} 
For each $\hat\gb^F_\kappa$-module $M$ we 
define the $\Rb$-module of twisted affine coinvariants of
$T(M)$ by the following formula
$$\CIloc(M)=H_0\bigl(\gen[\CC]^F,T(M)_{\Rb}\bigr).\leqno(3.6)$$
\endproclaim 

\noindent
By formula (3.5) it is quite clear that the $\gen[\CC]^F$-action on $T(M)_{\Rb}$ centralizes the
$\Bb$-action. Thus the representation of $\Bb$ factors to
$\CIloc(M)$, yielding a functor
$$\CIloc\,:\,\hat\gb_\kappa^F\text{-}\modb\to \Bb\text{-}\modb.
\leqno(3.7)$$

\subhead 3.5.~The $\Hb_{h,H}$-action on the space 
of twisted affine coinvariants of $T(M)$\endsubhead 
Now, let $M$ be an almost smooth $\hat\gb^F_\kappa$-module. 
Set $S=A\cup\{n+1\}$.
Let us define new operators on $T(M)_\Rlocnl$. 
For each $i\in S$ we have the $\CC$-linear operator
of $T(M)$ given by
$$e^{(a)}_{k,l,(i)}=(e_{k,l}^{(a)})_{(i)},$$ 
where $e_{k,l}^{(a)}$
acts on $\Vb^*_\Rlocn$ as the operator $e_{k,l}\otimes x_i^{a}$. 
See Section 2.19.
Assume that
$$p\in\Lambda,\quad i\in A,\quad
j\in S,\quad (j,p)\neq(i,0).$$ If $j\neq n+1$ we set
$$\gathered
\g_{i,j}^{(p)}=\sigma_{i,j}^{(p)}/\kappa(x_i-\eps^{-p}x_j),\vspace{2mm}
\g_{i,n+1}^{(p)}={\eps^p\over
\kappa}\sum_{a>0}\sum_r(-1)^{a-1}(F')^{-p}(e_{k,l,(i)}^{(a-1)})e^{(a)}_{l,k,(n+1)}.
\endgathered$$
The sum converges because $M$ is almost smooth. Note that
the action of $$e^{(a)}_{l,k,(n+1)}$$ here does not involves any twist
by the automorphism $\sharp$, contrarily to the
$\gen[\CC]^F$-action on $T(M)_{\Rb_\nl}$ in (3.5). Next, if
$i\neq j$ we set
$$
\g^F_{i,j}=\sum_p\g^{(p)}_{i,j},\quad
\g^F_{i,i}=\sum_{p\neq 0}\g^{(p)}_{i,i}.
$$Finally we consider the
$\CC$-linear operator $\bar y_i$ on $T(M)_\Rlocnl$ given by
$$\bar y_i=\partial_{x_i}-
\sum_{j\neq
i,n+1}\sum_p\g_{i,j}^{(p)}(\underline{s_{i,j}^{(p)}}-1)-\sum_{p\neq
0}\g_{i,i}^{(p)}(\underline{\eps_i^{p}}-1)+\g^F_{i,n+1}.
\leqno(3.8)$$ By (3.4) we have a representation of $\Bnl$ on
$T(M)_\Rlocnl$. This representation can be extended to
a representation of $\Hb_{h,H,\nl}$ as follows.

\proclaim{3.6.~Proposition} Let $M\in\Cc^{F}_\kappa$, $k=-1/\kappa$ and
$\g_p=-\tr(g^{p})/\kappa$.

(a)  The
assignment $y_i\mapsto \bar y_i$ defines a $\Hb_{h,H,\nl}$-action on
the $\Bnl$-module $T(M)_\Rlocnl$ which normalizes the
$\gen[\CC]^F$-action.

(b)  The operator $\bar
y_i-\g^F_{i,n+1}$ vanishes on the subspace $T(M)\subset
T(M)_\Rlocnl$. The representation of $\Hb_{h,H}$ on $T(M)_\Rlocnl$
preserves the subspace $T(M)_{\Rb}\subset T(M)_\Rlocnl$. It factors to
a representation of $\Hb_{h,H}$ on $\CI(M)$.

\endproclaim

Note that the $\Hb_{h,H,n,\ell}$-action on $T(M)_\Rlocnl$ in $(a)$
yields a representation of $\Hb_{h,H,n,\ell}$ for any parameter
$h,H$, since the algebra $\Hb_{h,H,\nl}$ does not depend on the
choice of $h,H$. However the formula (3.8) for the action of $y_i$
does not hold for arbitrary $h,H$.
Hence $(b)$ is false for arbitrary $h,H$.
We'll prove Proposition 3.6 from Section 3.9 onwards.

\subhead 3.7.~
The image by $\CI$ of the parabolic Verma modules
\endsubhead 
Conjugacy classes of elements $g\in G$ such
that $g^\ell=1$ are labeled by integral weights of level $\ell$ in
the dominant  alcove. More precisely, since $g^\ell=1$ there is a
cocharacter $\CC^\times\to G$ such that $\eps\mapsto g^{-1}$. We may
assume that this cocharacter maps into the torus
$T$, so it is identified with
a coweight. We may also assume that its coordinates in the
basis $(\check\epsilon_i)$ are $
(-1)^{\nu_1},(-2)^{\nu_{2}},\cdots(-\ell)^{\nu_\ell}$
for some $\nu\in\Cc_{m,\ell}$ because $\eps^\ell=1$. Then the Lie algebra
$\hen$ in (3.1) is equal to $\hen_\nu$ and
we have
$$g=\bigoplus_p\eps^{p}\,\Id_{\Vb^*_{p}}. $$
From now on we'll always assume that $\hen$, $\nu$, $g$
are as above.
Now we compute the image by $\CI$ of the parabolic Verma modules.

\proclaim{3.8.~Proposition} 
(a) Given a $\nu$-dominant weight $\l$ we have $\Xen(L(\hen_\nu,\l))=0$ if
$|\l|\notin\Pc_{n,\nu}^\ell$ and $\Xen(L(\hen_\nu,\l))\simeq \Xen_{\l^\circ}$
else. If $\nu\in\Cc_{m,\ell,n}$ then
$\Xen\bigl(\Irr(U(\hen_\nu))\bigr)=\Irr(\CC W)$.

(b) For each $\hat\gb^{F}_\dom$-module $M$ the canonical inclusion $M\subset
M_\kappa^F$ yields a $\Bb$-module isomorphism $\Gamma_W^\Bb
(\Xen(M))\to\CIloc(M_\kappa^F).$

\endproclaim

Note that in $(b)$ the symbol $\Xen(M)$ means the functor $\Xen$ in (3.3)
applied to the restriction of $M$ to the Lie subalgebra
$\hen_\nu\subset\hat\gb^{F}_\dom$ and that in $(a)$ the $W$-module
$\Xen_{\l^\circ}$ is as in (1.3).

\vskip3mm

\noindent{\sl Proof :} 
For each $\mu\in\Cc_{n,\ell}$ the subspace
$$\bigotimes_p(\Vb^*_{p})^{\otimes\mu_p}\subset 
(\Vb^*)^{\otimes n}$$ is preserved by
the action of the parabolic subgroup $W_\mu\subset W$ in Section 1.5.
Fix a $W\times U(\hen_\nu)$-module isomorphism
$$(\Vb^*)^{\otimes n}\simeq\bigoplus_{\mu\in\Cc_{n,\ell}}\Ind_{W_\mu}^{W}
\Bigl(\bigotimes_p(\Vb^*_{p})^{\otimes\mu_p}\Bigr).$$
It yields a $W$-module isomorphism
$$\Xen(L(\hen_\nu,\l))\simeq\bigoplus_\mu\Ind_{W_\mu}^W
\Bigl(\bigotimes_pH_0(\gen_{J_p}, (\Vb^*_{p})^{\otimes\mu_p}\otimes
L(\l_p))\Bigr).$$
If $\l_p\in\Pc_{\mu_p}$ then  the $W_{\mu_p}$-module
$H_0(\gen_{J_p},(\Vb^*_{p})^{\otimes\mu_p}\otimes L(\l_p))$
is the tensor product of the
$\Sen_{\mu_p}$-module $\Xen_{\l_p}$ in (2.1) and the
one-dimensional $(D_\ell)^{\mu_p}$-module $\chi^{\otimes\mu_p}_{-p}$.
Else it is zero.
Now fix $\mu=(|\l_p|).$
Then we have
$$\Xen(L(\hen_\nu,\l^\circ))=\Gamma_{W_{\mu^\circ}}^W\bigl(\Xen_{\l_\ell}\chi_{\ell-1}^{\otimes{|\l_\ell|}}\otimes
\Xen_{\l_{\ell-1}}\chi_{\ell-2}^{\otimes{|\l_{\ell-1}|}}\otimes\cdots
\Xen_{\l_1}\chi_{\ell}^{\otimes{|\l_1|}}\bigr).$$
By (1.3) we have also $$\Xen_\l=\Gamma_{W_\mu}^W\bigl(
\Xen_{\l_1}\chi_{\ell}^{\otimes{|\l_1|}}\otimes
\Xen_{\l_2}\chi_{1}^{\otimes{|\l_2|}}\otimes\cdots
\Xen_{\l_\ell}\chi_{\ell-1}^{\otimes{|\l_\ell|}}\bigr).$$ Thus we
have $\Xen(L(\hen_\nu,\l^\circ))\simeq \Xen_{\l}$ as  $W$-modules. Part $(a)$
is proved.
Now, let us prove $(b)$. Set
$$N=M^F_\kappa,\quad \aen=\gen[\CC]^F,\quad
\aen'=(z\gen[\CC])^F,\quad \ben=\hen_{\nu}.$$ The linear map
$$\gen[\CC]\to\hat\gb,\quad\xi\otimes z^b\mapsto\sharp\xi^{(b)}$$
yields a representation of $\aen_\Rb$ on $N_{\Rb}$. View $M_\Rb$ as a
$\ben_\Rb$-module via the restriction to the Lie subalgebra
$\ben_\Rb\subset(\hat\gb_{\Rb,\dom})^F$. We have an isomorphism
of $\aen_\Rb$-modules 
$$N_{\Rb}\simeq\Ind^{\aen_\Rb}_{\ben_\Rb} (M_{\Rb}).$$ 
The
assignment $\xi\otimes z^b\mapsto \sum_ix_i^b\xi_{(i)}$ yields a
representation of $\aen_\Rb$ on $(\Vb^*)^{\otimes n}_{\Rb}$. Since the
induction datum $(\aen_\Rb,\ben_\Rb,\aen'_\Rb,M_{\Rb})$ is $\Rb$-split, 
the tensor
identity yields an isomorphism of $\aen_\Rb$-modules
$$T(N)_{\Rb}\simeq (\Vb^*)^{\otimes n}_\Rb\otimes_\Rb N_\Rb\simeq
(\Vb^*)^{\otimes n}_\Rb\otimes_\Rb \Gamma_{\ben_\Rb}^{\aen_\Rb}(M_\Rb)\simeq
\Ind^{\aen_\Rb}_{\ben_\Rb}(T(M)_{\Rb}),$$
see Section A.3. Therefore we have an isomorphism of $\Bb$-modules
$$\CIloc(N)=H_0(\aen_\Rb,T(N)_\Rb))
\simeq H_0(\ben_\Rb,T(M)_\Rb))
\simeq\Gamma_W^\Bb(\Xen(M)).$$

\qed

\vskip3mm

\subhead 3.9.~Twisted affine coinvariants\endsubhead 
In the rest of Section 3 we'll prove
Proposition 3.6. First we recall what twisted affine coinvariants
are, following \cite{FS}.
The group $D_\ell$ acts faithfully on $\PP^1$ by multiplication, yielding
a cyclic cover $\pi:\PP^1\to\PP^1/D_\ell$ which is ramified at
$0,\infty$. Fix an $S$-tuple of distinct points
$y_i\in\PP^1/D_\ell$ and pick points $x_i$ of
$\pi^{-1}(y_i)$ for each $i\in S$. Let $\ell_i$ be the number of
points in the $D_\ell$-orbit of $x_i$.
Next, fix $S=A\cup\{n+1\}$ and $x_{n+1}=y_{n+1}=\infty$.
Put $$\PP^1_y=\PP^1-\pi^{-1}(\{y_i;i\in S\}).$$
Let $\flat\in\Aut(\CC[\PP^1_y])$ be the comorphism of the action of the
generator of $D_\ell$. We set
$$z_i=z-x_i,\quad z_{n+1}=-z^{-1},\quad z_{i,p}=z-\eps^px_i,\quad
i\neq n+1.$$ We'll abbreviate
$z_{i,p}^a=(z_{i,p})^a$ for each $a$. Note that
$\flat(z)=\eps^{-1}z$, $\flat(z_{n+1})=\eps z_{n+1}$ and $\flat
(z_{i,p})=\eps^{-1}z_{i,p+1}.$ We get the following automorphism
of $\gen[\PP^1_y]$
$$F=g\otimes\flat:\gen[\PP^1_y]\to\gen[\PP^1_y].$$
Let $\Gc_S^F$ be the Lie subalgebra of $\Gc_S$ spanned by
$\Gc_{[n+1]}^F$ and $\Gc_{[i]}$ for $i\neq n+1$. Let
$\hat\Gc_S^F$ be the central extension of $\Gc_S^F$ by $\CC$
associated with the cocycle $(\xi\otimes f,\zeta\otimes g)\mapsto
\la\xi:\zeta\ra\sum_{i\in S}\ell_i\Res_{t_i=0}(gdf).$ Consider the
algebra homomorphism
$$\iota_x:\CC[\PP^1_y]\to\CC((t_S)),\quad
f(z)\mapsto\sum_{i\neq n+1}f(t+x_i)_{[i]}+f(-t^{-1})_{[n+1]}.$$ Here
$f(t+x_i)$, $f(-t^{-1})$ are identified with the corresponding
formal series in $t$.
If $\xi\otimes f, \zeta\otimes g\in \gen[\PP^1_y]^F$
then we have
$$\sum_{i\in S}\ell_i\Res_{z_i=0}(\la\xi:\zeta\ra gdf)=0,\leqno(3.9)$$
because the residue of the one-form
$\la\xi:\zeta\ra gdf$ is the same at each point of the
$D_\ell$-orbit of $x_i$. Thus the linear map
$$\gen[\PP^1_y]^F\to\Gc_S^F,\quad\xi\otimes f\mapsto
\xi\otimes\iota_x(f)$$ lifts to a Lie algebra embedding
$$\gen[\PP^1_y]^F\to\hat\Gc_S^F. \leqno(3.10)$$
Now, assume that $M_1,M_2,\dots M_n\in\Cc_\kappa$ and
$M_{n+1}\in\Cc^F_\kappa.$ The assignment $$\xi\otimes
f(t)_{[i]}\mapsto\xi\otimes f(t)_{(i)}\leqno(3.11)$$ yields a
representation of $\hat\Gc_S^F$ on the tensor product
$$\s M=\bigotimes_{i\in S}M_i.$$
So $\s M$ is also a $\gen[\PP^1_y]^F$-module via the map (3.10).

\proclaim{3.10.~Definition}The space of the twisted affine coinvariants 
of the $M_i$'s is
$$\la M_i;i\in S\ra_x=H_0(\gen[\PP^1_y]^F,\s M).$$
\endproclaim

\noindent 
Now, let the point $x$
vary in $\Alocnl$. View the $x_i$'s as coordinates on
$\Alocnl$. Let $R\subset\Rb_\nl(z)$ be the $\Rnl$-submodule
spanned by
$$z_{i,p}^{-a},\ z^b\ \roman{with}\  a>0,\ i\in A,\ p\in\Lambda,
\  b\geqslant 0.$$
It is
closed under multiplication. Let $R'\subset R$ be the
$\Rlocnl$-subalgebra consisting of the functions which vanish on
$\{z=\infty\}$. As a $\Rlocnl$-module $R'$ is spanned by
$$z_{i,p}^{-a}\ \roman{with}\ a>0,\ i\in A,\ p\in\Lambda.$$
There is an unique $\Rlocnl$-algebra automorphism $\flat$ of $R$,
$R'$ such that $f(z)\mapsto f(\eps^{-1} z)$. It yields the automorphism
$F=g\otimes\flat$ of the $\Rlocnl$-Lie algebras $\gen_R$,
$\gen_{R'}$.
The construction above yields the
$\Rlocnl$-Lie algebra $\hat\Gc^F_{\Rb_\nl,S}$
and the morphism of $\Rlocnl$-Lie algebras
$$(\gen_R)^F\to\hat\Gc^F_{\Rb_\nl,S}\leqno(3.12)$$
which is analogous to (3.10). Assume that $M_1,M_2,\dots
M_n\in\Cc_\kappa$ and $M_{n+1}\in\Cc^F_\kappa.$
The assignment (3.11) yields a representation of
$\Gc^F_{\Rb_\nl,S}$ on $\s M_\Rlocnl$. Thus (3.12) yields a
representation of $(\gen_R)^F$ on $\s M_\Rlocnl$. Taking the
coinvariants with respect to the Lie algebras $(\gen_R)^F$,
$(\gen_{R'})^F$ we get the following two $\Rlocnl$-modules :
$$\la M_i;i\in S\ra,\quad\la M_i;i\in S\ra'.$$

We'll also use a more general version of twisted affine coinvariants
obtained by inserting a new module
$M_0\in\Cc^{F'}_\kappa$  at the point 0.
Set $$\hat S=\{0\}\cup S=A\cup\{0,n+1\}.$$
Let $\hat x=(0,x)$ and $\hat y=(0,y)$.
The Lie algebra $\gen[\PP^1_{\hat y}]^F$ acts on the tensor product
${}_{\hat S}M$ in the obvious way. We define 
$$\la M_i;i\in \hat S\ra_{\hat x}=
H_0(\gen[\PP^1_{\hat y}]^F,{}_{\hat S} M).$$

\proclaim{3.11.~Proposition}
(a) If $M_1,\dots
M_n\in\hat\Oc_{\dom,\kappa}^\f$ and
$M_{n+1}\in\hat\Oc_{\dom,\kappa}^{F,\f}$ the
$\Rb_\nl$-module $\la M_i;i\in S\ra$ is projective of finite rank.
Its fiber at the point
$x\in C_\nl$ is equal to the $\CC$-vector space
$\la M_i;i\in S\ra_x$.

(b) Let $\hat x=(0,x)\in C_{n+1,\ell}$, $M_1,\dots
M_n\in\Cc_\kappa$ and $M_{n+1}\in\Cc^F_\kappa$.
Set $M_0=M(c\o_0)^{F'}$. The obvious
inclusion $L(\hen_\nu,0)\subset M_0$ yields an isomorphism 
$\la M_i;i\in \hat S\ra_{\hat x}=\la M_i;i\in S\ra_x$.
\endproclaim

\noindent{\sl Proof :} The proof is quite standard, so we'll be very brief.
Part $(a)$ is the twisted version of Proposition 2.16. 
The proof that the $\Rb_\nl$-module $\la M_i;i\in S\ra$ 
is finitely generated is the same as in the untwisted case, 
see e.g., \cite{KL, prop.~9.12}. 
Indeed, observe first that
for each $i\in A$ there is an integer $r_i>0$ such that $M_i$ 
is generated by $M_i(r_i)$ as a $\hat\gb$-module.
Similarly there is an integer $r_{n+1}>0$ such that $M_{n+1}$ 
is generated by $M_{n+1}(r_{n+1})$ as a $\hat\gb^F$-module,
where $M_{n+1}(r_{n+1})$ is the set of elements of $M_{n+1}$ killed by
any product of $r_{n+1}$ elements of $\gb_\ddom^F$. 
Next, set 
$$f_i=z_{i,0}^{-1},\quad f_{n+1}=-z,\quad i\in A.$$
Then $f_1,f_2,\dots,f_{n+1}\in R$ and for each $i,j\in S$
the expansion of $f_i$ at $x_j$ is $t^{-1}$ if $i=j$ and it belongs to
$\Rb_\nl[[t]]$ else.
The rest of the proof is by induction, using the spaces
$M_i(r_i)$ and the functions $f_i$ as in loc.~cit.
Next, the $\Rb_\nl$-module $\la M_i;i\in S\ra$ is projective,
because it is finitely generated and admits a connection, 
compare \cite{KL, prop.~12.12}.
This connection is called the {\it orbifold KZ connection}. 
It is constructed in a quite general setting in \cite{Sz}.
We have given the construction of the orbifold KZ connection in the 
second step of the proof of Proposition 3.6, see Section 3.13 below.
Finally, the third claim in $(a)$ is also
proved as in the untwisted case, see e.g., \cite{KL, sec.~9.13}.
It is a twisted analogue of the commutation of affine 
coinvariants with base change.

Part $(b)$ is the twisted analogue of the {\it propagation of vacua} recalled in
Proposition 2.14$(c)$.
It is proved as in the untwisted case in \cite{KL, prop.~9.18}.
More precisely, the expansion at $\hat y$ yields a Lie algebra embedding
$$\gen[\PP^1_{\hat y}]^F\to\hat\Gc_{\hat S}^F,$$
where $\hat\Gc_{\hat S}^F$
is a central extension of the Lie subalgebra  
$\Gc_{\hat S}^F\subset\Gc_{\hat S}$
spanned by $\Gc_{[0]}^{F'}$, $\Gc_{[i]}$, $i\in A$, and 
$\Gc_{[n+1]}^F$.
Then $\la M_i;i\in \hat S\ra_{\hat x}$ is the space of coinvariants
of the representation of $\gen[\PP^1_{\hat y}]^F$ on the space
$${}_{\hat S}M=M(c\o_0)^{F'}\otimes{}_SM.$$
Note that the $\hat\Gc^F_{\hat S}$-module ${}_{\hat S}M$ is induced from the
$\hat\Gc_{\hat S}^{F,+}$-module $\CC\otimes{}_SM,$
where
$\hat\Gc_{\hat S}^{F,+}\subset\hat\Gc^F_{\hat S}$
is the central extension of the Lie subalgebra  
$\Gc_{\hat S}^{F,+}\subset\Gc_{\hat S}^F$
spanned by $(\gen\otimes\CC[[t_0]])^{F'}$, $\Gc_{[i]}$, $i\in A$, and 
$\Gc_{[n+1]}^F$.
Therefore the claim follows from the tensor identity as in loc.~cit., and from
the following equality
$$
\hat\Gc_{\hat S}^{F}=\hat\Gc_{\hat S}^{F,+}+
\gen[\PP^1_{\hat y}]^F,$$
which is left to the reader.
Compare Proposition 4.5 below.

\qed

\vskip3mm

\subhead 3.12.~Remark\endsubhead
The $\Rb_\nl$-module of twisted affine coinvariants
$\la M_i;i\in S\ra$ and the $\Rb$-module 
$\Cen(M)$ of twisted affine coinvariants of $T(M)$ are related by the following 
$\Rb_\nl$-module isomorphism, see (3.15) below 
$$M_1=\cdots=M_n=\Vb^*_\kappa,\quad M_{n+1}=M\in\Cc_\kappa^F
\ \Rightarrow\ 
\la M_i;i\in S\ra\simeq\Cen(M)_\nl.$$

\subhead 3.13.~Proof of Proposition 3.6\endsubhead 
Recall that $M\in\Cc^F_\kappa$. 
Let $R'\subset R\subset\Rlocnl(z)$ be as in Section 3.9.
To avoid
cumbersome notation we write
$$
\aen=(\gen_R)^F,\quad\aen'=(\gen_{R'})^F,\quad
\ben=(\gen_{\Rb_\nl[z]})^F.$$ We have
$\aen=\aen'\oplus\ben$, and the
$\Rb_\nl$-modules $\aen'$, $\ben$ are respectively spanned by
\vskip1mm
\itemitem{$\bullet$}
$\xi\otimes z^p(z^\ell-x_i^\ell)^{-a}\ \roman{with}\
\xi\in\gen_p,\ a>0,\ p=0,\dots\ell-1,\
i\in A,$
\vskip1mm
\itemitem{$\bullet$}
$\xi\otimes z^b\ \roman{with}\ \xi\in\gen_b,\ b\geqslant 0.
$
\vskip1mm

\noindent The proof is long and technical. We'll split it into several steps.

\vskip2mm

\noindent{\sl Step 1 :}
First, we prove the formula
(3.15) below, which identifies $T(M)_\Rlocnl$ with a set of twisted affine
coinvariants. To do so we put 
$$M_1=\cdots=M_n=\Vb^*_\kappa,\quad M_{n+1}=M.$$ 
To simplify the notation we'll set $n=1$ in this part
of the proof. The case $n>1$ is identical. Therefore we have
$$T(M)=\Vb^*\otimes M,\quad \s M=\Vb^*_\kappa\otimes M$$
and the inclusion $\Vb^*\subset
\Vb^*_\kappa$ yields a $\Rb_{1,\ell}$-linear embedding
$$T(M)_{\Rb_{1,\ell}}\subset
\s M_{\Rb_{1,\ell}}.\leqno(3.13)$$ Recall that $\Rb_{1,\ell}=\CC[x_1^{\pm 1}]$ and that
$\aen$ acts on $\s
M_{\Rb_{1,\ell}}$ via the map (3.12). Thus, by definition of the map
$\iota_x$ the element $\xi\otimes z^b\in\ben$ acts as the sum
$$\xi\otimes(t+x_1)^b_{[1]}+\xi\otimes(-t^{-1})^b_{[2]}.$$ Since the
subspace $\Vb^*\subset \Vb^*_\kappa$ is killed by $\gb_{\sss >0}$, we get
that $\xi\otimes z^b$ acts on the subspace $T(M)_{\Rb_{1,\ell}}$ as the sum
$$(x_1^b\xi)_{(1)}+(\sharp\xi^{(b)})_{(2)}.$$ This is precisely the
$\ben$-action on $T(M)_{\Rb_{1,\ell}}$ given in (3.5). Therefore,
the inclusion (3.13) is an embedding of $\ben$-modules.

Now we prove that there is an isomorphism of $\aen$-modules
$$\s M_{\Rb_{1,\ell}}
\simeq\Ind^{\aen}_{\ben}
 (T(M)_{\Rb_{1,\ell}}).
\leqno(3.14)$$ Let $\hat\ben=\ben\oplus\Rb_{1,\ell}$ (the trivial
central extension) and let $\hat\aen$ be the central extension of
$\aen$ by $\Rb_{1,\ell}$ associated with the cocycle $(\xi\otimes
f,\zeta\otimes g)\mapsto \la\xi:\zeta\ra\Res_{z=x_1}(gdf). $ 
By the residue theorem there is a $\Rb_{1,\ell}$-Lie algebra homomorphism
$$\hat\aen\to\hat\gb^F_{\Rb_{1,\ell}},\quad\xi\otimes f(z)\mapsto \xi\otimes f(-t^{-1}),
\quad\un\mapsto -\ell\un,$$
compare (3.9).
Thus, since $M$ is a $\hat\gb^F$-module of level $c$, the assignment
$$\xi\otimes f(z)\mapsto \xi\otimes f(-t^{-1})$$
yields a representation of $\hat\aen$ on $M_{\Rb_{1,\ell}}$ of level $-c$.
Similarly, since $\Vb^*_\kappa$ is a $\hat\gb^F$-module of level $c$, the assignment
$$\xi\otimes f(z)\mapsto \xi\otimes f(t+x_1)$$ 
yields a representation of $\hat\aen$ of level $c$
on $(\Vb^*_\kappa)_{\Rb_{1,\ell}}$. 
The representation of $\aen$ on $\s M_{\Rb_{1,\ell}}$
is the restriction of the tensor product of the $\hat\aen$-modules
$M_{\Rb_{1,\ell}}$ and $(\Vb^*_\kappa)_{\Rb_{1,\ell}}$. 
The $\Rb_{1,\ell}$-submodule
$\Vb^*_{\Rb_{1,\ell}}\subset(\Vb^*_\kappa)_{\Rb_{1,\ell}}$ is preserved by
the $\hat\ben$-action. So the quadruple $(\hat\aen, \hat\ben, \aen',
\Vb^*_{\Rb_{1,\ell}})$ is a $\Rb_{1,\ell}$-split induction datum. See
Section A.3 for details.

We claim that the representation of $\hat\aen$ on
$(\Vb^*_\kappa)_{\Rb_{1,\ell}}$ is isomorphic to the induced module
$\Ind^{\hat\aen}_{\hat\ben}(\Vb^*_{\Rb_{1,\ell}})$. Then (3.14) follows
from the tensor identity, because we have
$$\s M_{\Rb_{1,\ell}}
=\Gamma_{\hat\ben}^{\hat\aen}(\Vb^*_{\Rb_{1,\ell}})\otimes_{\Rb_{1,\ell}}
M_{\Rb_{1,\ell}}=
\Gamma_{\hat\ben}^{\hat\aen}(\Vb^*_{\Rb_{1,\ell}}\otimes_{\Rb_{1,\ell}}
M_{\Rb_{1,\ell}})=\Gamma_{\ben}^{\aen}(T(M)_{\Rb_{1,\ell}}).
$$
Now we prove the claim. The $\hat\aen$-action on $(\Vb^*_\kappa)_{\Rb_{1,\ell}}$
yields a map
$$\gathered
\Ind^{\hat\aen}_{\hat\ben}(\Vb^*_{\Rb_{1,\ell}})\to
(\Vb^*_\kappa)_{\Rb_{1,\ell}},\quad
(\xi\otimes f(z))\otimes v\mapsto (\xi\otimes f(t+x_1)) v.
\endgathered $$
We must prove that it is invertible.
Recall that 
$$\gathered
\Ind^{\hat\aen}_{\hat\ben}(\Vb^*_{\Rb_{1,\ell}})=U(\aen')\otimes_{\Rb_{1,\ell}}\Vb^*_{\Rb_{1,\ell}},
\vspace{2mm}
(\Vb^*_\kappa)_{\Rb_{1,\ell}}=\Ind^{\hat\gb}_{\hat\gb_\dom}(\Vb^*_{\Rb_{1,\ell}})
=U(\gb_{_{\Rb_{1,\ell},{\sss <0}}})\otimes_{\Rb_{1,\ell}}\Vb^*_{\Rb_{1,\ell}}.
\endgathered$$
Further, for each $a>0$ the assignment
$z\mapsto t+x_1$ takes $z^p(z^\ell-x_1^\ell)^{-a}$ to a formal
series in $$(\Rb_{1,\ell})^\times t^{-a}+\Rb_{1,\ell}[[t]]t^{1-a}.$$ 
Therefore the claim is obvious.

Note that (3.14) yields the second of the two isomorphisms below
$$
\la M_i;i\in S\ra'\simeq T(M)_\Rlocnl,\quad \la M_i;i\in
S\ra\simeq\CIloc(M)_{n,\ell}.\leqno(3.15)$$ Indeed, with the
notation used in this proof the definition (3.6) yields
$$\CIloc(M)_{n,\ell}=H_0(\ben,T(M)_\Rnl).$$
The proof of the first isomorphism is identical because the
quadruple $(\hat\aen,\hat\ben, \aen', T(M)_{\Rb_{1,\ell}})$ is also
a $\Rb_{1,\ell}$-split induction datum, hence (3.14) and Section A.3
imply that
$$\s M_{\Rb_{1,\ell}}\simeq\Ind^{\aen'}(T(M)_{\Rb_{1,\ell}}).$$

\vskip2mm

\noindent{\sl Step 2 :}
Now we define an integrable
connection on the $\Rb_\nl$-modules $T(M)_\Rlocnl$ and $\CI(M)_{n,\ell}$.
For each $i$ we consider the differential operator on $\s M_\Rlocnl$ given by
$$\nabla_i=\partial_{x_i}+\Lb_{-1,(i)}.$$ Next we set
$\nabla=\sum_i\nabla_i\,dx_i$. It is an integrable connection on the
$\Rlocnl$-module $\s M_\Rlocnl$. By (2.2), for each $\xi\in\gen$ and
$f\in\Rlocnl((t_S))$ we have
$$[\nabla_i,\xi\otimes f]=\xi\otimes(\partial_{x_i}f-\partial_{t_i}f).
\leqno(3.16)$$ We claim that $\nabla$ normalizes the $\aen$-action
given by (3.12). Hence the isomorphisms (3.15) yield an integrable
connection on the $\Rb_\nl$-module $T(M)_\Rlocnl$ which factors to a
connection on $\CI(M)_{n,\ell}$. Write $\nabla$ again for these
connections.

Now we prove this claim. For each integers $a$, $b$ with $a> 0$ we
define the constant $c_a^b$ by the following formula
$$\partial_t^{a-1}(t^b)/(a-1)!=c_a^bt^{b-a+1}.$$ A direct computation
yields the following relations in $\Rlocnl((t_S))$

$$
\iota_x (z_{j,p}^{-a})= \delta_{p,0}\,t^{-a}_{[j]}
-\sum_k\sum_{b\geqslant 0}c_{a}^{-b-1}
(\eps^px_j-x_k)^{-b-a}t^b_{[k]}-\sum_{b\geqslant
0}(-1)^bc^b_{a}(\eps^px_j)^{b-a+1}t^{b+1}_{[n+1]},$$ $$
\iota_x(z^b)= \sum_k(x_k+t)_{[k]}^b+(-1)^bt^{-b}_{[n+1]},$$ where
$k=1,2,\dots n$ and $(k,p)\neq(j,0)$. Thus the derivation
$\partial_{x_i}-\partial_{t_i}$ annihilates $\iota_x (z^b)$ and
takes $\iota_x (z_{j,p}^{-a})$ to $ \delta_{i,j}\,\eps^pa\,\iota_x
(z_{j,p}^{-a-1})$. So
$(\partial_{x_i}-\partial_{t_i})\circ\iota_x=\iota_x\circ\partial_{x_i}$
on $R$.
Hence (3.16) yields the following relations
$$\gathered
[\nabla_i,\xi\otimes\iota_x (z^b)]=0,\vspace{2mm}
[\nabla_i,\xi\otimes\iota_x (z^p(z^\ell-x_j^\ell)^{-a})]=
\delta_{i,j}\,a\ell x_i^ {\ell-1}\,\xi\otimes\iota_x
(z^p(z^\ell-x_j^\ell)^{-a-1}).\endgathered$$ The claim is proved.

\vskip2mm

\noindent{\sl Step 3 :} Now we claim that the connection $\nabla$ on
$T(M)_{\Rnl}$ is $W$-equivariant and
$$\bar y_i= \nabla_i+k\sum_{j\neq i}
\sum_{p}{1\over x_i-\eps^{-p} x_j}s_{i,j}^{(p)}+\sum_{p\neq
0}{\g_p\over x_i-\eps^{-p}x_i}\eps_i^p.\leqno(3.17)$$ Thus part (a)
of the proposition follows from Proposition 1.8$(a)$. Note that the
$\Hb_{h,H}$-action on $T(M)_{\Rb_\nl}$ we have just constructed
normalizes the $\ben$-action given by (3.5), because the connection
$\nabla$ normalizes the $\aen$-action on $\s M_\Rlocnl$ given by
(3.12) by the previous claim.

Now, let us prove the latest claim. We'll compute the connection
$\nabla$ on $T(M)_{\Rb_\nl}$. Since $\nabla_i$ is a derivation, it
is enough to consider its restriction to the subspace $T(M)\subset
T(M)_{\Rb_\nl}$. Under the quotient map $\s M_\Rlocnl\to
T(M)_\Rlocnl$ in (3.15), the obvious inclusion $T(M)\subset \s
M_\Rlocnl$ is taken to the obvious inclusion $T(M)\subset
T(M)_\Rlocnl$. So it is enough to compute the action of the operator
$\Lb_{-1,(i)}$ on the subspace $T(M)\subset\s M_\Rlocnl$.

For $\xi\in\gen$ we consider the following elements
$$\sigma(\xi,i)=\sum_p\eps^{p}g^p\xi\otimes z_{i,p}^{-1}\in\aen',\quad
\xi_{[j]}^{(b)}=\xi\otimes t^b_{[j]}\in\hat\Gc_{S,\Rb_\nl}.$$ Note
that
$$\iota_x(\sigma(\xi,i))=\xi_{[i]}^{(-1)}-\sum_{j,p,b}
\eps^{p}(\eps^px_i-x_j)^{-b-1}(g^p\xi)_{[j]}^{(b)}
-\sum_{b,p}(-1)^{b}\eps^{p}(\eps^px_i)^{b}(g^p\xi)^{(b+1)}_{[n+1]},$$
where $j\in A$, $b\geqslant 0$, $p\in\Lambda$ and
$(j,p)\neq(i,0)$ in the first sum. Let $v\in T(M)$, viewed as a
subspace of $\s M_{\Rb_\nl}$. Then, we have
$$\xi_{[i]}^{(-1)}v=
\sum_{j,p}{g^p\xi_{(j)}\over
x_i-\eps^{-p}x_j}v+\sum_{b,p}(-1)^bx_i^bF^p(\xi^{(b+1)})_{(n+1)}v\
\mod\ \aen'(\s M_{\Rb_\nl}).$$ 
Note that $\xi_{(i)}v=0$ for each $\xi\in\hat\gb_\ddom$, $i\in A$.
We have
$$\aligned
&\Lb_{-1,(i)}v={1\over\kappa}\sum_{k,l}e^{(-1)}_{k,l,(i)}e_{l,k,(i)}v,
\quad i\in A,\vspace{2mm}
&\g_{i,j}^{(p)}= \cases {\ds{1\over \kappa}\sum_{a>0}\sum_{k,l}
(-1)^{a-1}x_i^{a-1}F^p(e^{(a)}_{k,l})_{(n+1)}e_{l,k,(i)}}&\roman{if}\
j=n+1\hfill\cr
{\ds\sum_{k,l}(g^pe_{k,l})_{(j)}e_{l,k,(i)}/\kappa(x_i-\eps^{-p}x_j)}
&\roman{else}.\endcases \endaligned$$ Setting $\xi=e_{k,l}$ in the
formula above, we get
$$\Lb_{-1,(i)}v=
\g^F_{i,1}v+\cdots\g^F_{i,n+1}v\ \mod\ \aen'(\s M_{\Rb_\nl}).$$ In
other words, the connection $\nabla$ on $T(M)_{\Rb_\nl}$ is given by
$$\nabla_i=
\partial_{x_i}+{1\over\kappa}\sum_{j\neq i}\sum_{p}{\sigma_{i,j}^{(p)}\over x_i-\eps^{-p} x_j}
+{1\over\kappa} \sum_{p\neq 0}{\sigma_{i,i}^{(p)}\over
x_i-\eps^{-p}x_i}+\g^F_{i,n+1},$$ with $j\in A$. Now, for
$p\neq 0$ the operator $\sigma_{i,i}^{(p)}$ acts on $T(M)_{\Rb_\nl}$
as the operator $\tr(g^{p})g^p_{i}= -\kappa\g_pg^p_{i}$.
Thus we have
$$\nabla_i=
\partial_{x_i}+{1\over\kappa}\sum_{j\neq i}\sum_{p}{\sigma_{i,j}^{(p)}\over x_i-\eps^{-p} x_j}
-\sum_{p\neq 0}{\g_p\,g_i^p\over x_i-\eps^{-p}x_i} +\g^F_{i,n+1}.$$
Since $k=-1/\kappa$ the right hand side of (3.17) is
$$\gathered
\partial_{x_i}-{1\over\kappa}\sum_{j\neq i}\sum_{p}{1\over x_i-\eps^{-p} x_j}
(s_{i,j}^{(p)}-\sigma_{i,j}^{(p)})+ \sum_{p\neq 0}{\g_p\over
x_i-\eps^{-p}x_i}(\eps_i^p-g_i^{p})+\g^F_{i,n+1}=\cr
\partial_{x_i}-
\sum_{j\neq i}\sum_p\g_{i,j}^{(p)}
(\underline{s_{i,j}^{(p)}}-1)-\sum_{p\neq
0}\g_{i,i}^{(p)}(\underline{\eps_i^{p}}-1)+\g^F_{i,n+1}=\bar y_i.\endgathered$$
So (3.17) is proved. The $W$-equivariance of $\nabla$ is obvious
from the formula above.

Let us  prove $(b)$. If $v\in T(M)$ then (3.8) gives $\bar
y_iv=\g^F_{i,n+1}v.$ The subspace $T(M)_{\Rb}\subset
T(M)_{\Rb_{n,\ell}}$ is preserved by the $\Bb$-action in (3.4). By
(3.8) it is also preserved by the operators $\bar y_1,\bar y_2,\dots
\bar y_n$, i.e., the denominator in $\g_{i,j}^{(p)}$ simplifies.
Thus $T(M)_\Rb$ is a $\Hb_{h,H}$-submodule. The representation of
$\Hb_{h,H}$ on $T(M)_\Rb$ normalizes the $\gen[\CC]^F$-action by
part $(a)$. Hence $\Hb_{h,H}$ acts on $\CI(M)$.

\qed

\vskip3mm

\subhead 3.14.~Remarks\endsubhead 
$(a)$ The connection $\nabla$ is a particular case of the orbifold
KZ connection in \cite{Sz}.

$(b)$ In Section 4.4 we'll also use a more general $\Hb_{h,H}$-module
$\Cen(M',M)_\nl$ obtained by inserting also an almost smooth
$\hat\gb_\kappa^{F'}$-module $M'$ at the point $0\in\PP^1$ before
taking twisted affine coinvariants. See also Section 3.9. Then
Proposition 3.8 generalizes in the following way : if
$M'=(E')^{F'}_\kappa$ and $M=E^F_\kappa$ then we have an isomorphism
of $\Bnl$-modules
$$\Cen(M',M)_\nl\simeq\Gamma^\Bnl_W(\Xen(E'\otimes E)).$$


\vskip2cm

\head 4. Untwisting the space of twisted affine
coinvariants\endhead

Any twisted affine Lie algebra associated with an inner automorphism
is isomorphic to a non-twisted one. In this section we prove that,
in a similar way, the space of twisted affine coinvariants can
expressed in terms of non-twisted affine coinvariants. Recall,
however, that twisted affine coinvariants yield a local system
over the configuration space of the stack $[\PP^1/D_\ell]$
while affine coinvariants yield a local system over the
configuration space of  $\PP^1$.

\subhead 4.1.~Notation\endsubhead Given a weight $\pi$ we may 
shift the origin in $\ten^*$ and write
$$\l_\pi=\l+\pi,\quad
\hat\l_\pi=\l_\pi+c\omega_0, \quad
\forall\l\in\ten^*.$$
We'll choose a new origin
$\pi$ as follows
$$\pi=c\g/\ell,\quad
\g=(-1,-1,\dots,-1,-2,\dots,-\ell),$$ 
where the integer $-p$ has
multiplicity $\nu_p$. 
The coweight $\check\g$ 
associated with $\g$ can be viewed as a cocharacter of the torus $T$.
So it yields the element $\check\g(z)\in T$ for each
$z\in\CC^\times$. Note that $\check\g(\eps)=g^{-1}$.
Note also that accordingly to (3.2) we have $\hat\l=\l+(c/\ell)\omega_0$.

\proclaim{4.2.~Proposition} (a) There is an algebra isomorphism
$\varkappa:\hat\gb^F_\kappa\to\hat\gb_\kappa$ yielding equivalences
of categories $\Cc_\kappa^F\to\Cc_\kappa$,
$\hat\Oc^F_{\dom,\kappa}\to\hat\Oc_{\nu,\kappa}$
such that $M(\hat\l)^F\mapsto M(\hat\l_\pi)_\nu$,
$L(\hat\l)^F\mapsto L(\hat\l_\pi)$.

(b) There is an algebra isomorphism
$\varkappa:\hat\gb^{F'}_\kappa\to\hat\gb_\kappa$ yielding equivalences
of categories $\Cc_\kappa^{F'}\to\Cc_\kappa$,
$\hat\Oc_{\dom,\kappa}^{F'}\to
\hat\Oc_{\nu,\kappa}$ such that $M(\hat\l)^{F'}\mapsto
M(\hat\mu_\pi)_\nu$, $L(\hat\l)^{F'}\mapsto L(\hat\mu_\pi)$ with
$\mu=-w_\nu(\l)$.

(c) If $M\in\Oc_\kappa^{F,\f}$ then $\dg D(\vk M)=\vk\, (\dg D M)$. If
$M\in\Cc_\kappa^F$ then $\vk M=\vk\,(\dg M)$.

\endproclaim

\noindent{\sl Proof :} $(a)$ We have
$$\gb^F=\bigoplus_{a\equiv b\,\mod\,\ell}\gen(b)\otimes t^a, 
\quad
\gen(b)=\{\xi\in\gen;[\check\g,\xi]=b\xi\}.\leqno(4.1)$$ 
Conjugating by $\check\g(-t)^{-1}$ takes $\gen(b)\otimes
t^a$ onto $\gen(b)\otimes t^{a-b}$. Composing this map with the
assignment $(-t)^{\ell}\mapsto -t$ yields a Lie algebra isomorphism
$\gb^F\to \gb$.
It lifts to a Lie algebra isomorphism
$${\gathered
\varkappa:\hat\gb^F\to\hat\gb,\quad\un\mapsto\un/\ell,\quad
\xi^{(a)}\mapsto
\check\g(-t)^{-1/\ell}(\xi^{(a/\ell)})-\delta_{a,0}\la\check\g:\xi\ra\un/\ell,
\endgathered}\leqno(4.2)
$$ see \cite{K1, thm.~8.5}.
The Lie
algebra $\varkappa(\hat\gb_{\dom}^F)$ is spanned by $\un$ and the
elements $\xi^{(a)}$ with $\xi\in\gen(b)$ and $a\ell+b\geqslant 0$.
Therefore we have
$$\gathered
\varkappa(\hat\gb_{\dom}^F)=\hat\qb_\nu,\vspace{2mm}
\varkappa(\xi+z\un)=\xi+(z-\la\check\g:\xi\ra)\un/\ell,\quad
\forall\xi\in\hen_\nu.
\endgathered$$
Part $(b)$ is the same, using the conjugation by $\check\g(t)$
instead of $\check\g(-t)^{-1}$. Indeed, we get the Lie algebra
isomorphism
$$\varkappa':\hat\gb^{F'}\to\hat\gb,\quad
\un\mapsto\un/\ell,\quad \xi^{(a)}\mapsto
\check\g(t)^{1/\ell}(\xi^{(a/\ell)})+\delta_{a,0}\la\check\g:\xi\ra\un/\ell.\leqno(4.3)$$
Set $\varkappa=\dag\circ\varkappa'$. Note that $\dg
L(\hat\l)=L(\widehat{-\l})$ where, in the right hand side, 
the highest weight is
relative to $\hat\bb'$. Note also that twisting by $\varkappa'$
takes $L(\hat\l)^{F'}$ to the simple module
$L(w_\nu(\hat\l_{-\pi}))$ in $\hat\Oc'_{\nu,\kappa}$.
Part $(c)$ is obvious. 

\qed

\vskip3mm

To unburden notation, we'll write $M\mapsto\vk M$ both for the
equivalence $\Cc_\kappa^F\to\Cc_\kappa$ and its
inverse, hoping it will not create any confusion.

\vskip3mm

\subhead 4.3.~Notation\endsubhead 
Let $S=A\cup\{n+1\}$.
Let $x_i$, $y_i$, $i\in S$, be as in Section 3.9. 
Recall that $x_{n+1}=y_{n+1}=\infty$.
Recall  $n$-tuples $(x_i;i\in A)$, $(y_i;i\in A)$ belong
to $\Alocnl$, $\Aloc_{n,1}$ respectively. View
$x_i$, $y_i$ as coordinates on $\Alocnl$, $\Aloc_{n,1}$.
Then the assignment $y_i\mapsto x_i^\ell$ yields an inclusion
$$\Rloc_{n,1}\subset\Rlocnl.$$
The group $D_\ell$ acts on $\PP^1_y$ by multiplication. We have
$$\CC[\PP^1_y/D_\ell]=\CC[\PP^1_y]^\flat.$$ 
Next, there is an obvious isomorphism
$\PP^1/D_\ell\to\PP^1$ which gives rise to an isomorphism
$\PP^1_y/D_\ell\to\PP^1-\{x_i^\ell;i\in S\}$. 
So we may identify
$$\CC[\PP^1_y]^\flat=\CC[z,(z-y_i)^{-1}].$$

\subhead 4.4.~Twisted affine coinvariants versus untwisted ones\endsubhead 
We want to compare twisted affine
coinvariants with untwisted ones. To do that we first generalize the
construction of the functor $\Cen$ in (3.6) by inserting an almost smooth
$\hat\gb_\kappa^ {F'}$-module at the point $x_0=0$. More precisely,
given $M\in\Cc^F_\kappa$ and $M'\in\Cc^{F'}_\kappa$
we consider the vector space 
$$\aligned
T(M',M)&=M'\otimes T(M)\vspace{2mm}
&=M'\otimes(\Vb^*)^{\otimes n}\otimes M.
\endaligned$$ 
Equip the
$\Rb_\nl$-module $T(M',M)_\Rlocnl$ with the unique $\Rlocnl$-linear
representation of the Lie algebra $\gen[\CC^\times]^F$ such that
$$\xi\otimes z^a\mapsto
(\xi^{(a)})_{(0)}+\sum_{i=1}^nx_i^a\xi_{(i)}+(\sharp\xi^{(a)})_{(n+1)}.$$
Then we set
$$\CIloc(M',M)_{n,\ell}=H_0\bigl(\gen[\CC^\times]^F,T(M',M)_\Rlocnl\bigr).$$
The $\gen[\CC^\times]^F$-action on $T(M',M)_\Rlocnl$ centralizes the
$\Bb_{n,\ell}$-action. We get a bifunctor
$$\hat\gb_\kappa^{F'}\text{-}\modb\times\hat\gb_\kappa^F\text{-}\modb\to
\Bb_{n,\ell}\text{-}\modb,\quad (M',M)\mapsto\CIloc(M',M)_{n,\ell}.$$
The $\Bb_\nl$-modules $\Cen(M)_\nl$ and $\Cen(M',M)_\nl$ are related
in the following way.

\proclaim{4.5.~Proposition}
There is a natural isomorphism of $\Bnl$-modules
$\Cen(M)_\nl\to\Cen(M(c\o_0)^{F'},M)_\nl$.
\endproclaim

\noindent{\sl Proof :} This is a consequence of Proposition
3.11$(b)$. To unburden the notation we give a proof for $n=1$ and we set
$R=\Rnl$. The case $n>1$ is the same. Set $$\aen=\gen[\CC^\times]^F,\quad
\ben=\gen[\CC]^F,\quad\aen'=(z^{-1}\gen[z^{-1}])^F.$$ Let
$\hat\ben=\ben\oplus\CC$, the trivial central extension of $\ben$
by $\CC$, and let $\hat\aen$ be the central extension of $\aen$ by
$\CC$ associated with the cocycle $(\xi\otimes f,\zeta\otimes
g)\mapsto\la\xi:\zeta\ra\Res_{z=0}(gdf)$. The assignments
$$\xi\otimes z^b\mapsto\xi^{(b)},\quad\xi\otimes z^b\mapsto
x_1^b\xi_{(1)}+(\sharp\xi^{(b)})_{(2)}$$ yield representations of
$\hat\aen_R$ on $M'_R$ and $\Vb^*_R\otimes_R M_R$ of level $c$
and $-c$ respectively. We have an isomorphism of $\hat\aen_R$-modules
$$M'_R\simeq\Ind^{\hat\aen_R}_{\hat\ben_R}(R).$$
Here $\ben_R$ acts trivially on $R$. The induction datum
$(\hat\aen_R,\hat\ben_R,\aen'_R,R)$ is $R$-split. Thus the tensor
identity yields an isomorphism of $\aen_R$-modules $$
T(M',M)_R\simeq\Ind^{\aen_R}_{\ben_R} (T(M)_R).$$ Thus we have
$$\CIloc(M',M)_\nl\simeq H_0\bigl(\ben,T(M)_R\bigr)=
\Cen(M)_\nl.$$

\qed

\vskip3mm

\noindent Next, assume that $M\in\Cc^F_\kappa$ and
$M'\in\Cc^{F'}_\kappa$. Set $N=\vk M$ and $N'=\vk M'$. The
$\hat\gb_\kappa$-modules $N$, $N'$ are almost smooth by Proposition 4.2.
They are canonically isomorphic to $M$, $M'$ respectively as vector
spaces. 
Let $\Bb_\nl$ act on $(\Vb^*)^{\otimes n}_{\Rlocnl}$ as in (3.4). 
Setting $\ell=1$ we get a $\Bb_{n,1}$-action on $(\Vb^*)^{\otimes n}_{\Rloc_{n,1}}$. 
The inclusion $\Rb_{n,1}\subset\Rb_{n,\ell}$ in Section 4.3 yields
an embedding $W\ltimes\Rb_{n,1}\subset\Bb_{n,\ell}$, where
the $\eps_i$'s act trivially on $\Rb_{n,1}$ in the left hand side.
Thus, by induction, the algebra $\Bb_{n,\ell}$ acts on
$$(\Vb^*)^{\otimes n}_{\Rloc_{n,1}}\otimes_{\Rloc_{n,1}}\Rlocnl.$$
Now, consider the $\Rb_\nl$-linear map
$$\gathered
(\Vb^*)^{\otimes n}_{\Rlocnl}\to
(\Vb^*)^{\otimes n}_{\Rloc_{n,1}}\otimes_{\Rloc_{n,1}}\Rlocnl,
\vspace{2mm}
v=v_1\otimes v_2\otimes\cdots\otimes v_n\mapsto 
\check\g(x_1)v_1\otimes \check\g(x_2)v_2\otimes\cdots\otimes 
\check\g(x_n)v_n,
\vspace{2mm}
\check\g(x_i)=(x_i^{-1},x_i^{-1},\dots,x_i^{-1},x_i^{-2},\dots,x_i^{-\ell}).
\endgathered$$
It is indeed a $\Bb_{n,\ell}$-linear map (the $\Sen_n$-linearity is obvious, and to prove
the $W$-linearity it is enough to observe that
$\check\g(\eps)=g^{-1}$, see Section 4.1, and that $\eps_iv=g_{(i)}v$, see Section 3.1).
This yields also a
$\Bb_{n,\ell}$-linear map $$\varkappa: T(M',M)_\Rlocnl\to T(\dg
N',N)_{\Rloc_{n,1}}\otimes_{\Rloc_{n,1}}\Rlocnl.
\leqno(4.4)$$

\proclaim{4.6.~Proposition} Let $M\in\Cc^F_\kappa$,
$M'\in\Cc^{F'}_\kappa$, $N=\vk M$ and $N'=\vk M'$. The map
$\varkappa$ gives a $\Bb_{n,\ell}$-module isomorphism
$\CIloc(M',M)_{n,\ell}\to \CIloc(\dg
N',N)_{n,1}\otimes_{\Rloc_{n,1}}\Rlocnl.$
\endproclaim

\noindent{\sl Proof :} We must prove that the $\Bb_\nl$-linear map
$\varkappa$ in (4.4) factors to an isomorphism of
$\Bb_{n,\ell}$-modules
$$H_0(\gen[\CC^\times]^F,T(M',M)_\Rlocnl)
\to H_0(\gen[\CC^\times],T(\dg N',N)_{\Rloc_{n,1}})\otimes_{\Rloc_{n,1}}\Rlocnl.$$ The vector space
$\gen[\CC^\times]^F$ is spanned by the elements
$$\xi\otimes z^a\ \roman{with}\ \xi\in\gen(b),\ a\equiv -b\ \mod\ \ell.$$
Here $\gen(b)$ is as in (4.1). The conjugation by $\check\g(z)$
takes $\xi\otimes z^a$ to 
$$\check\g(z)(\xi\otimes z^{a})=\xi\otimes z^{a+b}.$$
Composing it with
the linear map $\xi\otimes z^{a+b}\mapsto\xi\otimes z^{(a+b)/\ell}$
yields a Lie algebra isomorphism
$$\gen[\CC^\times]^F\to\gen[\CC^\times].\leqno(4.5)$$ 
We'll regard $\Rloc_{n,1}$ as the subalgebra
of $\Rlocnl$ generated by the elements $y_i=x_i^\ell$.
Under the map (4.5) the element $\xi\otimes z^a\in\gen[\CC^\times]^F$ 
above acts on $T(\dg N',N)_{\Rloc_{n,1}}$ as 
$$\gathered
\varkappa'(\xi^{(a)})_{(0)}+\sum_{i=1}^ny_i^{(a+b)/\ell}\xi_{(i)}+
\varkappa(\sharp\xi^{(a)})_{(n+1)}=\vspace{2mm}
\varkappa'(\xi^{(a)})_{(0)}+\sum_{i=1}^n
x_i^a(\check\g(x_i)\xi)_{(i)}+
\varkappa(\sharp\xi^{(a)})_{(n+1)}.
\endgathered$$ 
Here the maps
$\varkappa':\hat\gb^{F'}\to\hat\gb$ and
$\varkappa:\hat\gb^{F}\to\hat\gb$ are as in (4.2), (4.3). 
Note that
$$T(\dg N',N)_{\Rlocnl}=T(\dg N',N)_{\Rloc_{n,1}}
\otimes_{\Rloc_{n,1}}\Rlocnl.$$
Under
the canonical isomorphisms $$\dg N'={}^{\varkappa'}M'\simeq M',\quad
N=\vk M\simeq M$$ the actions of the elements
$\varkappa'(\xi^{(a)})$, $\varkappa(\sharp\xi^{(a)})$ on $\dg N'$,
$N$ are the same as the actions of the elements $\xi^{(a)}$,
$\sharp\xi^{(a)}$ on $M'$, $M$ respectively. Therefore, we have
$$\xi\otimes z^a\cdot\varkappa(v)=\varkappa((\xi\otimes
z^a)\cdot v),\quad\forall v\in T(M',M)_\Rlocnl.$$

\qed

\vskip2cm

\head 5. Complements on the category $\hat\Oc_{\nu,\kappa}$\endhead

This section is a reminder on the structure of the category $\Oc$
for (untwisted) affine Lie algebras. This section does not contain new
results. Proposition 5.8 is proved in the appendix. It is standard,
but we have not found any proof in the literature.
Recall that we'll always assume that
$\kappa\in\CC^\times$.


\subhead 5.1.~Complements on the affine Weyl group\endsubhead 
When it creates no confusion we'll abbreviate $\Sen=\Sen_m$.
Let $\hat\Sen$ be the affine Weyl group of
$\Sen$, i.e., the semi-direct product $\Sen\ltimes\ZZ\Pi$. For any real affine root $a$ let
$s_\a\in\hat\Sen$ be the reflection associated with $\a$.There
is a unique linear representation of $\hat\Sen$ on $\tb^*$ such
that $w\in \Sen$, $\tau\in\ZZ\Pi$ act in the following way :
$$w(\epsilon_i)=\epsilon_{v(i)}, \quad
w(\o_0)=\o_0, \quad w(\delta)=\delta, \quad \tau(\delta)=\delta,$$
$$\tau(\epsilon_i)=\epsilon_i-\la\tau:\epsilon_i\ra\delta, \quad
\tau(\o_0)=\tau+\omega_0-\la\tau:\tau\ra\delta/2.$$ 
We abbreviate
$$\rhoaf=\rho+m\o_0.$$
The dot-action of $\hat\Sen$ on
$\tb^*$ is given by
$$w\bullet\l=w(\l+\rhoaf)-\rhoaf.$$
Both actions factor to
representations of $\hat\Sen$ on the vector space $\tb'$.  Recall
the notation $$\hat\l=\l+c\o_0,\quad
\tilde\l=\hat\l+z_\l\delta,\quad z_\l=-\la\l:2\rho+\l\ra/2\kappa,\quad\forall\l\in\ten^*.
\leqno(5.1)$$ We have $\hat\l=w\bullet\hat\mu$ iff
$\tilde\l=w\bullet\tilde\mu$.

For each $\l\in\tb^*$ we set
$$\aligned
\hat\Pi(\l)=\{\a\in\hat\Pi;
2\la\l+\hat\rho:\a\ra\in\ZZ\la\a:\a\ra\}.\endaligned$$
Let $\tb^*_0=\{\l\in\tb^*\,;\,\la\l+\rhoaf:\delta\ra\neq 0\}$.
Note that if $\hat\l$ is as in (5.1) then it lies in $\tb^*_0$ iff
$\kappa\neq 0$.
For each $\l\in\tb^*_0$ we have
$$\aligned
\hat\Pi(\l)=\{\a\in\hat\Pi_\re;
\la\l+\hat\rho:\a\ra\in\ZZ\},\endaligned$$
a root system with the set of positive roots
$\Pib(\lambda)^+=\Pib^+\cap\Pib(\l).$
Let
$\hat\Sen(\l)$ be its Weyl group. 
We call $\hat\Sen(\l)$ the {\it integral Weyl group} associated with $\l$.

For each $\l\in\tb^*_0$ we consider also the group
$$\hat\Sen_{\l}=\{w\in\hat\Sen\,;\,w\bullet\l=\l\}.$$
It is a (finite) subgroup isomorphic to the Weyl group of the root system
$$\hat\Pi_\l=\{\a\in\Pib_\re;\la\l+\rhoaf:\a\ra=0\}.$$
See \cite{KT, sec.~2} for details.

Finally, let $\Pi_\nu$, $\Pi^+_\nu$, $\Sen_\nu$ be the root system,
the set of positive roots and the Weyl group of $\hen_\nu$. Note that
$\Pi_\nu\subset\Pi$ and $\Pi^+_\nu=\Pi^+\cap\Pi_\nu$.

\subhead 5.2.~The partial order on $\Delta_{\hat\Oc_\nub}$\endsubhead
We'll say that an affine weight $\l\in\tb^*$ is {\it $\nu$-regular} (resp.~{\it $\nu$-integral})
if we have $\la\l:\a\ra\neq 0$ (resp.~$\la\l:\a\ra\in\ZZ$) for each $\a\in\Pi_\nu$. Note that if $\l+\hat\rho$ is $\nu$-regular and integral then there exists a unique element $w\in\Sen_\nu$ such that
$$\l_+:=w\bullet\l$$
is $\nu$-dominant, and we set $$sn(\l)=(-1)^{l(w)}.$$
We equip the set $\tb^*$ with the partial order given by
$$\l\geqslant\mu\iff\l-\mu\in\NN\hat\Pi^+.$$
Now, we define the following partial orders on the set of $\nu$-dominant affine weights.

\vskip1mm
\itemitem{$(a)$}
Let $\preccurlyeq$ be the transitive closure of the binary relation
such that $\l\Uparrow\mu$ iff the simple
$\tilde\gb$-module $L(\l)$ is a
Jordan-H\"older factor of the parabolic Verma module $M(\mu)_\nu$.

\vskip1mm
\itemitem{$(b)$}
Let $\trianglelefteq$ be the transitive and reflexive closure of the
binary relation such that
$\l\uparrow\mu$ iff there is a root $\a$ in
$\hat\Pi^+_\re\setminus\Pi_\nu^+$ 
such that $\la\mu+\rhoaf:\a\ra\in\ZZ_{>0}$, $
s_\a\bullet\mu+\rhoaf\roman{\ is\ \nu\text{-}regular}, $ and we have
$$\gathered
\l=(s_\a\bullet\mu)_+<\mu.
\endgathered$$

\proclaim{5.3.~Proposition} (a) The partial order $\trianglelefteq$
refines the partial order $\preccurlyeq$.

(b) If $\Hom_{\hat\gb}(M(\hat\l_1)_\nu,M(\hat\l_2)_\nu)\neq 0$ then $\hat\l_1\trianglelefteq\hat\l_2$.
\endproclaim

\noindent{\sl Proof :}
Part $(b)$ follows from $(a)$, 
because if $\phi\in\Hom_{\hat\gb}(M(\hat\l_1)_\nu,M(\hat\l_2)_\nu)$
is non zero then $\phi(M(\hat\l_1)_\nu)$ is a submodule of
$M(\hat\l_2)_\nu$ whose top contains $L(\hat\l_1)$. Hence $L(\hat\l_1)$ is a Jordan-H\"older factor of
$M(\hat\l_2)_\nu$.
Now, we prove $(a)$. We must prove that if $\mu\preccurlyeq\l$ then $\mu\trianglelefteq\l$.
Given a $\nu$-dominant affine weight $\l\in\tb^*$, the Kac-Kazhdan formula for the
Shapovalov determinant of the contravariant form on $M(\l)_\nu$ restricted to its weight
$\mu$-subspace is given, up to a non-zero scalar, by the expression
$$\gathered
\det(\l)_{\nu,\mu}=\prod_{n>0}\prod_{\a\in\hat\Pi^+\setminus\Pi_\nu^+}\bigl(2\la\l+\hat\rho:\a\ra-n\la\a:\a\ra\bigr)^{\chi(\l-n\a)_{\nu,\mu}},\vspace{2mm}
\chi(\l)_{\nu,\mu}=\sum_{w\in\Sen_\nu}(-1)^{l(w)}\dim\,M(w\bullet\l)_{\nu,\mu},
\endgathered$$
see \cite{KK, p.~107}.
Let us consider the sets
$$
\aligned
\Sc_\l&=\{\a\in\hat\Pi^+\setminus\Pi_\nu^+\,;\,\exists n\in\ZZ_{>0},\,2\la\l+\hat\rho:\a\ra=n\la\a:\a\ra\},
\vspace{2mm}
&=\{\a\in\hat\Pi^+_\re\setminus\Pi_\nu^+\,;\,\la\l+\hat\rho:\a\ra\in\ZZ_{>0}\},
\vspace{2mm}
\Sc^0_\l&=\{\a\in\Sc_\l\,;\,s_\a\bullet\l+\hat\rho\ \roman {is\ }\nu\text {-regular}\}.
\endaligned$$
Now, suppose that $L(\mu)$ is a subquotient of $M(\l)_\nu$ with $\l,\mu$ both $\nu$-dominant affine weights
and $\mu<\l$. Then $L(\mu)$ must be a subquotient of the maximal submodule $M(\l)^1_\nu$ of $M(\l)_\nu$.
The Jantzen filtration yields a decreasing sequence of submodules 
$M(\l)^k_\nu$, $k>0,$ of $M(\l)_\nu$ such that
$$\sum_{k>0}\ch(M(\l)^k_\nu)=
\sum_{\a\in\Sc_\l^0}sn(s_\a\bullet\l)\,\ch(M((s_\a\bullet\l)_+)_\nu),$$
where the symbol $\ch$ denotes the formal character, 
see e.g., \cite{J, sec.~4.1}.
Thus there is a real affine root $\a\in\Sc_\l^0$ such that $L(\mu)$ is a subquotient of
$M((s_\a\bullet\l)_+)_\nu$. Since $(s_\a\bullet\l)_+\uparrow\l$, an obvious induction implies the proposition.
Compare the discussion
in \cite{JLT, sec.~9.4-10.6} for instance.

\qed

\vskip3mm

If $\tilde\l,\tilde\mu$ are as in (5.1)
we may write $\mu\trianglelefteq\l$  for
$\tilde\mu\trianglelefteq\tilde\l$.
To define the order relation we embed
$\hat\Oc_\nub$ in $\tilde\Oc_\nub$ and 
we equip the set
$\Delta_{\hat\Oc_\nub}$ with the order $\trianglelefteq$ of the highest
weights of the parabolic Verma modules.

\subhead 5.4.~Remark\endsubhead 
Note that if  $\l\trianglelefteq\mu$ then $\mu-\l\in\NN\Pi^+$ and
$\l\in\hat\Sen\bullet\mu$.

\vskip3mm

\subhead 5.5.~The truncated category\endsubhead
Recall that $\ub_\nu$ denotes the pronilpotent radical of $\hat\qb_\nu$.
Let $\Pib^\nu\subset\Pib^+$ be the set of roots of $\ub_\nu$.
Set $\zb=\CC\partial\oplus\zen\oplus\CC\un$ where $\zen\subset\hen_\nu$
is the central Lie subalgebra such that
$\hen_\nu=[\hen_\nu,\hen_\nu]\oplus\zen$.
Note that $\zb\subset\tb$.
Let $z:\tb^*\to\zb^*$ be
the restriction of linear forms.
We equip the set $\zb^*$ with the
partial order such that 
$$z_2\leqslant z_1\iff z_1-z_2\in z(\NN\Pib^\nu).$$
Given a finite subset $B\subset\zb^*$ we put
$$\lub\Lambda=\{\l\in\tb^*;\exists\b\in B,\,
z(\l)\leqslant\b\}.$$
Let $\lub\hat\Oc_{\nub}\subset\hat\Oc_{\nub}$ be the Serre
subcategory consisting of the modules whose $\l$-weight space
vanishes if $\l\notin\lub\Lambda$. Note that any object of $\hat\Oc_{\nub}^{\f}$ lies in
$\lub\hat\Oc_{\nub}$ for some $B$. The following is well-known.

\proclaim{5.6.~Proposition}
The category $\lub\hat\Oc_{\nub}^\f$ is a highest weight category. The
poset of standard modules is the set of parabolic Verma modules
which belong to $\lub\hat\Oc_{\nub}^\f$ with the order relation
$\trianglelefteq$. 
\endproclaim

\noindent{\sl Proof :}
By Proposition 2.9$(c)$ the category $\hat\Oc^\f_{\nub}$ is Abelian, any object has a finite length,
and Hom sets are finite dimensional. Thus te category $\lub\hat\Oc_{\nub}^\f$ is Abelian and Artinian.
The axioms $(a)$ and $(c)$ of quasi-hereditary categories, see Section 0.1, are obvious. 
The axiom $(b)$ follows from Proposition 5.3$(b)$. 
Thus it is enough to check that any parabolic Verma module $M$ 
has a projective cover in $\lub\hat\Oc_{\nub}^\f$ whose
kernel is standardly filtered with subquotients $>M$. 
This is well-known, and proved in \cite{RW, cor.~10, cor.~13, thm.~4}. 
See also \cite{FKM}, \cite{S} for a more recent exposition. 

\qed

\vskip3mm

\subhead 5.7.~Kazhdan-Lusztig polynomials\endsubhead 
Set 
$$\alc=\{\l\in\tb^*_0\,;\,\la\l+\rhoaf:\a\ra\leqslant0,\,\forall\a\in\hat\Pi(\l)^+\}.$$ 
Recall that for each $\l\in\tb^*_0$ such that $\hat\Pi(\l)\neq\emptyset$
and $\la\l+\rhoaf:\delta\ra\notin\QQ_\ddom$ the
set $(\hat\Sen(\l)\bullet\l)\cap\alc$ consists exactly of one element, see \cite{KT, lem.~2.10}.
Fix $\l\in\tb^*_0$.  Let
$\hat\Sen^{\l}\subset\hat\Sen$ be the set of minimal length
representatives in the left cosets relative to $\hat\Sen_{\l}$.  Since
$\hat\Sen_\l$ is also a parabolic subgroup of the integral Weyl group
$\hat\Sen(\l)$, the set $\hat\Sen(\l)^{\l}$ is well-defined. We'll use the
symbol $P^{\l,-1}_{v,w}$ for Deodhar's parabolic Kazhdan-Lusztig
polynomials of type $\hat\Sen(\l)/\hat\Sen_{\l}$, see \cite{D}. 

\proclaim{5.8.~Proposition} Let $\l\in\alc$ be
such that $\hat\Pi(\l)\neq\emptyset$. If $w\in\hat\Sen(\l)^{\l}$ is such that $w\bullet\l$
is $\nu$-dominant, then we have the following formula in
$[\hat\Oc_{\nub}^\f]$
$$[L(w\bullet\l)]=
\sum_{v\leqslant w}(-1)^{l(w)-l(v)}P^{\l,-1}_{v,w}(1)
[M(v\bullet\l)_\nu].$$ The sum is over all $v\in \hat\Sen(\l)^{\l}$ such
that $v\bullet\l$ is $\nu$-dominant.
The symbol $\leqslant$ denotes the
Bruhat order on $\hat\Sen(\l)$.
\endproclaim

\subhead 5.9.~Remark\endsubhead Note that if $\l\in\alc$ is
such that $\hat\Pi(\l)\neq\emptyset$ then we have
$\kappa:=\la\l+\rhoaf:\delta\ra\notin\QQ_\dom$, 
see e.g., \cite{KT, lem.~2.10}.

\vskip2cm

\head 6.~ Definition of the functor $\Een$\endhead

We can now construct our main functor $\Een$. It takes a module from
the affine parabolic category $\Oc$ to a module in $\Hc_{h,H}$. We'll
prove that the functor $\Een$ preserves the posets of standard
modules. In this section we'll make the following assumption
$$\gathered
\nu\in\Cc_{m,\ell},\quad 
\kappa\notin\QQ_{\sss\geqslant 0},\quad
h=1/\kappa,\vspace{2mm}
h_p=\nu_p^\bullet/\kappa-m/\ell\kappa,\quad
\forall p\in\Lambda.
\endgathered$$

\subhead 6.1.~Notation and definition of $\Een$\endsubhead
Recall that by Proposition 3.6$(b)$ we have a functor
$$\Cen:\Cc^F_\kappa\to\Hb_{h,H}\text{-}\modb.$$ Composing it with the functor $\varkappa$ in
Proposition 4.2$(a)$ we get the functor
$$\VVkap:\Cc_\kappa\to\Hb_{h,H}\text{-}\modb,\qquad
M\mapsto H_0\bigl(\gen[\CC]^F,T(\vk M)_{\Rb}\bigr).\leqno(6.1)$$
Recall that
$$\pi=c(-1,-1,\dots,-1,-2,\dots,-\ell)/\ell,$$ 
where the integer $-p$ has
multiplicity $\nu_p$. 
Set
$$\l_\pi=\l+\pi,\quad
\hat\l_\pi=\l_\pi+c\omega_0, \quad
\tilde\l_\pi=\hat\l_\pi+z_{\l_\pi}\delta,\quad
\forall\l\in\ten^*.$$
If $\l\in\ten^*$ is a $\nu$-dominant weight we'll abbreviate
$$\Delta_{\l,\nub}=M(\hat\l_\pi)_\nu,\quad
S_{\l,\nub}=L(\hat\l_\pi).$$
Given a finite subset $B\subset\zb^*$ let
$\lub P_{\l,\nub}$ denote the projective cover of
$S_{\l,\nub}$ in $\lub\hat \Oc_{\nub}$.

\vskip3mm

\subhead 6.2.~The functor $\Een$ and the standard modules\endsubhead
First, let us compute the image by $\Een$ of the standard modules.

\proclaim{6.3.~Proposition} 
(a) The functor $\Een$ is right exact and takes $\hat\Oc_{\nub}$ to $\Hc_{h,H}$.

(b) We have $\Een(\Delta_{\l,\nub})=\Delta_{\l^\circ,h,H}$ if
$\l\in\Pc_{n,\nu}^\ell$, and $0$ else.

(c) The functor $\Een$ takes 
$\hat\Oc_{\nub}^\f$ to $\Hc_{h,H}^\f$.

\endproclaim

\noindent{\sl Proof :} 
First, we prove $(a)$.
It is enough to check
that $T(M)_\Rb$ is a locally nilpotent $\Rb^{*}$-module for
$M\in\hat\Oc_{\nub}$. We'll identify $M$ with the $\hat\gb^F$-module
$\vk M$. By Proposition 3.6$(b)$ the operator $\bar
y_i-\g^F_{i,n+1}$ vanishes on the vector space $T(M).$ Fix a finite
dimensional $\hat\gb^F_\dom$-submodule $E\subset M$. Formula (3.8)
implies that
$$(\bar y_i-\g_{i,n+1}^F)(T(E)_{\Rb_{\sss\leqslant a+1}})\subset
T(E)_{\Rb_{\sss\leqslant a}},$$ where $\Rb_{\sss\leqslant a}\subset \Rb$
is the subspace of the polynomials of degree $\leqslant a$. We have
also
$$\g^F_{i,n+1}(T(E)_{\Rb})\subset
T(\gb_\domm^F E)_{\Rb}.\leqno(6.2)$$ Thus, since the action of
$\gb_\domm^F $ on $E$ is nilpotent and $E$ is finite
dimensional, there is an integer $d>0$ such that
$(\g^F_{i,n+1})^{d}(T(E)_{\Rb})=0.$ Therefore, if $b$ is large enough
then $\bar y_i^b(T(E)_{\Rb_{\sss\leqslant a}})=0$. The right
exactness of $\VVkap$ is obvious, because taking coinvariants is a
right exact functor.

Now, we prove $(b)$. 
Given a $\hat\gb^F_\dom$-module $M$ we consider the induced $\hat\gb^F$-module
$M^F_\kappa$. We have an isomorphism of $\Bb$-modules
$\Cen(M^F_\kappa)=\Ind_W^\Bb(\Xen(M))$ by Proposition 3.8$(b)$. The $y_i$-action on an
element $v\in \Xen(M)$ is given by $y_iv=\g_{i,n+1}^Fv$ by Proposition 3.6$(b)$.
 If $\gb_\domm^F$ annihilates $M$ then (6.2) yields $y_iv=0$. Thus
the subalgebra $\Bb^*\subset\Hb_{h,H}$ acts trivially on the subspace
$\Xen(M)\subset\Cen(M^F_\kappa)$. This yields a $\Bb^*$-module isomorphism
$$\Cen(M^F_\kappa)\simeq\Ind_{\Bb^*}^{\Hb_{h,H}}(\Xen(M)).$$
Now, Proposition 4.2$(a)$ yields
$$\VVkap(\Delta_{\l,\nub})=
\Een(M(\hat\l_\pi)_\nu)=\Cen(L(\hen_\nu,\l)_\kappa^F).$$
We are done, because $\Xen(L(\hen_\nu,\l))=\Xen_{\l^\circ}$
by Proposition 3.8$(a)$.

Finally, we prove $(c)$. Given $M\in\hat\Oc_{\nub}^\f$ we must check that
$\Een(M)$ has a finite length. This follows from an easy induction on the length
of $M$. First, if $M$ is simple then it is a quotient of a module $\Delta_{\l,\nu,\kappa}$.
Thus, since $\Een$ is right exact, $\Een(M)$ is a quotient of $\Delta_{\l^\circ,h,H}$.
Thus $\Een(M)$ lies in $\Hc_{h,H}^\f$. Next, there is an exact sequence
$$\Delta_{\l,\nu,\kappa}\to M\to N\to 0$$
such that the length of $N$ is strictly less than the length of $M$. Since $\Een$ is right exact this yields
the exact sequence
$$\Delta_{\l^\circ,h,H}\to \Een(M)\to \Een(N)\to 0.$$
Thus $\Een(M)$ has a finite length by the induction hypothesis.
 
\qed

\vskip3mm

Now, we can compare the partial orders on the set of standard modules in
$\hat\Oc_{\nu,\kappa}$ and in $\Hc_{h,H}$.

\proclaim{6.4.~ Proposition} If $\l,\mu\in\Pc_{n,\nu}^\ell$ and
$\Delta_{\mu,\nub}\trianglelefteq\Delta_{\l,\nub}$ then
$\Delta_{\mu^\circ,h,H}\preccurlyeq\Delta_{\l^\circ,h,H}$.
\endproclaim

\noindent{\sl Proof :} Recall the decomposition
$J=J_1\sqcup J_2\sqcup\dots\sqcup J_\ell$ in Section 1.4.
Let $\g,\tau$ be the weights given by 
$$\gathered
\g_j=h_1+h_2+\cdots h_{\ell-p},\quad
\tau_j=m+\nu_\ell,
\quad\forall j\in J_p.
\endgathered$$
For each weights $\l,\mu$ we set
$$a(\l,\mu)=\la\l:\l+2\tau+2\rho\ra/2\kappa-\la\mu:\mu+2\tau+2\rho\ra/2\kappa.$$
A computation yields
$$\aligned
\la\tilde\l_\pi-\tilde\mu_\pi:\pi+c\omega_0\ra
&=cz_{\l_\pi}-cz_{\mu_\pi}+\la\l-\mu:\pi\ra\cr
&=m\la\l-\mu:\pi\ra/\kappa-c\la\l:\l+
2\rho\ra/2\kappa+c\la\mu:\mu+2\rho\ra/2\kappa\cr
&=\la\l-\mu:c\tau+m\pi\ra/\kappa-ca(\l,\mu).
\endaligned$$
On the other hand, we have
$$\aligned
\kappa a(\l,\mu)=\sum_{p}
\bigl(n({}^t\!\l_p)-n({}^t\!\mu_p)-n(\l_p)+n(\mu_p)\bigr).
\endaligned$$
So we get
$$\aligned
\theta_{\mu^\circ}-\theta_{\l^\circ}&=-\ell\sum_{p}(h_1+\cdots
+h_{p-1})(|\l_p^\circ|-|\mu_p^\circ|)-\ell a(\l,\mu),\cr
&=-\ell\la\l-\mu:\g\ra-\ell a(\l,\mu).
\endaligned$$
We claim that we have
$$c(\theta_{\mu^\circ}-\theta_{\l^\circ})/\ell=
\la\tilde\l_\pi-\tilde\mu_\pi:\pi+c\omega_0\ra.\leqno(6.3)$$
It is enough to check that $c\kappa\g+c\tau+m\pi=0$. Since
this tuple has the same entries on each segment $J_p$,
it is enough to prove that its
$i_p$-th entry is zero. The latter is
$$\aligned
&
\sum_{r=1}^{\ell-p}
c\kappa h_r+c(\nu_\ell+\nu_1+\nu_2+\cdots\nu_{p-1})-cmp/\ell=\cr
&
\sum_{r=1}^{\ell-p}
c(\kappa h_r-\nu^\bullet_r)+
cm(\ell-p)/\ell=
\sum_{r=1}^{\ell-p}
c(\kappa h_r-\nu^\bullet_r+m/\ell)=0.
\endaligned$$
Now we prove the proposition using formula (6.3).
Assume that $\mu_\pi\triangleleft\l_\pi$. 
By definition of the
partial order $\trianglelefteq$, a simple induction allows us to assume that
$$\tilde\mu_\pi=(s_\a\bullet\tilde\l_\pi)_+<\tilde\l_\pi,\quad
\a\in\Pib^+_\re\setminus\Pi^+_\nu.$$
For each simple affine root $\a_i$ we have
$\la \a_i:\pi+c\o_0\ra=c/\ell$ or 0. Therefore
$$
\gathered
\a\in\Pib^+_\re\setminus\Pi_\nu^+\Rightarrow \la\a:\pi+c\o_0\ra
\in\ZZ_\domm\, c/\ell.
\endgathered
$$
This implies that
$$\la\tilde\l_\pi-\tilde\mu_\pi:\pi+c\o_0\ra=
\la\tilde\l_\pi-s_\a\bullet\tilde\l_\pi:\pi+c\o_0\ra\in\ZZ_\ddom\,c/\ell.$$
Using (6.3) this yields
$$\theta_{\mu^\circ}-\theta_{\l^\circ}\in\ZZ_\ddom.$$
Thus, the definition of the order in (1.4) implies that
$$\Delta_{\mu^\circ,h,H}\prec\Delta_{\l^\circ,h,H}.$$

\qed

\vskip3mm

\subhead 6.5.~Remarks\endsubhead $(a)$ We have $\Een(M)\in\Hc_{h,H}$
for any $M\in\hat\Oc_\kappa$, but $\Een(M)=0$ if $M$ is a simple
module which does not belong to $\hat\Oc_{\nub}$.

\vskip1mm

$(b)$ The Lie algebra $\hat\gb^F$ is $\ZZ$-graded by letting
$\xi^{(a)}$ be of degree $-a$. We can consider the category of
$\ZZ$-graded modules which belong to
$\hat\Oc^F_{\dom,\kappa}$.
The algebra
$\Hb_{h,H}$ is $\ZZ$-graded as in Section 1.6 and we can also consider
the category of $\ZZ$-graded modules which belong to $\Hc_{h,H}$.
The functor $\Cen$ 
lifts to a functor between these categories of graded
modules.

\vskip3mm

\subhead 6.6.~Examples\endsubhead 
Write $1_p$ for the $\ell$-partition whose
$p$-th partition is $(1)$ and all other are zero.

\vskip1mm

$(a)$ Set $\kappa=-2$, $n=1$, $m=10$, $\ell=4$, $\nu=(2,1,6,1)$. We
have $\Delta_{1_3^\circ,h,H}\succ\Delta_{1_2^\circ,h,H}\succ
\Delta_{1_1^\circ,h,H}\succ\Delta_{1_4^\circ,h,H}$ and
$\Delta_{1_3,\nub}\triangleright\Delta_{1_2,\nub}\ntriangleright
\Delta_{1_1,\nub}\triangleright\Delta_{1_4,\nub}$.
Thus the implication in Proposition 6.4 is not an equivalence.



\vskip1mm

$(b)$ 
A direct computation shows that the module $\Delta_{0,\nub}$ may not
be simple (f.i., set $\kappa=-1$, $m=7$, $\ell=4$, and $\nu=(1,1,4,1)$).

\vskip2cm

\head 7. The functor $\Een$ is exact on standardly filtered
modules\endhead

If $\ell=1$ the functor $\Een$ is an equivalence of quasi-hereditary
categories. Since there is no proof in the literature we have given
one in Section A.5. For an arbitrary positive integer $\ell$ this is not true
anymore. However we expect $\Een$ to be an important tool to prove
the dimension conjecture. We'll prove that $\Een$ is exact on
standardly filtered modules. We conjecture that it preserves the set
of indecomposable projective modules. In this section we'll assume
once again that 
$$\gathered
\nu\in\Cc_{m,\ell},\quad 
\kappa\notin\QQ_{\sss\geqslant 0},\quad
h=1/\kappa,\vspace{2mm}
h_p=\nu_p^\bullet/\kappa-m/\ell\kappa,\quad
\forall p\in\Lambda.
\endgathered$$

\subhead 7.1.~Reminder on the Kazhdan-Lusztig tensor product\endsubhead 
A {\it monoidal category} is a tuple $(\Ac,\otimes,a,\un)$ 
consisting of a category
$\Ac$, a functor $\otimes:\Ac\times\Ac\to\Ac$, a natural isomorphism
$$
a_{L,M,N}:(L\otimes M)\otimes N\to L\otimes(M\otimes N),\quad
L,M,N\in\Ac,$$ and a unit object $\un$
satisfying the triangle and pentagon axioms, see e.g., 
\cite{O, sec.~2.2}.
The isomorphism $a$ is called the
{\it associativity isomorphism}. 
A category $\Mc$ is a {\it left module category} over  the monoidal category $(\Ac,\otimes,a,\un)$ iff there exists a functor $\otimes: \Ac\times\Mc\to\Mc$
together with natural associativity and unit isomorphisms 
$$(M\otimes N)\otimes X\to M\otimes(N\otimes X),\quad \un\otimes X\to X,\quad X\in\Mc,\quad M,N\in\Ac,$$ which satisfy appropriate pentagon and triangle axioms. In other words, a left module category is the same as the datum of a monoidal functor
$\Ac\to\roman{Fun}(\Mc,\Mc)$ to the monoidal category of endofunctors of $\Mc$, see \cite{O, prop.~2.2}. 
Similarly one defines the structure of a {\it right module category}. 
Finally, a category $\Mc$ is a {\it bimodule category} over  $(\Ac,\otimes,a,\un)$ iff $\Mc$ has
both left and right $(\Ac,\otimes,a,\un)$-module category structures and
a natural family of isomorphisms
$$(X\otimes M)\otimes Y \to X\otimes(M\otimes Y),\quad X,Y\in\Ac,\quad M\in\Mc$$
satisfying the obvious pentagon axioms, see e.g., \cite{Gr, prop.~2.10}.

Now, recall that the Kazhdan-Lusztig tensor product 
$$\dot\otimes:\hat\Oc_{\dom,\kappa}^\f\times\hat\Oc_{\dom,\kappa}^\f\to
\hat\Oc_{\dom,\kappa}^\f$$ 
equips $\hat\Oc_{\dom,\kappa}^\f$ with the structure of a monoidal category
$$(\hat\Oc^\f_{\dom,\kappa},\dot\otimes,a,M(c\o_0)),$$
see \cite{KL, sec.~14.6, 18.2, 31}. 
The space of affine coinvariants 
is defined as in Remark A.2.2.
It depends on the choice of a point $x$ in the set $\Cc$ defined in Remark A.2.2.
When no confusion is possible we'll abbreviate
$$\la M_i;i\in S\ra=\la M_i;i\in S\ra_x.$$ 
The following is folklore. See Section A.2 for details.
Note that we'll note use Corollary 7.3 in this paper. It is given here for the sake of completeness.

\proclaim{7.2.~Proposition} (a) There are right biexact bifunctors
$\dot\otimes : \hat\Oc_{\dom,\kappa}^\f\times\hat\Oc_{\nub}^\f\to
\hat\Oc_{\nub}^\f$ and $\dot\otimes :
\hat\Oc_{\nub}^\f\times\hat\Oc_{\dom,\kappa}^\f\to \hat\Oc_{\nub}^\f$ yielding the structure
of a bimodule category over $(\hat\Oc^\f_{\dom,\kappa},\dot\otimes,a,M(c\o_0))$
on $\hat\Oc_{\nub}^\f$.

(b) For $M_1,\dots M_n\in\hat\Oc_{\dom,\kappa}^\f$ and
$M_0,M_{n+1}\in\hat\Oc_{\nub}^\f$ there are natural isomorphisms of
finite dimensional $\CC$-vector spaces
$$\aligned
\la M_0\dot\otimes\cdots\dot\otimes M_n,\dg
M_{n+1}\ra&\simeq\la M_0,M_1,\dots \dg M_{n+1}\ra\cr
&\simeq\Hom_{\hat\gb}(M_0\dot\otimes\cdots\dot\otimes
M_{n},\dg DM_{n+1})^*.
\endaligned
$$


(c) The bifunctor $\dot\otimes$ 
takes $\hat\Oc^\Delta_{\dom,\kappa}\times\hat\Oc^\Delta_{\nu,\kappa}$
and
$\hat\Oc^\Delta_{\nu,\kappa}\times\hat\Oc^\Delta_{\dom,\kappa}$
into $\hat\Oc^\Delta_{\nu,\kappa}$.
\endproclaim

\proclaim{7.3.~Corollary} 
For $M\in\hat\Oc_{\dom,\kappa}^\f$
the functor $R_M:\hat\Oc_{\nub}^\f\to\hat\Oc_{\nub}^\f$, 
$N\mapsto N\dot\otimes M$ is exact.
If the module $M$ is standardly filtered then $R_M$ preserves the subcategory
$\hat\Oc^\Delta_{\nu,\kappa}$.
The functors $R_M$, $R_{DM}$ are
adjoint (left and right) to each other.
\endproclaim

In what follows we will omit from notation the associativity and unit isomorphisms, as is justified by the Mac Lane coherence theorem. 

\subhead 7.4.~The KZ-functor\endsubhead
Fix $q\in\CC^\times$ and a tuple 
$Q=(q_1,\dots q_\ell)$ in $(\CC^\times)^{\ell}$.
The Ariki-Koike algebra
$\Ab_{q,Q}$ is the $\CC$-algebra with 1
generated by $T_0,T_1,\dots T_{n-1}$ modulo the
defining relations
$$\gathered
(T_0-q_1)\cdots(T_0-q_\ell)=(T_i+1)(T_i-q)=0,
\vspace{2mm}
T_0T_1T_0T_1=T_1T_0T_1T_0, \quad T_iT_{i+1}T_i=T_{i+1}T_iT_{i+1},
\quad T_iT_j=T_jT_i.
\endgathered$$
Assume that
$$q=\exp(2i\pi h),
\quad q_p=\exp\bigl(2i\pi (h_1+h_2+\cdots
h_{p-1}+(p-1)/\ell)\bigr),\quad \forall p\in\Lambda.\leqno(7.1)$$
Let $\reg:\Hc_{h,H}\to\Hc_{h,H,\reg}$ be the
quotient by the Serre subcategory generated by the modules $M$ such
that $M_{n,\ell}=0$. 
By \cite{GGOR, sec.~5.1,5.3} the
Riemann-Hilbert correspondence yields an equivalence of categories
$$\KZ : \Hc_{h,H,\reg}\to\Ab_{q,Q}\text{-}\modb.$$
Composing $\Een$, $\heartsuit$ and KZ yields the functor
$$\VVKZ:\hat\Oc_{\nub}\to\Ab_{q,Q}\text{-}\modb.$$

\subhead 7.5.~Comparison of $\VVKZ$ with the Kazhdan-Lusztig tensor product
\endsubhead 
Recall that $\Vb$ is the dual of the vectorial representation of $\gen$
and that $\Vb^*_\kappa=D(\Vb_\kappa)$.
Consider the following $\hat\gb_\kappa$-module
$$\Vb_{n,\nub}=\Delta_{0,\nub} \,
\dot\otimes (\Vb_\kappa)^{\dot\otimes n}.$$
It lies in the category $\hat\Oc_{\nub}^\f$ by
Proposition 7.2$(a)$. So the $\CC$-algebra 
$$\Ab_{n,\nub}=\End_{\hat\gb}(\Vb_{n,\nub})$$
is finite dimensional. The duality of $\CC$-vector spaces $M\mapsto M^*$ 
yields a functor
$$\delta:\Ab_{n,\nub}^\op\text{-}\modb\to\Ab_{n,\nub}\text{-}\modb.$$
Composing the functor
$$\Fen^\op:\hat\Oc_{\nub}\to \Ab_{n,\nub}^\op\text{-}\modb,\
M\mapsto\Hom_{\hat\gb}(\Vb_{n,\nub},M)$$ 
with $\delta$ yields the functor
$$\Fen:\hat\Oc_{\nub}\to \Ab_{n,\nub}\text{-}\modb.$$

\proclaim{7.6.~Proposition} There is an algebra homomorphism
$\phi:\Ab_{q,Q}\to \Ab_{n,\nub}$ such that
$\VVKZ=\phi\circ\Fen\circ \dg D.$
\endproclaim

\noindent{\sl Proof :}  Let $\VVCC$ be the
composition of $\VVKZ$ and the forgetful functor
$$\Ab_{q,Q}\text{-}\modb\to\CC\text{-}\modb.$$ Let
$\W_{\CC}$ be the composition of $\Fen$ and the forget
functor
$$\Ab_{n,\nub}\text{-}\modb\to\CC\text{-}\modb.$$ The proof consists of
two steps : first we construct an isomorphism of functors
$$\psi:\VVCC\circ \dg D\to\W_{\CC},$$ then we define an algebra
homomorphism $$\phi:\Ab_{q,Q}\to \Ab_{n,\nub}$$ such that $\psi$
lifts to an isomorphism of functors $$\VVKZ\circ
\dg D\to\phi\circ\Fen.$$
Now fix a module $N$ in $\hat\Oc^\f_{\nub}$. Proposition 7.2$(b)$
yields
$$\aligned
\Fen_\CC(N)= 
\la M(\pi+c\o_0)_\nu,\Vb_\kappa,\dots \Vb_\kappa, DN\ra.
\endaligned$$
Consider the $\hat\gb^F$-module
$M=\vk N$. By Proposition 4.2 we
have $$M(c\o_0)^{F'}=\vk(\dg M(\pi+c\o_0)_\nu).$$
Recall that we have 
$$\Vb^*_\kappa=\dg\,\Vb_\kappa.$$ Fix tuples $x\in C_\nl$ and $y\in
C_{n,1}$ as in Section 4.3. Assume that $x,y$ belong to the set $\Cc$ defined in Remark A.2.2.
Applying successively
Propositions 2.14$(d)$, Propositions 4.2$(b)$ and 4.6,
and Proposition 3.11$(b)$, we get natural isomorphisms
$$\aligned
\Fen_\CC(N)&\simeq\la \dg M(\pi+c\o_0)_\nu,\Vb^*_\kappa,\dots
\Vb^*_\kappa,\dg DN\ra_{y}\cr &\simeq\la M(c\o_0)^{F'},\Vb^*_\kappa,\dots
\Vb^*_\kappa, \dg DM\ra_x\cr &\simeq\la \Vb^*_\kappa,\dots \Vb^*_\kappa,
\dg DM\ra_x.\endaligned$$
By Proposition 3.11$(a)$ the vector space $\la \Vb^*_\kappa,\dots
\Vb^*_\kappa, \dg DM\ra_x$ is the fiber at the point $x$ of the
$W$-equivariant locally free sheaf over $\Alocnl$ associated with
the $\Rlocnl$-module $\la \Vb^*_\kappa,\dots \Vb^*_\kappa,\dg DM\ra.$ 
Remark 3.12 yields
$$\la \Vb^*_\kappa,\dots \Vb^*_\kappa,\dg DM\ra
=\Cen(\dg DM)_\nl.$$ Further the $\Rnl$-module $\Cen(\dg DM)_\nl$ is equipped
with a flat $W$-equivariant connection $\nabla$ which comes, via
Proposition 1.8, from the representation of $\Hb_{h,H}$ on $\Cen(\dg DM)$
in Proposition 3.6$(b)$.

Let $p:\tilde C_\nl\to C_\nl$ be the universal cover
of the complex manifold associated with the
$\CC$-scheme $C_\nl$. By definition of the functor $\Een$ 
we have $$\Cen(\dg DM)=\V(\dg DN).$$ Since the vector space $\VVCC(\dg DN)$ is
obtained from the $\Ab_{q,Q}$-module $\V_\KZ(\dg DN)$ by forgetting the
$\Ab_{q,Q}$-action, it is canonically identified with the vector
space of holomorphic horizontal sections of the connection $\nabla$
over $\tilde C_\nl$. Recall that $\Cc$ is a contractible subset of
$C_\nl$ containing the point $x$. Fix once for all a
contractible subset $\tilde \Cc\subset\tilde C_\nl$ such that $p$
restricts to an isomorphism $\tilde \Cc\to \Cc$. Restricting
functions on $\tilde C_\nl$ to $\tilde \Cc$ and taking the fiber at
$x$, viewed as a point of $\tilde \Cc$ via the map $p$, yields a
natural isomorphism of vector spaces
$$\psi(N):\VVCC(\dg DN)\simeq\la
\Vb^*_\kappa,\dots \Vb^*_\kappa,\dg DM\ra_x\simeq\Fen_\CC(N).$$
Further, the $\Ab_{q,Q}$-action on the functor $\VVCC$ given by $\VVKZ$
gives, via the isomorphism $\psi$, a $\Ab_{q,Q}$-action on the
functor $\W_{\CC}$. This action comes from an algebra homomorphism
$\phi$ as above, by the following consequence of Yoneda's lemma.

\proclaim{7.7.~Lemma} Let $M$ be an object of an additive category $\Ac$.
If a ring $\Ab$ acts on the functor $\Hom_\Ac(M,-)$,
then there is a ring homomorphism $\phi:\Ab\to\End_\Ac(M)^\op$ such that
the action of any element $a\in\Ab$ on $F$ is the composition by
$\phi(a)$.
\endproclaim

\qed

\vskip3mm

\proclaim{7.8.~Corollary}
The functor $\Een$ restricts to an exact functor
$\hat\Oc^\Delta_{\nub}\to \Hc_{h,H}^\Delta$.
\endproclaim

\noindent{\sl Proof :} First, let us note the following basic fact
whose proof is left to the reader.

\proclaim{7.9.~Lemma} Let $\Een:\Ac\to\Bc$ be a right exact
functor of quasi-hereditary categories such that
$\Een(\Delta_\Ac)\subset\Delta_{\Bc}\cup\{0\}$. Let
$\heartsuit:\Bc\to\Bc_\heartsuit$ be an exact functor such
that $\Hom_{\Bc}(M,N)=0$ for each $M\in\Ker(\heartsuit)$,
$N\in\Bc^\Delta$. If $\heartsuit\circ\Een$ is exact on
$\Ac^\Delta$ then $\Een$ restricts to an exact functor
$\Ac^\Delta\to\Bc^\Delta$.
\endproclaim

\noindent There are no nonzero homomorphisms $M\to N$ for each
$M\in\Ker(\heartsuit)$ and $N\in\Hc_{h,H}^\Delta$, because a standard
$\Hb_{h,H}$-module is torsion free as a $\Rb$-module. To prove that
$\heartsuit\circ\Een$ is exact on $\hat\Oc_{\nub}^\Delta$, it is enough
to check it for $\Een_\KZ$. Thus the claim follows from Proposition
7.6, because the module $\Vb_{n,\nub}$ is standardly filtered
by Proposition 7.2$(c)$.

\qed

\vskip3mm

\proclaim{7.10.~Conjecture} If $\lub\Lambda$ is large enough then we have
$\Een(\lub P_{\l,\nub})=P_{\l^\circ,h,H}$ for each $\l\in\Pc_{n,\nu}^\ell$.
\endproclaim

\subhead 7.11.~Remarks\endsubhead $(a)$ For $\ell>1$ the functor
$\Een$ may be not exact. It may also take a simple object to a
non-simple non-zero $\Hb_{h,H}$-module. Indeed, set $\kappa=-1$,
$m=7$, $\ell=4$, $\nu=(1,1,4,1)$, $n=1$, as in Example 6.5$(b)$. We
have $\Een(\Delta_{1_p,\nub})=\Delta_{1^\circ_p,h,H}$ for each $p$.
We have also $\Een(\Delta_{\l_1,\nub})=0$ and there is an exact
sequence $\Een(\Delta_{\l_1,\nub})\to\Een(\Delta_{1_1,\nub})\to
\Een(S_{1_1,\nub})\to 0$. Thus
$\Een(S_{1_1,\nub})=\Delta_{1^\circ_1,h,H}$, which is not simple.
Further we have $S_{1^\circ_4,h,H}=\Delta_{1^\circ_4,h,H}$ and there
is an exact sequence $0\to S_{1^\circ_4,h,H}\to
\Delta_{1^\circ_1,h,H}\to S_{1^\circ_1,h,H}\to 0$. Thus the derived
functor $L^{-1}\Een$ is non-zero and it takes $S_{\l_1,\nub}$ to
$S_{1^\circ_4,h,H}$.

\vskip1mm

$(b)$ The morphism $\phi$ in Proposition 7.6 is injective.
We do not know if it is invertible.
To prove this we must check that the representation of $\Ab_{q,Q}$ on
$\Een_\KZ(M)$ is faithful for some module $M$. Since $\KZ$ is exact
and $\Hc_{h,H}$ has enough projective objects, there is a module
$P_\KZ\in\Hc_{h,H}^\proj$ which represents $\KZ$. We have
$\Ab_{q,Q}=\KZ(P_\KZ)$ by \cite{GGOR, thm.~5.15}. 
We claim that if $\lub\Lambda$
is large enough  there is a projective  module $P\in\lub\hat\Oc_{\nub}$
such that $P_\KZ$ is a direct summand of $\Een(P)$. Therefore the
representation of $\Ab_{q,Q}$ on $\Een_\KZ(P)$ is faithful. The proof
is omitted because we'll not use this result.

\vskip1mm

$(c)$ Recall that we have $[\hat\Oc_{\nub}^\f]=[\hat\Oc^\Delta_{\nub}]$ and
$[\Hc_{h,H}^\Delta]=[\Hc^\f_{h,H}]$. The derived functor $L^*\Een$
yields a group homomorphism $[\hat\Oc_{\nub}^\f]\to[\Hc^\f_{h,H}]$
such that $[\Delta_{\l,\nub}]\mapsto[\Delta_{\l^\circ,h,H}]$ if
$\l\in\Pc_{n,\nu}^\ell$ and $[\Delta_{\l,\nub}]\mapsto 0$ else.
Conjecture 7.10 implies that $[\lub P_{\mu,\nub}]$ maps to
$[P_{\mu^\circ,h,H}]$ for all $\mu\in\Pc_{n,\nu}^\ell$. So Brauer
reciprocity implies that
$$[\nabla_{\l^\circ,h,H}:S_{\mu^\circ,h,H}]=
[\nabla_{\l,\nub}:S_{\mu,\nub}],\leqno(7.3)$$ where
$\nabla_{\l^\circ,h,H}$, $\nabla_{\l,\nu}$ are the costandard
modules with socles $S_{\l^\circ,h,H}$, $S_{\l,\nu}$ in $\Hc_{h,H}$,
$\hat\Oc_\nub$ respectively. This dimension formula is not the same as
in the dimension conjecture. See Section 8 for a comparison of the
two formulas. We do not know if they are equivalent.

\vskip1mm

$(d)$ The module $\Een(\lub P_{\l,\nub})$ does not depend on the set
$B$ if $\l\in\Pc_{n,\nu}^\ell$ and $\lub\Lambda$ is large enough so that
$\tilde\mu_\pi\in\lub\Lambda$ for all $\mu\in\Pc_{n,\nu}^\ell$. Indeed,
if $B\subset B'$ and $M\in{}^{B'}\hat\Oc_{\nub}$ let $\lub M$ be the
maximal quotient of $M$ which belongs to $\lub\hat\Oc_{\nub}$. The
functor $M\mapsto\lub M$ maps ${}^{B'}P_{\l,\nub}$ to $\lub
P_{\l,\nub}$, see \cite{Do, sec.~A.1}. Brauer reciprocity implies
that the kernel of the obvious projection ${}^{B'}P_{\l,\nub}\to\lub
P_{\l,\nub}$ is filtered by parabolic Verma modules whose highest
weights belong to ${}^{B'}\Lambda\setminus\lub\Lambda$. The claim
follows.

\vskip2cm

\head 8.~The affine parabolic category $\Oc$ and the Fock space\endhead

Fix integers $m,\ell>0$, $e>1$ and fix  a composition
$s\in\Cc_{m,\ell}$. 
In this section we'll use freely the notation from the appendix. In
particular, we have defined there the integers
$\Delta^+_{\l,\mu,e,-s}$, $\nabla^-_{\l,\mu,e,s^\circ}$, the
$\widehat{\sen\len}_e$-module $\Lambda^s$, and the basis elements
$|\l,s,e,s^\circ\ra$, $\Gc(\l,s,e,s^\circ)^-$ of $\Lambda^s$.

\subhead 8.1.~The affine parabolic category $\Oc$ and the Fock space
\endsubhead
The level $\ell$ {\it Fock space} $Fock_{e,s^\circ}$ associated
with the multicharge $s^\circ$ is a $\widehat{\sen\len}_e$-module of
level $\ell$ which is the limit of a filtered inductive system of
vector spaces $\Lambda^r$, with $r$ a positive integer. Each
$\Lambda^r$ is equipped with a level zero representation of
$\widehat{\sen\len}_e$. See Section A.4 for details.

Now, let $\kappa=-e$.
For each weight $\l\in\ten^*$ we'll write
$\hat\l=\l+(\kappa-m)\o_0$.
Let $\pi\in\ten^*$ be given by
$$\pi+\rho=(s_1,s_1-1,\dots
1,s_2,s_2-1,\dots,s_{\ell},s_{\ell}-1,\dots 1).$$
Note that $\pi$ is NOT the same as in Section 6.1.
Given an $s$-dominant weight $\l$ we
abbreviate
$$\Delta_{\l,s,-e}=M(\hat\l_{\pi})_s,\quad S_{\l,s,-e}=L(\hat\l_{\pi}).$$ Let
$\Ac_{s,-e}\subset\hat\Oc_\kappa$
be the Serre category
generated by the modules $S_{\l,s,-e}$, $\l\in\ZZ^s_\dom$.
The proposition below identifies the vector space
$[\Ac_{s,-e}^\f]\otimes\CC$ with a submodule
$\Lambda^s$ of the $\widehat{\sen\len}_e$-submodule $\Lambda^{m}$.
Thus we can regard $[\Ac_{s,-e}^\f]\otimes\CC$ as a subspace of
$Fock_{e,s^\circ}$.

\proclaim{8.2.~Proposition} There is a vector space isomorphism
$$[\Ac_{s,-e}^\f]\otimes\CC\to \Lambda^s,\quad
[\Delta_{\l,s,-e}]\mapsto|\l,s,e,s^\circ\ra,\quad
[S_{\l,s,-e}]\mapsto\Gc(\l,s,e,s^\circ)^-.$$ We have
$\nabla^-_{\l^\circ,\mu^\circ,e,s^\circ}=[\Delta_{\l,s,-e}:
S_{\mu,s,-e}]$ for each $\l,\mu\in\Pc_{n,s}$.
\endproclaim

\noindent{\sl Proof :} We have $\l_{\pi}+\rho=\a(\l,s,s)$ by Remark
A.4.4$(a)$. Proposition 5.8 and Proposition A.4.3$(b)$ yield
$$[S_{\mu,s,-e}]=\sum_{\l}
(-1)^{l(v_\l)-l(v_\mu)}P^{\g,-1}_{v_\l,v_\mu}(1)
[\Delta_{\l,s,-e}],$$
$$\Gc(\mu,s,e,s^\circ)^-=\sum_{\l}
(-1)^{l(v_\l)-l(v_\mu)}P^{\g,-1}_{v_\l,v_\mu}(1)
|\l,s,e,s^\circ\ra,$$ where $\g\in\alc$, $v_\l,
v_\mu\in\hat\Sen^{\g}$ such that $v_\l\bullet\g=\hat\l_{\pi}$,
$v_\mu\bullet\g=\hat\mu_{\pi}$ and $v_\mu\geqslant v_\l$. This proves
the first claim.  The second one follows from Proposition A.4.3$(c)$.

\qed

\vskip3mm

\subhead 8.3.~Yvonne's conjecture\endsubhead
Fix an integer $n>0$, $q\in\CC^\times$ and a tuple 
$Q=(q_1,\dots q_\ell)$ in $(\CC^\times)^{\ell}$.
Let $\Ab_{q,Q}$ be defined as in Section 7.4.
The {\it cyclotomic $q$-Schur algebra} $\Sb_{q,Q}$ is the endomorphism
algebra of a particular projective $\Ab_{q,Q}$-module. It is a
quasihereditary algebra. In particular, for each $\l\in\Pc_n^\ell$
there is a standard module $\Delta_{\l,q,Q}$ with a simple top
$S_{\l,q,Q}$. See \cite{DJM, def.~6.13} for details.
Yvonne's conjecture
\cite{Yv, conj.~2.13, def.~2.5,4.4} is the following one
(note that our hypothesis differ slightly from the ones in loc.~cit.).

\proclaim{8.4.~Conjecture} Assume that $q=\exp(-2i\pi/e),$
$q_p=\exp(2i\pi s_p/e)$ and $s_{p+1}-s_{p}\geqslant n$ for
$p\neq\ell$. Then we have
$[\Delta_{\l,q,Q}:S_{\mu,q,Q}]=\Delta^+_{{}^t\!\mu,{}^t\!\l,e,-s}$
for all $\l,\mu\in\Pc_n^\ell$.
\endproclaim

\vskip3mm

\subhead 8.5.~The dimension conjecture\endsubhead  Let
$q$, $Q$ be as in (7.1). Assume that
$$(q+1)\prod_{p'\neq p''}(q_{p'}-q_{p''})\neq 0,\quad
h_p\geqslant(1-n)h,\quad h<0,\quad \forall p\neq\ell.\leqno(8.1)$$
We have the following \cite{R, thm.~6.8}.

\proclaim{8.6. ~Theorem} If $(7.1)$, $(8.1)$ hold there is an equivalence
of quasi-hereditary categories $\Hc_{h,H}\to\Sb_{q,Q}$-$\modb$ taking
$\Delta_{\l,h,H}$ to $\Delta_{\l,q,Q}$.
\endproclaim

\noindent Now, assume that (7.1), (8.1) hold and fix $\l,\mu$ in $\Pc^\ell_n$.
Theorem 8.6 yields
$$[\Delta_{\l,h,H}:S_{\mu,h,H}]=[\Delta_{\l,q,Q}:S_{\mu,q,Q}].$$ 
Assume further that the following hold
$$h=-1/e,\quad h_p=s_{p+1}/e-s_{p}/e-1/\ell,\quad \forall
p\neq\ell.\leqno(8.2)$$ 
The {\it dimension conjecture}
compares also $[\Hc_{h,H}^\f]\otimes\CC$ with a subspace of
$Fock_{e,s^\circ}$. We'll state a categorical analogue of this
conjecture. It simply claims that $\Hc_{h,H}$ should be equivalent to
a full subcategory $\Ac_{n,s,-e}$ of $\Ac_{s,-e}$ if
$s\in\Cc_{m,\ell,n}$. 
Taking the Grothendieck groups we recover the usual dimension conjecture
from its categorical version by Proposition 8.2.
More precisely, let first note that (7.1), (8.1) and (8.2) imply that
$$q=\exp(-2i\pi/e),\quad q_p=\exp(2i\pi s_p/e),\quad s_{p+1}-s_p\geqslant n-1+e/\ell.$$ Assume also that the hypothesis
in Conjecture 8.4 holds. Then we should have
$$[\Delta_{\l,h,H}:S_{\mu,h,H}]
=\Delta^+_{{}^t\!\mu,{}^t\!\l,e,-s}.\leqno(8.3)$$
The dimension conjecture says
that this equality should hold without the lower bound in (8.1)
on the parameters $h_p$, see \cite{R, sec.~6.5}.
If $s\in\Cc_{m,\ell,n}$ then $(A.6)$ and Proposition 8.2
yield the equality
$$\Delta^+_{{}^t\!\mu,{}^t\!\l,e,-s}=
[\Delta_{\l^\circ,s,-e}:S_{\mu^\circ,s,-e}].\leqno(8.4)$$
Therefore, composing the equalities (8.3) and (8.4) we get the following one
$$[\Delta_{\l,h,H}:S_{\mu,h,H}]=[\Delta_{\l^\circ,s,-e}:S_{\mu^\circ,s,-e}].$$
Before formulating the {\it categorical dimension conjecture}
note the following easy fact, see
Section A.6. Let $\Ac_{n,s,-e}$ be the Serre
subcategory of $\Ac_{s,-e}$ generated by the modules $S_{\l,s,-e}$ with
$\l\in\Pc_{n,s}$.

\proclaim{8.7.~Proposition} The category $\Ac_{n,s,-e}$ is
quasi-hereditary with respect to the order $\trianglelefteq$.
The standard modules are the
modules $\Delta_{\l,s,-e}$ with $\l\in\Pc_{n,s}$.
\endproclaim

\noindent We conjecture the following.

\proclaim{8.8.~Categorical dimension conjecture} 
If $s\in\Cc_{m,\ell,n}$ and $(8.2)$
holds there is an equivalence of categories
$\Ac_{n,s,-e}\to\Hc_{h,H}$ taking 
$\Delta_{\l^\circ,s,-e}$ to $\Delta_{\l,h,H}$.
\endproclaim

\subhead 8.9.~Comparison of the dimension conjecture with the functor
$\Een$\endsubhead
Fix an integer $m>0$ and a composition $\nu\in\Cc_{m,\ell}$.
Let $\Ac_{n,\nub}\subset\hat\Oc_\nub$ be the Serre subcategory
generated by the simple modules $S_{\l,\nub}$ with
$\l\in\Pc_{n,\nu}^\ell$ (with the notation from Section 6.1).
We define $h,H$ in the following way
$$\kappa=-e,\quad
h=1/\kappa,\quad h_p=\nu_p^\bullet/\kappa-m\ell/\kappa,\quad 
\forall p\in\Lambda.\leqno(8.5)$$
Then the functor $\Een$ restricts to a functor
$$\Een:\Ac_{n,\nub}\to\Hc_{h,H}.\leqno(8.6)$$ 
Note that the integer $m$ in Conjecture 8.8 and the integer $m$
in (8.6) are not the same.
If $\ell=1$ we may choose them to be equal. Then $\Een$ is Suzuki's functor
and it yields the equivalence 
$$\Ac_{n,s,-e}\to\Hc_{h,H}$$ 
which is conjectured in 8.8, by Theorem A.5.1. 
If $\ell>1$ we do not know how to get a functor
$\Ac_{n,s,-e}\to\Hc_{h,H}$ from $\Een$ because, given $h$, $H$,
the composition $s$ and the integer $m$ in (8.2) 
are different from the composition $\nu$
and the integer $m$ in (8.5).

\vskip3mm

\subhead 8.10.~Remarks\endsubhead
$(a)$
Note that $\Ac_{n,\nub}$ is not a quasi-hereditary category for the
order $\trianglelefteq$ on the parabolic Verma modules,
see Example 6.6$(b)$.

\vskip1mm

$(b)$
Let $\Sb'$ be the degenerate analogue of $\Sb_{q,Q}$ considered in
\cite{BK}. Assume that $s\in\Cc_{m,\ell,n}$ and (8.2) hold. Here we
allow $e$ to be any nonzero complex number. Now, assume also that
$e\notin\QQ$. Then, it is not difficult to prove that
the category $\Ac_{n,s,-e}$ is equivalent to a
parabolic subcategory $\Ac_{n,s}\subset\Oc$.
This is due to the fact that if $\kappa\notin\QQ$ then
the induction functor $\Gamma:\Oc_\nu^\f\to\hat\Oc_{\nub}^\f$ is
is an equivalence of categories and also to the fact that
there is a canonical isomorphism
$$\Gamma(M)\dot\otimes\Gamma(E)=\Gamma(M\otimes E)$$
for all $M\in\Oc^\f_\nu$ and all finite dimensional $\gen$-module $E$
(the proof is standard and will be given elsewhere).
It is probably not difficult to prove that, under this
assumption, the categories $\Hc_{h,H}$ and $\Sb'$-$\modb$ are equivalent.
We expect that, in this case, the equivalence in Conjecture 8.8 is
precisely the equivalence of categories $\Ac_{n,s}\to\Sb'$-$\modb$ in
\cite{BK, thm.~C}.

\vskip3mm

\subhead 8.11.~Example\endsubhead
Assume that $m=7$, $\ell=4$, $\kappa=-1$, $\nu=(1,1,4,1)$, $n=1$ and
$h=-1$, $H=(-9/4,3/4,3/4)$. Then we have $\nu\in\Cc_{m,\ell,n}$ and
(8.5) holds. Now set $m'=9$ and $e=1$, $s=(3,1,2,3)$. Then we have
$s\in\Cc_{m',\ell,n}$ and (8.2) holds. Note that for this choice of
$m$, $\ell$ and $\nu$ the equality (7.3) holds (a computation yields
$[\nabla_{\l,\nub}:S_{\mu,\nub}]=1$ if
$\Delta_{\mu,\nub}\trianglelefteq\Delta_{\l,\nub}$ and 0 else).

\vskip2cm

\head A. Appendix\endhead

\subhead A.1. Proof of Proposition 5.8\endsubhead
Fix $\l\in\Aen_\rho$ and fix $w\in\hat\Sen(\l)^{\l}$ such that $w\bullet\l$ is $\nu$-dominant.
Since $w\in \hat\Sen(\l)^{\l}$ and $\la\l+\rhoaf:\a\ra\leqslant 0$ for all
$\a\in\Pib(\l)^+$, we have the following formula in
$[\hat\Oc_\kappa^\f]$
$$[L(w\bullet\l)]=
\sum_{v\in\hat\Sen(\l)}(-1)^{l(w)-l(v)}P_{v,w}(1)[M(v\bullet\l)].$$ Here
the sum is over all $v$'s such that $w\geqslant v$ and $P_{v,w}$ is
the Kazhdan-Lusztig polynomial relative to $\hat\Sen(\l)$. See
\cite{KT, thm.~1.1} for details. Since $w\in \hat\Sen(\l)^{\l}$, for
each $u\in\hat\Sen(\l)^{\l}$ we have also
$$\sum_{x\in\hat\Sen_{\l}}(-1)^{l(x)}P_{ux,w}(1)=
P^{\l,-1}_{u,w}(1)$$ by definition of the parabolic Kazhdan-Lusztig
polynomial of type $\hat\Sen(\l)/\hat\Sen_{\l}$. This yields the following
formula
$$[L(w\bullet\l)]=
\sum_{u\in\hat\Sen(\l)^{\l}}(-1)^{l(w)-l(u)}P_{u,w}^{\l,-1}(1)
[M(u\bullet\l)]. \leqno(A.1)$$
Next, observe that $\Sen_\nu$ is a parabolic subgroup of $\hat\Sen(\l)$. Indeed, given a simple root $\a_i$ which belongs to 
$\Pi_\nu$
we write
$$\la\l+\rhoaf:\a_i\ra=
\la\l+\rhoaf-w(\l+\rhoaf):\a_i\ra+\la w\bullet\l:\a_i\ra+\la\rhoaf:\a_i\ra.$$ Then the first term
is an integer because $w\in\hat\Sen(\l)$, and the second one is also an
integer because $w\bullet\l$ is $\nu$-dominant. Thus we have
$$\Sen_\nu\subset\hat\Sen(\l).$$
We must check that if the simple root $\a_i$ lies in $\Pi_\nu$ then it belongs also to the
basis of $\hat\Pi(\l)^+$.  The
latter consists of the real affine roots $\a$ such that
$s_\a(\Pib(\l)^+\setminus\{\a\})\subset\Pib^+(\l)$. So the
claim is obvious, because
$$s_{\a_i}(\Pib^+\setminus\{\a_i\})\subset\Pib^+,\quad
s_{\a_i}(\Pib(\l))\subset\Pib(\l).$$ 
Now, consider the set
$$S=
\{u\in\hat\Sen(\l)^{\l}\,;\,u>s_{\a_i}u,\,\forall s_{\a_i}\in\Sen_\nu\}.$$ 
Since $\l\in\alc$, for each $u\in\hat\Sen(\l)^{\l}$ we have
$$u\bullet\l\ \roman{is}\ \nu\text{-dominant} \iff
u\in S.$$ Further, it is
well-known that if $u,v\in \hat\Sen(\l)^{\l}$ are such that $u>s_{\a_i}u$ and
$v>s_{\a_i}v$, then 
$$s_{\a_i}u, s_{\a_i}v\in\hat\Sen(\l)^{\l},\quad P^{\l,-1}_{s_{\a_i}u,v}(1)=P^{\l,-1}_{u,v}(1).
\leqno(A.2)$$  Therefore if $u,v\in\hat\Sen(\l)^{\l}$ and $u\bullet\l$, $v\bullet\l$ are both $\nu$-dominant
then we have
$$P^{\l,-1}_{s_{\a_i}u,v}(1)=P^{\l,-1}_{u,v}(1),\quad\forall\a_i\in\Pi_\nu.$$
Note also that $w\in S$ and that if $u\bullet\l$ is $\nu$-dominant then the BGG resolution yields
$$[M(u\bullet\l)_\nu]=
\sum_{x\in \Sen_\nu}(-1)^{l(x)} [M(xu\bullet\l)].$$ Thus the
contribution of the elements $u\in\hat\Sen(\l)^{\l}$ in the sum (A.1) which 	are of the
form $$u=xv,\quad x\in\Sen_\nu,\quad v\in S$$ is equal
to the sum
$$\sum_{v}
(-1)^{l(w)-l(v)}P_{v,w}^{\l,-1}(1) [M(v\bullet\l)_\nu]$$ over all
elements $v\in \hat\Sen(\l)^{\l}$ such that $w\geqslant v$ and
$v\bullet\l$ is $\nu$-dominant. Using (A.2) once again one can easily
check that the other $u$'s do not contribute to the right hand side of (A.1).

\qed

\vskip3mm

\subhead A.2. Proof of Proposition 7.2 and Corollary 7.3\endsubhead

\subhead A.2.1.~The Kazhdan-Lusztig tensor product\endsubhead
First, let us recall the definition of the Kazhdan-Lusztig tensor product. 
Let $R$ be a commutative $\CC$-algebra with 1.
Assume that $R$ is a Noetherian integral domain. Let $F$ be its fraction field.
Fix a unit $\kappa\in R^\times$, and set $S=\{1,2,\dots,n\}$.  
We'll use the same notation as in Sections 2.11, 2.12.
Recall that $z$ is a local coordinate on $\PP^1$ centered at 0 and that
$x=(x_i;i\in S)$ is a familly of distinct points of $\PP^1$.
Consider the Lie algebra
$$\Gamma_{R,x}=\gen[\PP^1_x]_R.$$
Fix a family $z=(z_i;i\in S)$ of local coordinates on
$\PP^1$ such that the
coordinate $z_i$ is centered at $x_i$. We can regard $z_i$ as an
isomorphism $\PP^1\to\PP^1$. To each rational function $f\in
F(\PP^1)$ and to each index $i$ we associate the power series
expansion in $F((z))$ of the rational function $f\circ z_i$.
Composing this formal series with the assignment $z\mapsto t_i$ we
get a field homomorphism $F(\PP^1)\to F((t_i))$.
Their sum gives a $R$-algebra homomorphism $R[\PP^1_{x}]\to
R((t_S))$ and a $R$-Lie algebra homomorphism 
$$\Gamma_{R,x}\to\hat\Gc_{R,S}.$$
Now, we set $\hat S=\{0,1,\dots n\}$. Fix a family $\hat x=(x_i;i\in\hat S)$ of distinct points of $\PP^1$. 
Let $\hat\Gamma_{R,\hat x}$ be the central extension of the Lie
algebra $\Gamma_{R,x}$ by $R$ associated with the cocycle
$(\xi_1\otimes f_1,\xi_2\otimes f_2)\mapsto\Res_{x_0}(f_2df_1)$. Let
$$U(\hat\Gamma_{R,\hat x})\to\hat\Gamma_{R,\hat x,\kappa}$$ 
be the quotient by the ideal
generated by the element $\un-c$. 
The expansion at the
points $x_i$, $i\in S$, gives a $R$-algebra homomorphism $R[\PP^
1_{\hat x}]\to R((t_S))$ and a $R$-algebra homomorphism
$$\hat\Gamma_{R,\hat x,2m-\kappa}\to\hat\Gc_{R,S,\kappa}.\leqno(A.3)$$ 
The expansion at the point
$x_0$ yields a $R$-algebra homomorphism $R[\PP^ 1_{\hat x}]\to
R((t))$ and a $R$-algebra homomorphism
$$\hat\Gamma_{R,\hat x,2m-\kappa}\to\hat\gb_{R,2m-\kappa}.\leqno(A.4)$$
See \cite{KL, sec.~4.6, 8.2} for details.
Now, for each $M_i\in\Cc(\hat\gb_{R,\kappa})$,
$i\in S$, the tensor product $$W=\bigotimes_{i\in S}M_i$$
(over $R$) has an obvious structure of
$\hat\Gc_{R,S,\kappa}$-module. We equip $W$ with the structure of a
$\hat\Gamma_{R,\hat x,2m-\kappa}$-module via the map $(A.3)$. 
For each integer $r>0$ let
$$G_{R,r}\subset U(\hat\Gamma_{R})$$
be the $R$-submodule generated by the products of $r$ elements in
$$\gen\otimes\{f;f(x_0)=0\}. $$
Using $(A.4)$ we equip the projective limit of $R$-modules
$$\widehat W=\pro_r W_r,\quad W_r=W/G_{R,r} W$$
with the structure of a $\hat\gb_{R,\kappa^\sharp}$-module. We define
$$M_1\dot\otimes_R M_2\dot\otimes_R\cdots \dot\otimes_R M_n=
\sp\,T(W), \quad T(W)=\widehat W(-\infty).$$ To summarize, if $n=2$
then for any almost smooth $\hat\gb_{R,\kappa}$-modules $M_1$, $M_2$
the (smooth) $\hat\gb_{R,\kappa}$-module
$M_1\dot\otimes_RM_2$ is given by
$$M_1\dot\otimes_R M_2={}^\sharp T(W), \quad T(W)=\widehat
W(-\infty),\quad\widehat W=\pro_rW/G_rW,\quad W=M_1\otimes_RM_2.$$
See \cite{KL, sec.~4.9, 8.4} for details.
As above, if $R=\CC$ we simply forget the subscript $R$ everywhere. 

\vskip3mm

\subhead A.2.2.~Remark\endsubhead  
The functor $\dot\otimes$ depends on the choice of the tuple 
of distinct points $x=(x_i;i\in S)$ 
and on the choice of the coordinate $z_i$ centered at $x_i$ for each $i$, 
see \cite{KL, sec.~9}. 
Unless mentioned otherwise, we'll choose once for all the coordinates 
as in Section 2.12, i.e., 
we set $z_i=z-x_i$ if $x_i\neq\infty$ and $z_i=-z^{-1}$ else. 
The systems of coordinates associated with the tuples in the set
$$\Cc=\{x\in\RR^n;\ 0<x_1<\cdots<x_n\}$$ 
belong to the contractible real manifold
introduced in \cite{KL, sec.~13.1}. 
So the space of affine coinvariants, see Definition 2.13,
is independent on the choice of $x\in \Cc$ by
\cite{KL, sec.~13.3}.
We may abbreviate
$$\la M_i;i\in S\ra=\la M_i;i\in S\ra_x,\quad x\in\Cc,$$
if this does not create any confusion.

\vskip3mm

\subhead A.2.3.~Proof of Proposition 7.2$(a)$\endsubhead
If $R=\CC$ and $\kappa\notin\QQ_\dom$
the bifunctor $\dot\otimes$ yields the Kazhdan-Lusztig's monoidal category
\cite{KL, sec.~31}
$$(\hat\Oc^\f_{\dom,\kappa},\dot\otimes,a,M(c\o_0)).$$


\proclaim{A.2.4.~Proposition} Assume that $R=\CC$ and
$\kappa\notin\QQ_\dom$. The functor $\dot\otimes$ takes
$\hat\Oc^{fg}_{\dom,\kappa}\times\hat\Oc^{fg}_{\nu,\kappa}$ and
$\hat\Oc^{fg}_{\nu,\kappa}\times\hat\Oc^{fg}_{\dom,\kappa}$
into $\hat\Oc^{fg}_{\nu,\kappa}$. The tuple 
$(\hat\Oc_{\nu,\kappa}^{fg},\dot\otimes,a,M(c\o_0))$ is a bimodule
over $(\hat\Oc_{\dom,\kappa}^{fg},\dot\otimes,a,M(c\o_0))$.

\endproclaim

\noindent{\sl Proof :} The first part is proved in
\cite{Y, thm.~1.6}, and the second one is claimed there. 
Since the construction of the associativity isomorphism is the same as in
\cite{KL, sec.~18.2} we'll not give more details there. 
Note that the axioms of a bimodule over a category imply that
there is a canonical isomorphism from
$M(c\o_0)\dot\otimes-$ to the identity functor
of $\hat\Oc^{fg}_{\nu,\kappa}$. It is given by the following chain of
isomorphisms, for each modules $M_1$, $M_2$
$$\aligned
\Hom_{\hat\gb}(M_1,M_2)
&=\la M_1,D M_2\ra^*\cr
&=\la M(c\o_0),M_1,D M_2\ra^*\cr
&=\la M(c\o_0)\dot\otimes M_1,D M_2\ra^*\cr
&=\Hom_{\hat\gb}(M(c\o_0)\dot\otimes M_1,M_2).
\endaligned$$
The second isomorphism is as in Proposition 2.14$(c)$, 
the other ones are as in Proposition A.2.6 below. 
A similar construction yields an isomorphism from
$-\dot\otimes M(c\o_0)$ to the identity functor. 

\qed

\vskip3mm

\subhead A.2.5.~Proof of Proposition 7.2$(b)$\endsubhead
Let us prove the following proposition.

\proclaim{A.2.6.~Proposition} 
$(a)$ Let $M_1,\dots M_{n}\in\hat\Oc^\f_{\dom,\kappa}$ and
$M_0,M_{n+1}\in\hat\Oc^\f_{\nu,\kappa}$.
We have  natural isomorphisms of vector spaces
$$\aligned
\la M_0\dot\otimes M_1\dot\otimes\cdots\dot\otimes M_n,\dg
M_{n+1}\ra &=\Hom_{\hat\gb}(M_0\dot\otimes
M_1\dot\otimes\cdots\dot\otimes M_n, \dg DM_{n+1})^*\vspace{2mm} &=\la
M_0,M_1,\dots\dg M_{n+1}\ra.\endaligned$$

$(b)$ If
$M_i=(N_i)_\kappa$ is a generalized Weyl module for each $i\in S$
then we have a natural isomorphism of vector spaces
$\la M_1,M_2,\dots, M_n\ra
=H_0(\gen,N_1\otimes N_2\otimes\dots\otimes N_n)$.
\endproclaim

\noindent{\sl Proof :} First, observe that the module
$M_0\dot\otimes M_1\dot\otimes\cdots \dot\otimes M_n$ lies in  the
category $\hat\Oc^{fg}_{\nu,\kappa}$ by Proposition A.2.4$(a)$.
Thus to prove the first equality in $(a)$ it is enough to check that for each
modules $M$, $N$ in $\hat\Oc^{fg}_{\nu,\kappa}$ we have a natural
isomorphism of (finite dimensional) vector spaces
$$\la M, \dg N\ra^*=\Hom_{\hat\gb}(M,\dg DN).
\leqno(A.5)$$
The right hand side of $(A.5)$ consists of the families of linear forms
$$(f_{\l,\mu})\in\prod_\l\bigoplus_{\mu}(M_\l\otimes N_{\mu})^*$$
which vanish on the set
$$\{(\xi m)\otimes n+m\otimes(\ddag\xi
n);\,\forall \xi\in\hat\gb,m\in M,n\in N\}.$$ We claim that the
obvious inclusion
$$\prod_\l\bigoplus_{\mu}(M_\l\otimes N_{\mu})^*\subset
\prod_{\l,\mu}(M_\l\otimes N_{\mu})^*$$ factors to a natural
isomorphism
$$\Hom_{\hat\gb}(M,\dg DN)\to\Hom_{\hat\gb}(M\otimes\ddg N,\CC).$$
The definition of the set of affine coinvariants  yields
$$\la M,\dg N\ra=H_0(\hat\gb,M\otimes\ddg N).$$
Therefore there is a natural
isomorphism
$$\Hom_{\hat\gb}(M\otimes\ddg N,\CC)=\la M, \dg N\ra^*.$$
Now we check our claim.
The injectivity is clear. To prove surjectivity, fix
$$(f_{\l,\mu})\in\Hom_{\hat\gb}(M\otimes\ddg N,\CC),\quad
f_{\l,\mu}\in(M_\l\otimes N_\mu)^*,
$$ and fix a finite
subset $\Sc\subset\tb^*$ such that $N$ is generated by
the subspace
$\bigoplus_{\mu\in\Sc}N_{\mu}$. Fix also $\l$ and $x\in M_\l$. For
each $\mu$ the weight space $N_\mu$ is spanned by elements of the
form $y=\xi z$ with $z\in N_{\g}$, $\g\in\Sc$ and $\xi\in
U(\hat\gb_\kappa)$ of weight $\mu-\g$. For all $\mu$ except a
finite number the weight space $M_{\l-\mu+\g}$ vanishes, hence
$$f_{\l,\mu}(x\otimes y)=-f_{\l,\mu}((\ddag\xi)x\otimes z)=0.$$
So $f_{\l,\mu}=0$ for all $\mu$ except a finite number. This proves
the claim.

The vector space $\la M_0,M_1,\dots \dg M_{n+1}\ra$ is
finite dimensional, because the modules
$M_0,M_1,\dots\dg M_{n+1}$ are quotient of generalized Weyl modules.
Therefore, the same argument as in \cite{KL, sec.~13.4}
yields the second isomorphism in $(a)$.

Now we concentrate on $(b)$.
The inclusions $N_i\subset M_i$, $i\in S$, yield an inclusion
$\bigotimes_{i\in S}N_i\subset\bigotimes_{i\in S}M_i$.
Taking the coinvariants we get a natural map
$$H_0(\gen,N_1\otimes N_2\otimes \dots\otimes N_n)\to
\la M_1,M_2,\dots,M_n\ra$$
which is invertible by \cite{KL, prop.~9.15}. 

\qed

\vskip3mm



\subhead A.2.7.~Proof of Proposition 7.2$(c)$\endsubhead

\proclaim{A.2.8.~Proposition} If $\kappa\notin\QQ_\dom$ the
functor $\dot\otimes$ takes
$\hat\Oc^\Delta_{\dom,\kappa}\times\hat\Oc^\Delta_{\nu,\kappa}$
and
$\hat\Oc^\Delta_{\nu,\kappa}\times\hat\Oc^\Delta_{\dom,\kappa}$
into $\hat\Oc^\Delta_{\nu,\kappa}$.

\endproclaim

\noindent{\sl Proof :} Fix a module $M_1$ in
$\hat\Oc^\Delta_{\dom,\kappa}$ and a module $M_2$ in
$\hat\Oc_{\nu,\kappa}^\Delta$. The module $M_1\dot\otimes M_2$
belongs to the category $\hat\Oc^{fg}_{\nu,\kappa}$ by
Proposition A.2.4. Fix a finite set $B$ such that $M_2$ and
$M_1\dot\otimes M_2$ belong to the subcategory
$\lub\hat\Oc_{\nu,\kappa}^\Delta$. To unburden notation we'll
write $\Bc=\lub\hat\Oc_{\nu,\kappa}^{fg}$. Recall that $\Bc$ is
a highest weight category with a weak duality functor $\dg D$. To
prove the claim it is enough to check that we have
$$\Ext^1_{\Bc}(M_1\dot\otimes M_2,\dg DM)=0,\quad\forall M\in\Delta_{\Bc}.$$
See \cite{D, prop.~A.2.2$(iii)$}. Fix a module $P\in\Bc^\proj$ which
maps onto $M$. Since $P$, $M$ both lie in $\Bc^\Delta$, the kernel
$K$ of the surjective map $P\to M$ lies also in $\Bc^{\Delta}$
\cite{D, prop.~A.2.2$(v)$}. 
Since $\dg DP$ is injective we have also
$$\Ext^1_\Bc(M_1\dot\otimes M_2,\dg DP)=0.$$
So the long exact
sequence of the Ext-group and the Proposition A.2.6$(a)$ yield
$$\dim\Ext^1_\Bc(M_1\dot\otimes M_2,\dg DM)=
\dim\la M_1,M_2,\dg P\ra -\dim\la M_1,M_2,\dg K\ra -\dim\la
M_1,M_2,\dg M\ra.$$ The right hand side is zero by Lemma A.2.10 below.

\qed

\vskip 3mm

\proclaim{A.2.9.~Definition} 
If $M\in\hat\Oc^{fg}_{\nu,\kappa}$ we write $(M:\l)$ for the
coefficient of $[M]$ along the element $[M(\hat\l)_\nu]$ of the
basis of the free Abelian group $[\hat\Oc^{fg}_{\nu,\kappa}]$
consisting of the standard modules.
If $M\in\Oc^{fg}_{\nu}$ we write $(M:\l)$ for
the coefficient of the element $[M]$ along the element $[M(\l)_\nu]$
of the basis of the free Abelian group
$[\Oc^{fg}_{\nu}]$ consisting of the standard modules.
\endproclaim

\proclaim{A.2.10.~Lemma} Let $\kappa\notin\QQ_\dom$. For
$M_1,M_2\in\hat\Oc^\Delta_{\nu,\kappa}$ and
$M\in\hat\Oc^\Delta_{\dom,\kappa}$ we have
$$\dim\la M_1,M,\dg M_2\ra=\sum_{\l_1,\l_2,\mu} (M_1:\l_1)
\,(M_2:\l_2)\,(M:\mu)\,
(M(\l_1)_\nu\otimes L(\mu):\l_2),$$
where $\mu$ runs over $\ZZ^m_\dom$ and $\l_1,\l_2$ run over $\ZZ^\nu_\dom$.
\endproclaim

\noindent{\sl Proof :} First assume that $M=M(\hat\mu)$,
$M_1=M(\hat\l_1)_\nu$ and $M_2=M(\hat\l_2)_\nu$. Note that 
$\dg M_2$ is again a generalized Weyl module. 
Thus Proposition A.2.6$(b)$ yields
$$\aligned
\la M_1,M,\dg M_2\ra &=H_0(\gen, M(\l_1)_\nu\otimes
L(\mu)\otimes\dg M(\l_2)_\nu)\vspace{2mm}
&=\Hom_{\gen}(M(\l_1)_\nu\otimes L(\mu),\dg D M(\l_2)_{\nu})^*.
\endaligned$$
Since the module
$M(\l_1)_\nu\otimes L(\mu)$
lies in $\Oc_{\nu}^\Delta$, the Hom space above has dimension
$$(M(\l_1)_\nu\otimes L(\mu): \l_2)$$
by \cite{D, prop.~A.2.2$(ii)$}. 
The proof is the same if $M$, $M_1$, $M_2$ are
generalized Weyl modules. The general case follows
as in \cite{KL, lem.~28.1}.

\qed

\vskip3mm

\subhead A.2.11.~Proof of Corollary 7.3\endsubhead
The duality functor $D$ equip
$(\Oc^{fg}_{\dom,\kappa},\dot\otimes,a)$ with the structure of
a rigid monoidal category, see \cite{KL}.
This means that there are natural morphisms in $\Oc^{fg}_{\dom,\kappa}$
$$i_M:M(c\o_0)\to M\dot\otimes DM,\quad
e_M: DM\dot\otimes M\to M(c\o_0)$$
such that the following hold

\vskip1mm

\itemitem{$(a)$}
for any $M_1,M_2\in\Oc^{fg}_{\dom,\kappa}$ the map
$\Hom_{\hat\gb}(M_1,DM_2)\to\Hom_{\hat\gb}(M_1\dot\otimes M_2,M(c\o_0))$
such that $f\mapsto e_{M_2}\circ (f\otimes 1)$ is an isomorphism,
\vskip1mm

\itemitem{$(b)$}
the compositions below are equal to the identity
$$\gathered
\xymatrix
{M=M(c\o_0)\dot\otimes M \ar[r]^{i_M\dot\otimes 1}& M\dot\otimes
DM\dot\otimes M\ar[r]^{1\dot\otimes e_M}&M\dot\otimes M(c\o_0)=M},\cr
\xymatrix{DM=DM\dot\otimes M(c\o_0)\ar[r]^{1\dot\otimes i_M}&
DM\dot\otimes M\dot\otimes DM\ar[r]^{e_M\dot\otimes 1}& M(c\o_0)\dot\otimes
DM=DM}.
\endgathered$$
\vskip1mm

\noindent Therefore, for any $M_1,M_2\in\hat\Oc_{\nu,\kappa}^{fg}$ and any
$M\in\hat\Oc_{\dom,\kappa}^{fg}$ the composed map
$$\aligned
\Hom_{\hat\gb}(M_1,M_2\dot\otimes M)\lra
\Hom_{\hat\gb}(M_1\dot\otimes DM,M_2\dot\otimes M\dot\otimes DM)
\lra\cr \lra\Hom_{\hat\gb}(
M_1\dot\otimes DM,M_2\dot\otimes M(c\o_0))=
\Hom_{\hat\gb}(M_1\dot\otimes DM,M_2)\endaligned $$ is
an isomorphism whose inverse map is the composition of the chain of maps
$$\aligned
\Hom_{\hat\gb}(M_1\dot\otimes DM,M_2)\lra
\Hom_{\hat\gb}(M_1\dot\otimes DM\dot\otimes M,M_2\dot\otimes M)
\lra\cr \lra\Hom_{\hat\gb}(M_1\dot\otimes M(c\o_0),M_2\dot\otimes M)=
\Hom_{\hat\gb}(M_1,M_2\dot\otimes M).\endaligned
$$
Thus the functors $-\dot\otimes M$,
$-\dot\otimes DM$ are adjoint (left and right) to each other.
Since every right (resp.~left) adjoint is right (resp.~left) exact
this implies that the functor $-\dot\otimes M$  is exact for each module
$M$ in $\hat\Oc_{\dom,\kappa}^{fg}$. The same holds for the functor
$M\dot\otimes-$.

\subhead A.3.~Reminder on induction\endsubhead
Let $R$ be a commutative ring with 1. A $R$-split induction datum
is a quadruple $(\aen,\ben,\cen,F)$ where $\aen$ is a $R$-Lie
algebra which is free as a $R$-module, $\ben$, $\cen$ are
supplementary $R$-Lie subalgebras of $\aen$ and $F$ is a
$\ben$-module. The corresponding induced module is
$$\Ind_\ben^\aen(F)=U(\aen)\otimes_{U(\ben)}F.$$
For each $\aen$-module $E$ the assignment $a\otimes e\otimes
f\mapsto\sum a_1e\otimes a_2\otimes f$, where $\sum a_1\otimes a_2$
is the coproduct of $a$, yields a $\aen$-module isomorphism
$$\Ind^\aen_\ben(E\otimes_R F)\to E\otimes_R\Ind^\aen_\ben(F).$$ This
isomorphism is called the {\it tensor identity}. There is also a
$\cen$-module isomorphism \cite{KL, sec.~II.A}
$$\Ind^\cen(F)=U(\cen)\otimes_R F\to\Ind^\aen_\ben(F).$$

\subhead A.4. Reminder on the Fock space\endsubhead
This section is a reminder from \cite{U}. Fix an integer $e >1$ and
a tuple $s\in\ZZ^\ell$. First, we
recall the definition of the Fock space $Fock_{s^\circ,e}$, a
$\widehat{\sen\len}_e$-module of level $\ell$ equipped with three
remarkable bases. As a vector space, it is the limit of a filtering
inductive system of vector spaces $\Lambda^k$ with $k>0$. Each
$\Lambda^k$ is equipped with a level zero representation of
$\widehat{\sen\len}_e$. We give the definition of $\Lambda^k$ in the
second section. Next, we introduce a particular submodule
$\Lambda^\nu\subset\Lambda^m$ for each integer $m>0$ and each
composition $\nu\in\Cc_{m,\ell}$. We do not give the construction of
the inductive system above, we'll not use it. However, in the
particular case where $\nu=s$ we define explicitly
the linear map $$\Lambda^s\to Fock_{s^\circ,e}.$$ This
section does not contain new results. Most proofs can be found in
loc.~cit.~and will be omitted.

\subhead A.4.1.~The fock space and its canonical bases\endsubhead 
Let $e_a$, $f_a$, $a\in\ZZ/e\ZZ$, be the
Chevalley generators of the affine Lie algebra
$\widehat{\sen\len}_e$. For each $s\in\ZZ^\ell$ the vector space
$$Fock_{s^\circ,e}=\bigoplus_{\l\in\Pc^\ell}\CC|\l,s^\circ,e\ra$$  is
equipped with a level $\ell$ representation of
$\widehat{\sen\len}_e$ and with two canonical bases
$$(\Gc(\l,s^\circ,e)^\pm;\l\in\Pc^\ell).$$ We define matrices
$\nabla^\pm=(\nabla^\pm_{\l,\mu,s^\circ,e})$ and
$\Delta^\pm=(\Delta^\pm_{\l,\mu,s^\circ,e})$ such that
$$\Gc(\l,s^\circ,e)^\pm=
\sum_\mu\Delta^\pm_{\l,\mu,s^\circ,e}|\mu,s^\circ,e\ra,\quad
\Delta^\pm=(\nabla^\pm)^{-1}.$$ We have
the following formula \cite{U, thm 5.15}
$$\Delta^+_{{}^t\!\mu,{}^t\!\l,-s,e}=\nabla^-_{\l,\mu,s^\circ,e}.\leqno(A.6)$$

%

\subhead A.4.2.~
The $\widehat{\sen\len}_e$-modules $\Lambda^m$, $\Lambda^\nu$
\endsubhead We define a representation of
$\widehat{\sen\len}_e$ on the vector space $U=\bigoplus_{a\in\ZZ}\CC
u_a$ by the following formulas : for $b\in\ZZ/e\ZZ$ and $a\in\ZZ$ we
set
$$\aligned
&e_b(u_{a+1})=u_{a},\ f_b(u_a)=u_{a+1}\ \roman{if}\ b\neq 0,\,a\in
b,\cr &e_0(u_{a+1})=u_{a+e-e\ell},\ f_0(u_a)=u_{a+1-e+e\ell} \
\roman{if}\ a\in e\ZZ,\cr &e_b(u_a)=f_b(u_a)=0\ \roman{else.}
\endaligned$$
It yields a representation of $\widehat{\sen\len}_e$ on the $m$-th
exterior power $\Lambda^m=\wedge^mU$ for each $m>0$. Write
$$|\underline{\a}\ra=u_{a_1}\wedge u_{a_2}\wedge\cdots u_{a_m},\quad\forall
\underline{\a}=(a_1,a_2,\dots a_m)\in\ZZ^m_\ddom.$$ The vectors
$|\underline{\a}\ra$ with $\a\in\ZZ^m_\ddom$ form a basis of
$\Lambda^m$. Let $(\Gc(\underline{\a})^+)$,
$(\Gc(\underline{\a})^-)$ denote the {\sl canonical bases} of
$\Lambda^m$ introduced in loc.~cit.

Now, fix a composition
$\nu\in\Cc_{m,\ell}$. We'll define a
$\widehat{\sen\len}_e$-submodule $\Lambda^\nu\subset\Lambda^m$.
First, for each integer $a$ let $$p_a\in\Lambda,\quad
c_a\in\ZZ/e\ZZ,\quad r_a,\phi_a\in\ZZ$$ be defined by
$$a=c_a+e(p_a-1)+e\ell r_a,\quad\phi_a=c_a+e r_a.$$ This yields the
bijection
$$\ZZ^m_\ddom\to\bigsqcup_{\nu\in\Cc_{m,\ell}}
\ZZ^{\nu^\circ}_\ddom\times\{\nu\}, \quad
\underline{\a}\mapsto(\a,\nu),\leqno(A.7)$$ where
$\underline{\a}=(a_1,a_2,\dots a_m),$
$\a=(\phi_{a_{w(1)}},\dots\phi_{a_{w(m)}})$ and $w\in\Sen$ is
minimal with $$(p_{a_{w(1)}}, \dots
p_{a_{w(m)}})=(\ell^{\nu_{\ell}},\dots 2^{\nu_{2}}, 1^{\nu_{1}}).$$

For example take $e=2$, $\ell=3$, $m=7$,
$\underline\a=(3,1,0,-2,-4,-6,-7)$. Then we get $\nu=(2,2,3)$ and
$\a= (0,-2,-3,1,0,1,0)$.

Let $\Lambda^\nu\subset\Lambda^m$ be the subspace spanned by the
basis elements $|\underline{\a}\ra$ such that $\underline{\a}$ maps
to $\ZZ^{\nu}_\ddom\times\{\nu^\circ\}$ under $(A.7)$. One can
prove that it is a $\widehat{\sen\len}_e$-submodule. 
Let us described explicitly this module.
Fix a tuple $s=(s_p)\in\ZZ^\ell$. We'll
define three bases of $\Lambda^\nu$ whose elements are labeled by
$\ZZ^\nu_\dom$. The $\widehat{\sen\len}_e$-module
$\Lambda^\nu$ does not depend on $s$. The tuple $s$ enters only in the 
labelling of the bases elements.
Let $J_p=J_{\nu,p}$ be as in Section 1.4.
Given a $m$-tuple $\l\in\ZZ^\nu_\dom$ let $$\a(\l,\nu,s)\in\ZZ^\nu_\ddom$$
be the $m$-tuple whose $j$-th entry is $\l_j+i_p-j+s_p$ for $j\in J_p$, and
let $$\underline\a(\l,\nu,s)\in\ZZ^m_\ddom$$ be the unique $m$-tuple such that
$(A.7)$ maps $\underline\a(\l,\nu,s)$ to
$(\a(\l,\nu,s),\nu^\circ)$.

For example take $e=2$, $\ell=3$, $m=7$, $\nu=(2,2,3)$,
$\l=(2,0,1,-3,1,-2,-4)$ and $s=(1,1,4)$. Then we get
$\a(\l,\nu,s)=(3,0,2,-3,5,1,-2)$ and $\underline\a(\l,\nu,s)=
(13,11,4,1,0,-9,-10).$

We define the following elements of $\Lambda^m$
$$|\l,\nu,s^\circ,e\ra=|\underline\a(\l,\nu,s)\ra,\quad
\Gc(\l,\nu,s^\circ,e)^\pm=\Gc(\underline\a(\l,\nu,s))^\pm.$$ 
For each $\l, \mu\in\ZZ^m$ and each $b\in\ZZ$
we write $\l{\buildrel b\over\to}\mu$ if there exist
$j\in J$ such that $\l_j=b$, $\mu_j=b+1$, and $\l_i=\mu_i$
if $i\neq j$.

\proclaim{A.4.3.~Proposition} (a) The elements $|\l,\nu,s^\circ,e\ra$,
$\Gc(\l,\nu,s^\circ,e)^+$ and $\Gc(\l,\nu,s^\circ,e)^-$ are basis
vectors of $\Lambda^\nu$ when $\l$ runs over $\ZZ^\nu_\dom$.
The representation of $\widehat{\sen\len}_e$ in
$\Lambda^\nu$ is given by
$$e_a(|\mu,\nu,s^\circ,e\ra)=\sum_{b,\l}|\l,\nu,s^\circ,e\ra,\quad
f_a(|\l,\nu,s^\circ,e\ra)=\sum_{b,\mu}|\mu,\nu,s^\circ,e\ra,\quad
a\in\ZZ/e\ZZ,$$ summing over all integer $b$ and all 
$\l$, $\mu$ such that $b\in a$ and
$\a(\l,\nu,s)\,\tob\,\a(\mu,\nu,s).$

(b) For each tuple $\mu\in\ZZ^\nu_\dom$ we have 
$$\Gc(\mu,\nu,s^\circ,e)^-=\sum_{\l}
(-1)^{l(v_{\l})-l(v_{\mu})}P^{\g,-1}_{v_{\l},v_{\mu}
}(1)\,|\l,\nu,s^\circ,e\ra.$$ The sum is over all tuples $\l$ such
that $v_{\l}\bullet\g=\widehat {\a(\l,\nu,s)}-\rhoaf$ and
$v_{\mu}\bullet \g=\widehat{\a(\mu,\nu,s)}-\rhoaf$ with
$\g\in\alc$, $v_{\l}, v_{\mu}\in\hat\Sen^\g$ and $v_{\mu}\geqslant
v_{\l}$.

(c) If $\nu=s$ there is a unique linear map
$\Lambda^s\to Fock_{s^\circ,e}$ such that
$|\l,s,s^\circ,e\ra\mapsto|\l^\circ,s^\circ,e\ra$,
$|\mu,s,s^\circ,e\ra\mapsto 0$ and
$\Gc(\l,s,s^\circ,e)^\pm\mapsto\Gc(\l^\circ,s^\circ,e)^\pm$ for
$\l\in\NN^s_\dom$ and $\mu\in\ZZ^s_\dom\setminus\NN^s_\dom$.
\endproclaim

\vskip3mm

\subhead A.4.4.~Remarks\endsubhead $(a)$ If $\nu=s$
then we have
$$\a(\l,s,s)-\rho=\l_\pi,$$ where $\pi$ is as in Section 8.1.

\vskip1mm

$(b)$ Proposition A.4.3$(b)$ is \cite{U,thm.~3.25-3.26}. 
The tuple $\underline\a(\l,s,s)$ which is used in Proposition A.4.3$(c)$
differs from the combinatorial definition in loc.~cit. Let us
explain this and let us explain how to deduce the lemma from
loc.~cit.
For any $\ell$-partition $\l$ and any tuple $s\in\ZZ^\ell$ we write
$$\a(\l,s)= \{\c_s(i,\l_{p,i}+1,p);i>0,p\in\Lambda\},\quad
\c_s(i,j,p)=s_p+j-i.$$ Set
$\Ac=\{\a(\l,s);\l\in\Pc^\ell,s\in\ZZ^\ell\}.$ In the particular
case where $\ell=1$ we write $\underline\l$, $\underline\a$,
$\underline\Ac$ for $\l$, $\a$, $\Ac$.
Consider the bijection
$$\underline\Ac\to\Ac,\quad
\underline\a\mapsto\a=\{(\phi_a,p_a); a\in\underline\a\}.$$ Let
$\a\mapsto\underline\a$ denote the inverse map.

Next, fix an integer $m>0$ and a composition $s\in\Cc_{m,\ell}$. For
each $\l\in\Pc^\ell$ there is a unique $\underline\l\in\Pc$ such
that $\underline\a(\underline\l,m)=\underline{\a(\l,s)}$. This
yields a map $\Pc^\ell\to\Pc$, $\l\mapsto\underline\l$.

Now let $\l\in\Pc^\ell$. The set $\underline{\a(\l^\circ,s^\circ)}$
is well-defined, it belongs to $\underline\Ac$. For each integer
$k>0$ let
$\underline\a(\l,k,s^\circ)$ be the tuple consisting of the $k$
largest entries of $\underline{\a(\l^\circ,s^\circ)}$ arranged in
decreasing order. This tuple belongs to $\ZZ^{k}_\ddom$.

Assume further that $\l$ is indeed a tuple in $\NN^ s_\dom$ which is
viewed as a $\ell$-partition as in Section 1.5. Then a direct
computation yields $m\geqslant l(\underline\l)$. Using this
inequality, another computation yields
$$\underline\a(\l,m,s^\circ)=\underline\a(\l,s,s).$$ So Proposition 
A.4.3$(c)$ follows from \cite{U, sec.~4.4}.

For example take
$$e=2,\quad\ell=3,\quad m=6,\quad\l=((1,1),(2,1),(1)),\quad s=(2,3,1).$$
We can identify $\l$ with the tuple $(1,1,2,1,0,1)\in \NN^s_\dom$.
Since $\underline\l=(9,4,3^2,1^2)$ we have $m=6=l(\underline\l)$. We
have also $\a(\l,s,s)=(3,2,5,3,1,2)$. Thus
$\underline\a(\l,s,s)=(15,11,9,6,3,2)$ by definition of $(A.7)$. The
tuple $\underline\a(\l,s,s)$ coincides with the first six entries of
$\underline{\a(\l^\circ,s^\circ)}$ arranged in decreasing order.

\vskip3mm

\subhead A.5.~The functor $\Een$ for $\ell=1$ \endsubhead
In this section we'll set $\ell=1$, $\nu=(m)$ and $n\leqslant m$.
Recall that we have identified the set of partitions of $n$ with a set
of integral dominant weights in (1.1). Write
$$\Delta_{\l,\dom,\kappa}=\Delta_{\l,\nub},\quad 
S_{\l,\dom,\kappa}=S_{\l,\nub}$$
for each dominant weight $\l$.
See Section 6.1 for details.
Let $$\Ac_{n,\dom,\kappa}\subset\hat\Oc_{\dom,\kappa}$$ 
be the Serre subcategory
generated by the modules $S_{\l,\dom,\kappa}$ with $\l\in\Pc_n$. It
is a quasihereditary category with respect to the order
$\trianglelefteq$, see e.g., Proposition 8.7. Since
$\ell=1$ the algebra $\Hb_{h,H}$ is the
rational DAHA of $\GL_n(\CC)$ with the parameter $h$ and $\Ab_{q,Q}$
is the Hecke algebra of $\GL_n(\CC)$ with the parameter $q$. 
Set $q=\exp(2i\pi h)$. Note that 
$$\Een:\hat\Oc_{\dom,\kappa}\to\Hc_{h,H}$$ is Suzuki's functor. The aim of
this section is prove the following theorem.

\proclaim{A.5.1. Theorem} Assume that $\kappa\notin\QQ_\dom$,
$h=1/\kappa$ and $n\leqslant m$. Assume also that
$(q+1)(q^2+q+1)\neq 0$. The functor $\Een$ restricts to an
equivalence of quasi-hereditary categories
$\Ac_{n,\dom,\kappa}\to\Hc_{h,H}$ which takes 
$\Delta_{\l,\dom,\kappa}$ to $\Delta_{\l,h,H}$ for each $\l\in\Pc_{n}$.
\endproclaim

To prove Theorem A.5.1 we need some material.
Let $\Ab$ be a finite dimensional $\CC$-algebra with 1. Let
$(\Ac,\Delta_\Ac,F)$, $(\Bc,\Delta_\Bc,G)$ be $\CC$-linear 1-faithful
covers of $\Ab$. See \cite{R, def.~4.37} for the terminology. Hence
$$F:\Ac\to\Ab\text{-}\modb^\f,\quad
G:\Bc\to\Ab\text{-}\modb^\f$$ are functors which restrict to equivalences
of exact categories
$$\Ac^\Delta\to(\Ab\text{-}\modb)^{F(\Delta_\Ac)},\quad
\Bc^{\Delta}\to(\Ab\text{-}\modb)^{G(\Delta_\Bc)}$$  by \cite{R,
prop.~4.41(2)}. The following lemma follows easily from loc.~cit.

\proclaim{A.5.2.~ Lemma} Assume there is a functor $\phi:\Ac\to\Bc$
such that $F=G\circ\phi$ and $\phi(\Delta_\Ac)=\Delta_\Bc$. Then
$\phi$ yields an equivalence of exact categories
$\Ac^\Delta\to\Bc^\Delta$. In particular $\phi$ is
fully faithful on $\Ac^\proj$ and it takes a projective generator
of $\Ac$ to a projective generator of $\Bc$.
\endproclaim

\noindent
We'll also use the following lemma.

\proclaim{A.5.3.~Lemma} Let $\phi:\Ac\to\Bc$ be a functor of
Artinian Abelian categories and $P$ be a projective generator of
$\Ac$. Set $Q=\phi(P)$. Assume that

(a) the module $Q$ is a projective generator of $\Bc$,

(b) the functor $\phi$ is fully faithful on $\Ac^\proj$,

(c) the category $\Ac$ has finite projective dimension,

(d) the functor $\phi$ is right exact.

\noindent Then $\phi$ is an equivalence of categories.
\endproclaim

\noindent{\sl Proof :} Note that $(a)$, $(b)$ imply that the
categories $\Ac$, $\Bc$ are equivalent. We must prove that the functor $\phi$
is an equivalence.
Set
$$\Ab=(\End_\Ac P)^\op,\quad\Bb=(\End_{\Bc} Q)^\op.$$
The functor
$$\Hom_{\Bc}(Q,-):\Bc\to\Bb\text{-}\modb^\f$$
is an equivalence by $(a)$. Thus it is enough
to check that the functor
$$G=\Hom_{\Bc}(Q,\phi(-)):\Ac\to\Bb\text{-}\modb^\f$$
is an equivalence.
Note that $\phi$ yields a ring
homomorphism $\Ab\to \Bb$.
So we have the functor
$$F=\Bb\otimes_\Ab\Hom_\Ac(P,-):\Ac\to\Bb\text{-}\modb^\f,$$
and the morphism of functors
$$\Phi:
F\to G,\ r\otimes f\mapsto r\phi(f).$$
The ring homomorphism $\phi:\Ab\to\Bb$ is invertible by part $(b)$.
The functor
$$\Hom_\Ac(P,-):\Ac\to\Ab\text{-}\modb^\f$$
is an equivalence. Thus $F$ is also an equivalence.
Therefore it is enough to
prove that $\Phi$ is a isomorphism of functors.

First, assume that $M=P$. Then $F(M)=G(M)=\Bb$
and $\Phi(M)$ is the identity of $\Bb$.
Next, let $M$ be a projective object of $\Ac$. Then the morphism
in $\Bb$-$\modb^\f$
$$\Phi(M):F(M)\to G(M)$$ is invertible.
Indeed, we may assume that $M$ is
indecomposable. Then $M$ is a direct summand of $P$. 
So $\Phi(M)$ is the identity of the $\Bb$-module $\Bb\phi(a)$
for some idempotent $a\in \Ab$. 
Finally, let $M$ be any object of $\Ac$.
By $(c)$ the projective dimension of $M$ is
$e<\infty$. Fix an exact sequence
$$0\to M_2\to M_1\to M\to 0$$ with
$M_1\in\Ac^\proj$ and $M_2\in\Ac$ of projective dimension $<e$. Consider
the diagram
$$\matrix
F(M_2)&\to&F(M_1)&\to&F(M)&\to&0\cr
\downarrow&&\downarrow&&\downarrow&&\cr
G(M_2)&\to&G(M_1)&\to&G(M)&\to&0.
\endmatrix$$
Here the vertical maps are $\Phi(M_2)$, $\Phi(M_1)$ and $\Phi(M)$. The
functor $G$ is right exact by $(d)$. Thus both rows are exact. We
may assume that $\Phi(M_2)$, $\Phi(M_1)$ are both invertible by
induction on $e$. Thus $\Phi(M)$ is also invertible by the five
lemma. We are done.

\qed

\vskip3mm

\subhead A.5.4.~ Proof of Theorem A.5.1\endsubhead
First, observe that
$\VVkap(\Delta_{\l,\dom,\kappa})=\Delta_{\l,h,H}$ for each $\l\in\Pc_{n}$
by Proposition 6.2$(a)$. Thus it is enough to check that $\Een$ is
an equivalence of categories
$\Ac_{n,\dom,\kappa}^\f\to\Hc_{h,H}^\f.$
Set 
$$\Ac=\Ac_{n,\dom,\kappa}^\f,\quad\Bc=\Hc_{h,H}^\f,\quad\Ab=\Ab_{q,Q},\quad
\phi=\Een.$$ The hypothesis $(c)$, $(d)$ in Lemma A.5.3 are
obviously true. By \cite{R, thm.~5.3} the functor
$$G=\KZ\circ\heartsuit:\Bc\to\Ab\text{-}\modb^\f$$
is a 1-faithful
cover. We claim that the functor
$$F=G\circ\phi:\Ac\to\Ab\text{-}\modb^\f$$
is also a 1-faithful cover. Therefore $\phi$ satisfies also the
hypothesis $(a)$, $(b)$ in Lemma A.5.3, by Lemma A.5.2. Hence $\phi$
is an equivalence of categories.
Now we prove the claim. Write 
$$\Vb_{n,\dom,\kappa}=( \Vb_\kappa)^{\dot\otimes n},\quad
\Ab_{n,\dom,\kappa}=\End_{\hat\gb}(\Vb_{n,\dom,\kappa}).$$
Proposition 7.6 gives an algebra homomorphism
$$\Ab\to\Ab_{n,\kappa,\dom}\leqno(A.8)$$
such that
$$F=\Hom_{\hat\gb}(\Vb_{n,\kappa,\dom},-),$$
up to a twist by some duality functor that we omit to
simplify. Let $U_q(\gen)$ be the quantized enveloping algebra of
$\gen$ with the parameter $q$ and let $\Vb_q$ be its vectorial
representation. Under the Kazhdan-Lusztig tensor equivalence
\cite{KL, thm.~IV.38.1} the category $\hat\Oc^\f_{\dom,\kappa}$ is
equivalent to the category of finite dimensional $U_q(\gen)$-modules. 
Therefore the
ring homomorphism $(A.8)$ is invertible, because it is taken to the
isomorphism
$$\Ab\to\End_{U_q(\gen)}(\Vb_q^{\otimes n})$$
given by the Schur-Weyl duality. Thus the functor $F$ is taken to
the Schur functor $$M\mapsto\Hom_{U_q(\gen)}(\Vb_q^{\otimes
n},M).$$ It is well-known that the Schur functor is a 1-faithful
cover, see \cite{R, rem.~6.7} and the reference there for instance.
Hence $F$ is also a 1-faithful cover.

\qed

\vskip3mm

\subhead A.5.5.~Remark\endsubhead The idea to use the Kazhdan-Lusztig
equivalence to prove Theorem A.5.1 is not new. 
However we have not found any proof of Theorem
A.5.1 in the literature.

\subhead A.6.~Proof of Proposition 8.7\endsubhead
By \cite{CPS, thm.~3.5$(a)$} it is enough to
prove the following.

\proclaim{A.6.1.~Proposition} If $\mu+\pi\trianglelefteq\l+\pi$,
$\l\in\Pc_{n,s}$ and $\mu\in\ZZ^s_\dom$, then
$\mu\in\Pc_{n,s}.$
\endproclaim

\noindent{\sl Proof :} Recall that $\l\triangleleft\mu$ iff
there are $\mu_1,\mu_2,\dots,\mu_r\in\ZZ^s_\dom$ such that
$$\tilde\l=\tilde\mu_1<\tilde\mu_2<
\tilde\mu_3<\dots<\tilde\mu_r=\tilde\mu$$ and such that
$\tilde\mu_{i+1}=w_is_{\a_i}\bullet\tilde\mu_i$ for some 
$\a_i\in\hat\Pi_\re\setminus\Pi_s$ and some
$w_i\in\Sen_s$. Now, assume that $\mu+\pi\trianglelefteq\l+\pi$,
$\l\in\Pc_{n,s}$, and $\mu\in\ZZ^s_\dom$.
By an easy induction we may assume that
$$\tilde\mu+\pi=w
s_{\a}\bullet(\tilde\l+\pi),\quad
\la\tilde\l+\rhoaf+\pi:\a\ra\in\ZZ_{>0},
\quad\a\in\hat\Pi^+_\re\setminus\Pi_s,
\quad w\in\Sen_s.$$ So we have $|\mu|=n$, and we must prove that
$\mu_1, \mu_2,\dots, \mu_m$ are $\geqslant 0$. There is an unique
map
$$\CC^m\to\CC^\ZZ,\quad \l\mapsto\bar\l$$ such that
$\bar\l_1,\bar\l_2,\dots\bar\l_{m}$ are the entries of
$\l+\pi+\rho$ and $\bar\l_{j+m}=\bar\l_j-\kappa$ for all
$j\in\ZZ$. Under this map the dot action of the affine reflection
$s_{\a}$ is taken to the linear operator which switches the
$(a+km)$-th and the $(b+km)$-th entries of any sequence for each
$k\in\ZZ$ and some fixed integers $a\neq b$. We have
$\bar\l_j\geqslant\bar 0_j$ for all $j\in\ZZ$. We must check that
the same holds for the entries of $\bar\mu$.

Recall the partition $J=\bigsqcup_{p\in\Lambda} J_{s,p}$ 
with $J_{s,p}=[i_p,j_p]$. Since $\mu\in\ZZ^s_\dom$ it is enough to
prove that we have
$$\bar\mu_{j_1}\geqslant\bar 0_{j_1},\quad
\bar\mu_{j_2}\geqslant\bar
0_{j_2},\quad\dots\quad\bar\mu_{j_\ell}\geqslant\bar 0_{j_\ell}.$$
The $\ell$-tuples $(i_p)$, $(j_p)$ can be regarded as sequences of
integers such that
$$i_{p+\ell}=i_p+m,\quad j_{p+\ell}=j_p+m,\quad\forall p\in\ZZ.$$
For all $p\in\ZZ$ we set also $J_{s,p+\ell}=J_{s,p}+m$. Now fix
$p$, $q$ such that
$$a\in J_{s,p},\quad b\in J_{s,q}.$$  It is enough to prove that
$\bar\mu_{j_p}\geqslant\bar 0_{j_p}$ and that $
\bar\mu_{j_q}\geqslant\bar 0_{j_q}$. Assume that $b>a$. Then $q>p$
because $\a\notin\Pi_s$. Since $\bar\l_a>\bar\l_b$ we have $\bar
0_{j_p}-\bar 0_{j_q}\in\ZZ$. Since $\kappa\notin\RR_{\geqslant 0}$
we have $\bar 0_{j_q}\geqslant\bar 0_{j_p}$. Note that
$$\aligned
\{\bar\mu_i;i\in J_{s,p}\}= \{\bar\l_i;i\in
J_{s,p}\}\setminus\{\bar\l_a\}\cup\{\bar\l_b\}, \cr
\{\bar\mu_i;i\in J_{s,q}\}= \{\bar\l_i;i\in
J_{s,q}\}\setminus\{\bar\l_b\}\cup\{\bar\l_a\}.
\endaligned$$
Therefore we have
$$\aligned
\bar\mu_{j_p}&=\inf\{\bar\mu_i;i\in J_{s,p}\}
\geqslant\inf(\{\bar\l_i;i\in
J_{s,p}\}\cup\{\bar\l_b\})\geqslant\inf\{\bar 0_{j_p},\bar
0_{j_q}\}\geqslant \bar 0_{j_p},\cr
\bar\mu_{j_q}&=\inf\{\bar\mu_i;i\in J_{s,q}\}
\geqslant\inf(\{\bar\l_i;i\in J_{s,q}\}\cup\{\bar\l_a\})
\geqslant\inf\{\bar\l_i;i\in J_{s,q}\} \geqslant\bar 0_{j_q}.
\endaligned$$

\qed

\vskip3cm

\head Index of notation\endhead

\itemitem{0.1 :}
$[\Ac]$,
$[M]$,
$\Ac^\Delta$,
$\Delta_\Ac$,
$\Ac^\proj$,
$\Ac^\f$,

\vskip1mm
\itemitem{0.2 :}
${}^\phi M$,

\vskip1mm
\itemitem{0.3 :}
$M[X]$,
$M_\Rb$,
$M^F$,
$M^F_\Rb$,

\vskip1mm
\itemitem{1.1 :} 
$D_\ell$,
$\Sen_n$,
$W$,
$A$,
$W_A$,
$\Lambda$,
$\epsilon$,
$\epsilon_i$,
$s_{i,j}$,
$s_{i,j}^{(p)}$,

\vskip1mm
\itemitem{1.2 :} 
$x_i$,
$y_i$,
$\Hb_{k,\g}$,
$\Hb_{h,H}$,

\vskip1mm
\itemitem{1.3 :} 
$\Rb$, 
$\Rb^*$,
$\bar y_i$,

\vskip1mm
\itemitem{1.4 :} 
$\Cc_{m,\ell}$,
$\Cc_{m,\ell,n}$,
$J$,
$J_p=J_{\nu,p}$,
$i_p$,
$j_p$,
$\CC^\nu_\dom$,
$\ZZ^\nu_\dom$,
$\ZZ^\nu_\ddom$,
$\Pc_n$,
$n(\l)$,
$|\l|$,
$\Pc^\ell_n$,
$\Pc^\ell$,
$\Pc_{n,\nu}^\ell$,
${}^t\l$,
$\NN^\nu_\dom$,
$\nu^\circ$,
$\nu^\bullet$,

\vskip1mm
\itemitem{1.5 :}
$\chi_p$,
$\Irr(\CC \Sen_n)$,
$\Irr(\CC W)$,
$\Xen_\l$,
$A_{\mu,p}$,
$\Sen_\mu$,
$W_\mu$,
$w_\mu$,
$\Gamma$,

\vskip1mm
\itemitem{1.6 :}
$\Hc_{h,H}$,
$\Hc_{h,H}^\f$,
$\Delta_{\l,h,H}$,
$S_{\l,h,H}$,
$P_{\l,h,H}$,
$\eu$,
$\eu_0$,
$\theta_\l$,
$\succcurlyeq$,

\vskip1mm
\itemitem{1.7 :}
$\Rb_\nl$,
$C_\nl$,
$M_\nl$,
$\Bb$,
$\Bb_\nl$,
$\Hb_{h,H,\nl}$,

\vskip1mm
\itemitem{2.1 :}
$\gen$,
$G$,
$g\xi$,
$\ben$,
$\ten$,
$T$,
$\epsilon_i$,
$\check\epsilon_i$,
$\l$,
$\l_i$,
$\check\l$,
$\check\l_i$,
$\rho$,
$\a_i$,
$I$,
$\Pi$,
$\Pi^+$,
$e_{k,l}$,
$e_i$,
$f_i$,
$L(\l)$,
$\Xen_\l$,
$\Vb$,
$\Vb^*_p$,

\vskip1mm
\itemitem{2.2 :}
$\gb$,
$\gb_\dom$,
$\bb$,
$\hat\gb$,
$\hat\bb$,
$\tilde\gb$,
$\tilde\bb$,
$\hat\gb_\dom$,
$\un$,
$\partial$,
$\tb$,
$\gen_R$,
$\hat\gb_R$,
$\hat\Pi$,
$\hat\Pi^+$,
$\hat\Pi_\re$,
$\delta$,
$\o_0$,
$\hat I$,
$\hat\a_i$,
$\la\ :\ \ra$,

\vskip1mm
\itemitem{2.3 :}
$c=\kappa-m$,
$\hat\gb_{R,\kappa}$,

\vskip1mm
\itemitem{2.4 :}
$\Cc_{R,\kappa}$,
$Q_{R,\kappa}$,
$M(r)$,
$M(\infty)$,

\vskip1mm
\itemitem{2.5 :}
$\xi^{(r)}$,
$\Lb_s$,
$\Omega$,

\vskip1mm
\itemitem{2.6 :}
$\sp M$,
$\dg M$, 
$M^*$, 
$M^d$,
$D$,
$\dg D$,

\vskip1mm
\itemitem{2.7 :}
$\hat\qb$,
$\hat\lb$,
$\hat\Oc_\kappa$,
$\qen_\nu$,
$\hen_\nu$,
$\hat\qb_\nu$,
$\ub_\nu$,
$\hat\Oc_{\nu,\kappa}$,
$\hat\Oc_{\dom,\kappa}$,
$\Oc$,
$\Oc_\nu$,
$\Oc_\dom$,
$\qen'_\nu$,
$\hat\qb'_\nu$,
$\hat\bb'$,

\vskip1mm
\itemitem{2.8 :}
$M_\nu$,
$M_{\nu,\kappa}$,
$M_\kappa$,
$L(\hen_\nu,\l)$,
$M(\l)_\nu$,
$M(\hat\l)_\nu$,
$\hat\l=\l+c\o_0$,
$M(\hat\l)$,
$\hat\Oc'_{\nu,\kappa}$,

\vskip1mm
\itemitem{2.10 :}
$\tilde\qb$,
$\tilde\Oc_\kappa$,
$\tilde\Oc_{\nu,\kappa}$,
$M_\l$,
$\tilde\l$,
$z_\l$,

\vskip1mm
\itemitem{2.11 :}
$R((t_S))$,
$f(t)_{[i]}$,
$\Gc_R$,
$\Gc_{R,S}$,
$\hat\Gc_{R,S}$,
$\hat\Gc_{R,S,\kappa}$,
$\g_{(i)}$,
$\s M$,

\vskip1mm
\itemitem{2.12 :}
$z_i$, 
$x_i$,
$\PP^1_x$,
$\iota_x$,
$\la M_i;i\in S\ra_x$,
$\hat S$,
$\hat x$,
$C_n$,
$\Rb_n$,
$\la M_i;i\in S\ra$,

\vskip1mm
\itemitem{2.17 :}
$\Bb_n$,
$T(M)$,

\vskip1mm
\itemitem{2.19 :}
$\nabla_i$,
$\gamma_i$,
$\gamma_{i,j}$,

\vskip1mm
\itemitem{3.1 :}
$g$,
$\gen_p$,
$\hen$,
$F$,
$\hat\gb^F$,
$\hat\gb^F_\kappa$,

\vskip1mm
\itemitem{3.2 :}
$\hat\Oc_\kappa^F$,
$\hat\Oc_{\dom,\kappa}^F$,
$\Cc_\kappa^F$,
$F'$,
$M(\hat\l)^F$,
$\flat$,

\vskip1mm
\itemitem{3.3 :}
$\sigma_{i,j}^{(p)}$,
$\Xen(M)$,
$\Cen(M)$,

\vskip1mm
\itemitem{3.9 :}
$\PP^1_y$,
$z_{i,p}$,
$\Gc_S^F$,
$\hat\Gc_S^F$,
$\la M_i;i\in S\ra$,
$\la M_i;i\in S\ra'$,

\vskip1mm
\itemitem{4.1 :}
$\l_\pi$,
$\hat\l_\pi$,
$\pi$,
$\g$,

\vskip1mm
\itemitem{4.4 :}
$T(M',M)$,
$\Cen(M',M)_\nl$,

\vskip1mm
\itemitem{5.1 :}
$\hat\Sen$,
$w\bullet\l$,
$\hat\Pi(\l)$,
$\hat\Sen(\l)$,
$\hat\Pi_\l$,
$\hat\Sen_\l$,
$\pi_\nu$,
$\Sen_\nu$,
$\tb^*_0$,

\vskip1mm
\itemitem{5.2 :}
$\l_+$,
$sn(\l)$,
$\preccurlyeq$,
$\trianglelefteq$,

\vskip1mm
\itemitem{5.5 :}
$\hat\Pi^\nu$,
$\zb$,
$\zen$,
$\lub\Lambda$,
$\lub\hat\Oc_{\nu,\kappa}$,

\vskip1mm
\itemitem{5.7 :}
$\Aen$,

\vskip1mm
\itemitem{5.8 :}
$P^{\l,-1}_{v,w}$,

\vskip1mm
\itemitem{6.1 :}
$\Een$,
$\tilde\l_\pi$,
$\Delta_{\l,\nu,\kappa}$,
$S_{\l,\nu,\kappa}$,
$\lub P_{\l,\nu,\kappa}$,

\vskip1mm
\itemitem{7.1 :}
$(\hat\Oc^\f_{\dom,\kappa},\dot\otimes,a,M(c\o_0)),$

\vskip1mm
\itemitem{7.4 :}
$\Ab_{q,Q}$,
$\KZ$,
$\heartsuit$,
$\Een_\KZ$,

\vskip1mm
\itemitem{7.6 :}
$\Vb_{n,\nu,\kappa}$,
$\Ab_{n,\nu,\kappa}$,
$\Fen$,

\vskip1mm
\itemitem{8.1 :}
$\Ac_{s,-e}$,

\vskip1mm
\itemitem{8.3 :}
$\Sb_{q,Q}$,
$\Delta_{\l,q,Q}$,
$S_{\l,q,Q}$,

\vskip1mm
\itemitem{8.5 :}
$\Ac_{n,s,-e}$,

\vskip1mm
\itemitem{8.9 :}
$\Ac_{n,\nu,\kappa}$,

\vskip1mm
\itemitem{A.2 :}
$\dot\otimes$,
$\Cc$,

\vskip1mm
\itemitem{A.3 :}
$\Gamma^\aen_\ben$,

\vskip1mm
\itemitem{A.4 :}
$Fock_{s^\circ,e}$,
$\Gc(\l,s^\circ,e)^\pm$,
$|\l,s^\circ,e\ra$,
$\nabla^\pm_{\l,\mu,s^\circ,e}$,
$\Delta^\pm_{\l,\mu,s^\circ,e}$,
$\Lambda^m$,
$\Lambda^\nu$,
$\l{\buildrel b\over\to}\mu$,

\vskip1mm
\itemitem{A.5 :}
$\Delta_{\l,\dom,\kappa}$,
$S_{\l,\dom,\kappa}$,
$\Ac_{n,\dom,\kappa}$.

\enddocument